\documentclass[11pt,twoside]{article}

\usepackage{mathrsfs}
\usepackage{amssymb,amsmath,amsfonts,amsthm,color,mathrsfs}
\usepackage[Symbol]{upgreek}
\usepackage{txfonts}
\usepackage[nottoc,notlot,notlof]{tocbibind}
\usepackage[active]{srcltx}
\usepackage[top=1in, bottom=1in, left=1.25in, right=1.25in]{geometry}

\theoremstyle{plain}

\newtheorem{thm}{Theorem}[section]
\newtheorem{lem}[thm]{Lemma}

\newtheorem{cor}[thm]{Corollary}

\newtheorem{defn}[thm]{Definition}
\newtheorem{rem}[thm]{Remark}

\numberwithin{equation}{section}
\allowdisplaybreaks

\begin{document}
\title{ Regularity of Solutions of the Camassa-Holm Equations with Fractional Laplacian Viscosity}
\author{ Zaihui Gan\thanks{ Corresponding author:
ganzaihui2008cn@tju.edu.cn(Zaihui Gan)}, Yong He\thanks{Yong He: heyong80@tju.edu.cn}, Linghui Meng\thanks{Linghui Meng: mlh2017@tju.edu.cn},Yue Wang\thanks{Yue Wang: wangy2017@tju.edu.cn}
\\  {\small   Center for Applied Mathematics, Tianjin University,  Tianjin 300072, China}
 }
\date{}
\maketitle {
 {
\begin{quote}
{\bf Abstract.}
We study the existence, uniqueness and regularity of solutions to the $n$-dimensional ($n=2,3$) Camassa-Holm equations with fractional Laplacian viscosity with smooth initial data. It is a coupled system between the Navier-Stokes equations with nonlocal viscosity and a Helmholtz equation. The main difficulty lies in establishing some a priori estimates for the fractional Laplacian viscosity. To achieve this, we need to explore suitable fractional-power Sobolev-type estimates, and bilinear estimates for fractional derivatives. Especially, for the critical case $\displaystyle s=\frac{n}{4}$ with $n=2,3$, we will make extra efforts for acquiring the expected estimates as obtained in the case $\displaystyle  \frac{n}{4}<s<1$. By the aid of the fractional Leibniz rule and the nonlocal version of Ladyzhenskaya's inequality, we prove the existence, uniqueness and regularity to the Camassa-Holm equations under study by the energy method and a bootstrap argument, which rely  crucially on the fractional Laplacian viscosity. In particular, under the critical case $s=\dfrac{n}{4}$, the nonlocal version of Ladyzhenskaya's inequality is skillfully used, and the smallness of initial data in several Sobolev spaces is required to gain the desired results concernig existence, uniqueness and regularity.
\\[0.15cm]
\\[0.15cm]
\end{quote}}}
%\begin{abstract}
\allowdisplaybreaks
\section{Introduction}
Recently a great attention has been devoted to the study of nonlocal problems driven by fractional Laplacian type operators in the literature, which is not only for a pure academic interest, but also for the various applications in different fields. The Camassa-Holm equations with fractional dissipation naturally occur in hydrodynamics. Based on these known results for the equations with fractional Laplacian operator in \cite{ Caffarelli-1}, C\'{o}rdoba-C\'{o}rdoba-Fontelos \cite{Cordoba,Cordoba1}, Fujiwara-Georgiev-Ozawa \cite{FGO},  Kenig-Ponce-Vega \cite{KPV}, Musina-Nazarov \cite{Musina},  and Nezza-Palarucci-Valdinoci \cite{Nezza}, this paper is devoted to the study of the existence, uniqueness and regularity properties of weak solutions of the following Camassa-Holm equations with fractional Laplacian viscosity:
 \begin{equation}\label{1.1}
\left\{
\begin{array}{ll}
\mathbf{v}_{t}+\mathbf{u}\cdot \nabla \mathbf{v}+\mathbf{v}\cdot \nabla \mathbf{u}^{T}+\nabla p=-\nu (-\Delta)^{s}\mathbf{v},~(t,x)\in \mathbb{R}^{+}\times \mathbb{R}^{n},
\\[0.2cm]
\mathbf{ u}-\alpha^{2} \Delta \mathbf{u} =\mathbf{v},
\\[0.2cm]
\mbox{div}~ \mathbf{v}=\mbox{div}~ \mathbf{u}=0,
  \end{array}
\right.
\end{equation}
 subject to the prescribed initial condition
\begin{equation}\label{1.2}
 \mathbf{ v}(0,x)=\mathbf{v}_0(x),~~x\in \mathbb{R}^{n}.
\end{equation}
Here, $n=2,3$, $\displaystyle\frac{n}{4}\leq s< 1$, $\mathbf{v},\mathbf{u},p$ denotes three unknowns and the divergence free condition $\mbox{div} ~\mathbf{v}=\mbox{div} ~\mathbf{u}=0$ indicates the incompressibility of the fluid. In particular, $\mathbf{v}$ denotes the fluid velocity field, $\mathbf{u}$ the filtered fluid velocity, and $p$ the scalar pressure.  $\Lambda^{2s}:=(-\Delta)^{s}$ is the fractional Laplacian operator and  $\displaystyle 2^{*}_{s}=\frac{2n}{n-2s}$  the fractional critical Sobolev exponent. The parameter $\alpha$ is defined by solving the Helmholtz equations $\mathbf{u}-\alpha^{2} \Delta \mathbf{u} =\mathbf{v}$. \eqref{1.1} represents as the filter associated with the Camassa-Holm equations with fractional Laplacian viscosity, where $\alpha$ is a length scale parameter representing the width of the filter, and $\nu>0$ is the viscosity coefficient fixed in our discussions.\\
\indent Recall that the nonlocal Camassa-Holm equations \eqref{1.1} with $s=1$ read
 \begin{equation}\label{1.3}
\left\{
\begin{array}{ll}
\mathbf{v}_{t}+\mathbf{u}\cdot \nabla \mathbf{v}+\mathbf{v}\cdot \nabla \mathbf{u}^{T}+\nabla p= \nu  \Delta \mathbf{v},
 \\[0.2cm]
\mathbf{u}-\alpha^{2} \Delta \mathbf{u} =\mathbf{v},
 \\[0.2cm]
\mbox{div}~ \mathbf{v}=\mbox{div}~ \mathbf{u}=0.
 \end{array}
\right.
\end{equation}
 As it is well-known that the system \eqref{1.3} rose from works on shallow water equations \cite{Camassa-1,Holm-Marsden}. Specifically, it was introduced in \cite{Holm-Marsden} as a natural mathematical generalization of the integrable inviscid one-dimensional Camassa-Holm equations discovered in \cite{Camassa-1} through a variational formulation and with a lagrangian averaging. It could be used as a closure model for the mean effects of subgrid excitations, and be also viewed as a filtered Navier-Stokes equations with the parameter $\alpha$ in the filter, which obeys a modified Kelvin circulation theorem along filtered velocities \cite{Holm-Marsden}. Numerical examples that seem to justify this intuition were reported in \cite{Chen-holm}. The classical results on the existence, uniqueness and regularity for the equations \eqref{1.3} were established in \cite{Clayton-08,Jerrold-1}. Indeed, many different problems driven by the fractional Laplacian were considered in order to achieve existence, uniqueness and regularity, and  also to obtain qualitative properties of the solutions.\\
\indent For the nonlocal operator $(-\Delta)^{s}$, known as the fractional Laplacian of order $2s$ in the whole space, there are several ways to define it \cite{Caffarelli-1}. Let $\mathscr{S}(\mathbb{R}^{n}) $ be the Schwartz class. First of all, it is defined for any $g\in \mathscr{S}(\mathbb{R}^{n}) $ through the Fourier transform: if $(-\Delta)^{s}g=h$, then
\begin{equation}\label{1.4}
 \widehat{h}(\xi)=|\xi|^{2s}\widehat{g}(\xi).
\end{equation}
Secondly, if $0<s<1$ and  a function $f\in  \mathscr{S}(\mathbb{R}^{n}) $, using the rpresentation by means of a hypersingular kernel \cite{Landkof}, it can be defined as
\begin{equation}\label{1.5}
\left.
\begin{array}{ll}
 I_{s}f(x)=(-\Delta)^{s}f(x)&:=C_{n,s}~{\mbox P.V.}\displaystyle\int_{\mathbb{R}^{n} }\frac{f(x)-f(\xi)}{|x-\xi|^{n+2s}}d\xi\\
 \\
 &=C_{n,s}\lim\limits_{\varepsilon\rightarrow 0^{+}}\displaystyle\int_{|\xi|>\varepsilon }\frac{f(x+\xi)-f(\xi)}{|\xi|^{n+2s}}d\xi,
\end{array}
\right.
\end{equation}
 or
 \begin{equation}\label{1.6}
 I_{s}f(x)=(-\Delta)^{s}f(x):=\frac{C_{n,s}}{2}\int_{\mathbb{R}^{n} }\frac{2f(x)-f(x+y)-f(x-y)}{|y|^{n+2s}}dy,\end{equation}
 where
 the parameter $s$ is a real number with $0<s<1$, $\mbox{P.V.}$ is a commonly used abbreviation for "in the principle value sense" (as defined by the latter equation), and $C_{n,s}$  in \eqref{1.5} is a dimensional constant that depends on $n$ and $s$, precisely given by
  \begin{equation}\label{1.7}
\displaystyle C_{n,s}:=\left(\int_{\mathbb{R}^{n} }\frac{1-\cos(\zeta_{1})}{|\zeta|^{n+2s}}d \zeta\right)^{-1}=\frac{2^{2s}\Gamma\left((n+2s)/2\right)}{\pi^{\frac{n}{2}}\Gamma(1-s)}.\end{equation}
 In particular, $C_{n,s}$ is a normalization constant; see, for example, \cite{Landkof}.\\
 \indent In the whole space, if $f\in  \mathscr{S}(\mathbb{R}^{n}) $, let
  $\Lambda^{\gamma}=(-\Delta)^{s}$ with $\gamma=2s$, and
\begin{equation*}
\widehat{\Lambda^{2s}f} (\xi)=\widehat{(-\Delta)^{s}f} (\xi)=|\xi|^{2s}\widehat{f}(\xi),\end{equation*}
the domain of definition of the fractional Laplacian, $\mathcal{D}\left(\Lambda^{s}\right) $ is endowed with a natural norm $\|\cdot\|_{\mathcal{D}\left(\Lambda^{s}\right) }$ and is a Hilbert space. The norm of $u$ in $\mathcal{D}\left(\Lambda^{s}\right) $ is defined by
\begin{equation}\label{1.8}
 \|u\|_{\mathcal{D}\left(\Lambda^{s}\right)}:=
\left\|\Lambda^{s}u\right\|_{L^2( \mathbb{R}^n) }. \end{equation}
It should be pointed out that in the whole space, if any function $\psi\in  \mathscr{S}(\mathbb{R}^{n}) $, $\mathcal{D}\left(\Lambda^{s}\right) $ is equivalent to the fractional Sobolev space $\dot{H}^{s}(\mathbb{R}^n) $, defined as the completion of $C^{\infty}_{0}(\mathbb{R}^n)$ with the norm
\begin{equation}\label{1.9}
\|\psi\|_{\dot{H}^{s}(\mathbb{R}^n)}=\left(\int_{\mathbb{R}^n}\left|\widehat{\psi}
\right|^{2}d\xi\right)^{\frac{1}{2}}=\left\|(-\Delta)^{\frac{s}{2}}\psi\right\|_{L^{2}(\mathbb{R}^n)}.
\end{equation}
On the other hand, the norm $\|u\|_{H^{s}(\mathbb{R}^n)}$ in the fractioal Laplacian Sobolev space $H^{s}(\mathbb{R}^n)$ is represented as
\begin{equation}\label{1.10}
\|\mathbf{u}\|^{2}_{H^{s}(\mathbb{R}^n)}:=2C(n,s)^{-1}\left\|\Lambda^{s} u\right\|_{L^2( \mathbb{R}^n) }^2+\|u\|_{L^2(\mathbb{R}^n)}^2.
\end{equation}
In particular, the norm of $\mathcal{D}\left(\Lambda^{2}\right)=\mathcal{D}(-\Delta)$ is equivalent to the $H^{2}(\mathbb{R}^n)$ norm.\\
\indent The study of the nonlocal Camassa-Holm equations \eqref{1.1} is hindered by a lack of explicit information on the kernels of these nonlocal operators appearing in them. However, we can obtain various estimates by using these tools adapted to the Dirichlet boundary case: such as the C\'{o}rdoba-C\'{o}rdoba inequality \cite{Cordoba,Cordoba1}, a nonlinear lower bound in the spirit of \cite{P-V}, and commutator estimates \cite{Kato-Ponce}. In particular, we will utilize several fractional-type interpolation inequalities such as the nonlocal version of Ladyzhenskaya's inequalities \cite{Ladyzhenskaya-1,Ladyzhenskaya-2,Ladyzhenskaya-3,Ladyzhenskaya-4}, fractional Gagliardo-Nirenberg-Sobolev inequality and fractional Leibniz rule.\\
 \indent There has been a vast literature devoted to proving the existence, uniqueness and regularity issues for the Camassa-Holm equations \eqref{1.3} \cite{Clayton-08,Chen-Liu,de Lellis-Kapp,Escher-Yin,Hakkaev,Misiolek,Perrollaz-Vin,Tan-Yin,Wu-Yin,Yan-Li-Zhang}.  Bjorland and Schonbek \cite{Clayton-08} studied the decay and existence of solutions to \eqref{1.3}. Chen et al. in \cite{Chen-Liu} investigated the oscillation-induced blow-up to the modified Camassa-Holm equations with linear dispersion. de Lellis et al. in \cite{de Lellis-Kapp} considered the low-regularity solutions for the periodic Camassa-Holm equations. Escher and Yin in \cite{Escher-Yin} analyzed a kind of initial-boundary value problems of the Camassa-Holm equations.
  Hakkaev in \cite{Hakkaev} obtained the local well-posedness for a generalized Camassa-Holm equations.
  Misiolek in \cite{Misiolek} discussed classical solutions of a periodic Camassa-Holm equations. Perrollaz in \cite{Perrollaz-Vin} dealt with an initial boundary value problem for the Camassa-Holm equations on an interval. Tan and Yin in \cite{Tan-Yin} established the global periodic conservative solutions for a periodic modified two-component Camassa-Holm equations. Wu and Yin  in \cite{Wu-Yin} showed global existence and blow-up phenomena for the weakly dissipative Camassa-Holm equations. Yan et al. in \cite{Yan-Li-Zhang} took into account the Cauchy problem for a generalized Camassa-Holm equations in Besov space.\\
  \indent In stark contrast to those works for the Camassa-Holm equations \eqref{1.3} mentioned as above in recent decades, little has been known concerning these issues for the Camassa-Holm equations with space-fractional derivative viscosity \eqref{1.1} in the literature despite that non-standard diffusions are very natural also for these problems. Compared with the standard Laplacian, space-fractional derivatives, notably the fractional Laplacian, are more challenging for achiving the existence, uniqueness and regularity properties of solutions owing to the vector integral expression and nonlocal property. It is well-known that fractional Laplacian $(-\Delta)^{s}$ is a spatial integro-differential operator, which can be used to describe the spatial nonlocality and power law behaviors in various science and engineering problems. In the recent two decades, fractional Laplacian has been utilized to model energy dissipation of acoustic propagation in human tissue \cite{Chen-Holm}, turbulence diffusion \cite{Chen-W}, contaminant transport in ground water \cite{Pang-Chen}, non-local heat conduction \cite{Bobaru,Elia,Mongiovi}, and electromagnetic fields on fractals \cite{Tarasov}. In particular, it is expected that these results for the equations with fractional Laplacian viscosity would depend on $s$ and $\nu$. Some related problems have been previously considered in the literature motivated by some important equations appearing in fluid mechanics.\\
  \indent Since a fundamental problem for both Euler and Navier-Stokes equations is the regularity, in particular, proving global regularity for the 3D Navier-Stokes equations is one of the most challenging outstanding problems in nonlinear analysis, we consider in this paper a similar partial differential equation with fractional Laplacian viscosity. In stark contrast to the problem on the regularity for the Camassa-Holm equations without any nonlocal term \eqref{1.3}, it seems fair to say that extremely little is known about the regularity of the solutions to the nonlocal equations \eqref{1.1}. Indeed, to our best knowledge, the only work for the nonlocal Camassa-Holm equations established by Gui-Liu\cite{Gui-Liu}, in which some results of the nonlocal Camassa-Holm equations in one space dimension have been obtained is as follows:
   \\[0.30cm]
\indent $\bullet$ Global well-posedness and blow-up of solutions to the Camassa-Holm equations with fractional dissipation under the supercritical case: $\gamma\in \left[\frac{1}{2},1\right)$.
\\[0.2cm]
\indent $\bullet$ The zero filter limit of the Camassa-Holm equation with fractional dissipation, as well as the possible blow-up of solutions under the subcritical case: $0\leq \gamma <\frac{1}{2}$.
\\[0.30cm]
\indent The goal in this article is to
           investigate the effect of diffusion and its criticality for the system
         \eqref{1.1} in the same spirit as that in the study of Euler and Navier-Stokes equations. Specifically, the aim of this paper is to establish the regularity of solutions to the Camassa-Holm equations with fractional Laplacian viscosity  \eqref{1.1}. The main difficulty lies in proving some a prior estimates for the nonlocal viscosity. To achieve this, we need to establish some fractional Sobolev type estimates and bilinear estimates for fractional derivatives. Especially, for the critical case $\displaystyle s=\frac{n}{4}$ with $n=2,3$, we will make extra efforts for acquiring the expected estimates as obtained in the case $\displaystyle  \frac{n}{4}<s<1$. Fortunately, by the aid of the fractional Leibniz rule, we achieve the existence, uniqueness and regularity to the Camassa-Holm equations with nonlocal viscosity by the energy method and a bootstrap argument, which rely crucially on the fractional Laplacian viscosity. In particular, under the critical case $s=\dfrac{n}{4}$, the nonlocal version of Ladyzhenskaya's inequality is skillfully used, and the smallness of initial data in several Sobolev spaces is required to gain the desired results concernig existence, uniqueness and regularity.\\
    \indent Before going further, let us describe the notation we shall use in this paper.
    \\[0.3cm]
{\bf Notation}
 \\[0.3cm]
$\mathscr{S}(\mathbb{R}^{n})$ denotes the Schwartz calss. The $i^{th}$ component of $\mathbf{v}\cdot \nabla \mathbf{u}^{T}$ is denoted by $\left(\mathbf{v}\cdot \nabla \mathbf{u}^{T}\right)_{i}=\sum\limits_{j=1}^{n}v_{j}\partial_{i}u_{j}$. Let
     $\displaystyle\left \langle \mathbf{u},\mathbf{v}\right\rangle=\int_{ \mathbb{R}^{n} }\mathbf{u}\cdot \mathbf{v}dx$ and $\displaystyle\Sigma =\left\{\phi \in C^{\infty}_{0}( \mathbb{R}^{n} )|\nabla\cdot \phi=0\right\}.$ $H^{m}_{0}(\mathbb{R}^{n})$ denotes the completion of $C^{\infty}_{0}( \mathbb{R}^{n} )$ in the norm $\|\cdot\|_{H^{m}(\mathbb{R}^{n})}$. We denote by $L^{p}( \mathbb{R}^{n} )$ the standard Lebesgue space, and
    $L^{p}_{\sigma}( \mathbb{R}^{n} )$ the completion of $\Sigma$ in the norm $\|\cdot\|_{L^{p}( \mathbb{R}^{n} )}$. The completion of $\Sigma$ under the $\mathcal{D}\left(\Lambda^{s}\right)( \mathbb{R}^{n} )$-norm is denoted by $\mathcal{D}_{\sigma}\left(\Lambda^{s}\right)( \mathbb{R}^{n} )$ and $\left(\mathcal{D}_{\sigma}\left(\Lambda^{s}\right)\right)'( \mathbb{R}^{n} )$ is its dual space. The completion of $\Sigma$ under the $H^{m}( \mathbb{R}^{n} )$-norm will be denoted by  $H^{m}_{\sigma}( \mathbb{R}^{n} )$ and
      $(H^{m}_{\sigma}( \mathbb{R}^{n} ))'$  be the corresponding dual space. $\mathcal{F}(\phi)$ or $\hat{\phi}$ denotes the Fourier transform of a function $\phi$, with $\mathcal{F}^{-1}(\phi)$ or $ \breve{\phi} $  the inverse Fourier transform. For $a\lesssim b$, we mean that there is a uniform constant $C$, which may be different on different lines, such that $a\leq C b$. \hfill$\Box$
      \\[0.3cm]
    \indent The main results in the present paper  are the low-order regularity, high-order regularity and uniqueness of weak solutions to the Cauchy problem  \eqref{1.1}-\eqref{1.2}.\\
 \begin{thm}[Low-order regularity]\label{t1.1}\rm
  Let $ n=2,3 $. Assume that
  \\[0.2cm]
 \indent (1) for $\displaystyle\frac{n}{4}< s < 1$, the initial data  $\mathbf{v}_{0}\in L^{2}_{\sigma}(\mathbb{R}^{n})$,\\
  and \\
  \indent(2) for $\displaystyle s=\frac{n}{4}$, the initial data $\mathbf{v}_{0}\in L^{2}_{\sigma}(\mathbb{R}^{n})$, and in addition, there exists an $\varepsilon^{*}=\varepsilon^{*}(\alpha,\nu,n)$ sufficiently small such that $\|\mathbf{v}_{0}\|_{L^{2}(\mathbb{R}^{n})}\leq \varepsilon^{*}$.
  \\[0.2cm]
Then there exists a weak solution on $[0,T)$ to the Cauchy problem \eqref{1.1}-\eqref{1.2} in the sense of definition \ref{d2.3} satisfies the estimate \eqref{2.3}. In particular, the following two bounds hold:
\begin{equation}\label{1.11}
\|\mathbf{v} \|_{L^{\infty}\left([0,T],L^{2}_{\sigma}(\mathbb{R}^{n})\right)}+\left\|\Lambda^{s} \mathbf{v} \right\|_{L^{2}\left([0,T],L^{2}_{\sigma}(\mathbb{R}^{n})\right)}\leq C\left(n,s,\alpha,\nu,\|\mathbf{v}_{0}\|_{L^{2}(\mathbb{R}^{n})}\right),
\end{equation}
\begin{equation}\label{1.12}
\left\|\partial_{t}\mathbf{v} \right\|_{L^{2}\left([0,T],\mathcal{B}'\right)}\leq
C\left(n,s,\alpha,\nu,\|\mathbf{v}_{0}\|_{L^{2}(\mathbb{R}^{n})}\right),\end{equation}
\\[0.1cm]
where $\displaystyle \mathcal{B}=\left\{\begin{array}{ll}
\mathcal{D}_{\sigma}\left(\Lambda^{s}\right)(\mathbb{R}^{n})~~&\displaystyle\hbox{for}~~\frac{n}{4}< s < 1,
\\[0.3cm]
H_{\sigma}^{\frac{n}{4}}(\mathbb{R}^{n})~~&\displaystyle\hbox{for}~~s=\frac{n}{4}.\end{array}\right.$
$$\qquad\eqno\Box$$
 \end{thm}
\begin{thm}[High-order regularity]\label{t1.2}\rm
Let $n=2,3$.  Assume that
\\[0.2cm]
 \indent (1) for $\displaystyle\frac{n}{4}<s < 1$,  $\mathbf{v}_{0}\in \mathcal{D}_{\sigma}\left(\Lambda^{M} \right)(\mathbb{R}^{n} )$, $M\geq 0$,\\
  and\\
   \indent (2) for $\displaystyle s=\frac{n}{4}$, $\mathbf{v}_{0}\in H_{\sigma}^{M}(\mathbb{R}^{n})(\mathbb{R}^{n} )$, $M\geq 0$, and in addition, there exists an $\varepsilon^{**}=\varepsilon^{**}(\alpha,\nu,n)$ sufficiently small such that $\left\| \mathbf{v}_{0}\right\|_{H_{0}^{M}(\mathbb{R}^{n}) }\leq \varepsilon^{**}$.
   \\[0.2cm]
     Then the weak solutions to the Cauchy problem \eqref{1.1}-\eqref{1.2} constucted in Theorem \ref{t1.1} satisfy the following bound
\begin{equation}\label{1.13}
\left\|\partial^{k}_{t}\nabla^{m}\mathbf{v}\right\|^{2}_{L^{2}( \mathbb{R}^{n} )}+\nu\int_{0}^{T}
\left\|\partial^{k}_{t}\nabla^{m}\Lambda^{s} \mathbf{v}\right\|^{2}_{L^{2}( \mathbb{R}^{n} )}dt\leq C\left(n,s,\alpha,\nu,\| \mathbf{v}_{0}\|_{\mathcal{A}}\right)\end{equation}
for all $m+2ks\leq M$, where  $\displaystyle \mathcal{A}=\left\{\begin{array}{ll}
\mathcal{D}\left(\Lambda^{M}\right)(\mathbb{R}^{n})~~&\hbox{for}~~\frac{n}{4}< s < 1,
\\[0.3cm]
H_{0}^{M}(\mathbb{R}^{n})~~&\hbox{for}~~s=\frac{n}{4},\end{array}\right.$ and $m$ and $k$ are both non-negative integers \hfill$\Box$
 \end{thm}
 \begin{thm}[Uniqueness]\label{t1.3}\rm
Let $n=2,3$ and let $\displaystyle\frac{n}{4}\leq s < 1$. The weak solution to the Cauchy problem \eqref{1.1}-\eqref{1.2} constucted in Theorem \ref{t1.1} is unique.\hfill$\Box$
 \end{thm}
In Theorem \ref{t1.1} and Theorem \ref{t1.2}, the statements of Definition \ref{d2.3} and the energy estimate \eqref{2.3} will be given in Section 2. Thanks to Theorem \ref{t1.1} and Theorem \ref{t1.2}, one claims a Corollary concerning some estimates for the filter corresponding to the Camassa-Holm equations with fractional Laplacian viscosity \eqref{1.1}.
\\
\begin{cor}\label{c1.4}
 For $n=2,3$ and $\displaystyle\frac{n}{4}\leq s < 1$, let $(\mathbf{v},\mathbf{u})$ be the solution to the Cauchy problem \eqref{1.1}-\eqref{1.2} constructed in Theorem \ref{t1.1} and Theorem \ref{t1.2}. Then for all $m+2ks\leq M$, there holds\\
\begin{equation}\label{1.14}
\left\|\partial^{k}_{t}\nabla^{m}  \mathbf{u}\right\|^{2}_{L^{2}(\mathbb{R}^{n})}
+2\alpha^2\left\|\partial^{k}_{t}\nabla^{m+1} \mathbf{u}\right\|^{2}_{L^{2}(\mathbb{R}^{n})}
+\alpha^4\left\|\partial^{k}_{t}\nabla^{m+2} \mathbf{u}\right\|^{2}_{L^{2}(\mathbb{R}^{n})}=\left\|\partial^{k}_{t}\nabla^{m} \mathbf{v}\right\|^{2}_{L^{2}(\mathbb{R}^{n})},
\end{equation}
\begin{equation}\label{1.15}
 \left\|\partial^{k}_{t}\nabla^{m}\mathbf{u}\right\|^{2}_{L^{n}(\mathbb{R}^{n})}
+\left\|\partial^{k}_{t}\nabla^{m}\Lambda^{s} \mathbf{u}\right\|^{2}_{L^{n}(\mathbb{R}^{n})}
+\left\|\partial^{k}_{t}\nabla^{m+1}\mathbf{u}\right\|^{2}_{L^{n}(\mathbb{R}^{n})}
 \lesssim\left\|\partial^{k}_{t}\nabla^{m} \mathbf{v}\right\|^{2}_{L^{2}(\mathbb{R}^{n})},
 \end{equation}
\begin{equation}\label{1.16}
\displaystyle\left\|\partial^{k}_{t}\nabla^{m}\mathbf{u}\right\|^{2}_{L^{n}(\mathbb{R}^{n})}
+\nu\int^{t}_{0}\left\|\partial^{k}_{t}\nabla^{m}\Lambda^{s}\mathbf{u}\right\|^{2}_{L^{n}(\mathbb{R}^{n})}ds
\displaystyle\leq C\left(n,s,\alpha,\nu,m,k,\|\mathbf{v}_{0}\|_{\mathcal{D}\left(\Lambda^{M} \right)(\mathbb{R}^{n})}\right).
\end{equation}
Here, $m$ and $k$ are both non-negative integers.\hfill$\Box$
\end{cor}
\indent Before going further, we collect some facts on the fractional Sobolev spaces $W^{s,p}(\mathbb{R}^{n})$ and $H^{s}(\mathbb{R}^{n})$, as well as the definition of the fractional fractional Laplacian \cite{Nezza}.\\
\begin{defn}\label{d1.5}
Let $s\in (0,1)$. For any $p\in [1,\infty)$, we define $W^{s,p}(\mathbb{R}^{n})$ as follows
\begin{equation}\label{1.17}
W^{s,p}(\mathbb{R}^{n}):=\left\{u\in L^{p}(\mathbb{R}^{n}):~~\frac{|u(x)-u(y)|}{|x-y|^{\frac{n}{p}+s}}\in L^{p}(\mathbb{R}^{n}\times\mathbb{R}^{n})\right\},
\end{equation}
i.e., an intermediary Banach space between $L^{p}(\mathbb{R}^{n})$ and $W^{1,p}(\mathbb{R}^{n})$, endowed with the natural norm
\begin{equation}\label{1.18}
\|u\|_{W^{s,p}(\mathbb{R}^{n})}:=\left(\int_{\mathbb{R}^{n}}|u|^{p}dx
+\int_{\mathbb{R}^{n}}\int_{\mathbb{R}^{n}}\frac{|u(x)-u(y)|^{p}}{|x-y|^{n+sp}}
dxdy\right)^{\frac{1}{p}},
\end{equation}
where the term
\begin{equation}\label{1.19}
[u]_{W^{s,p}(\mathbb{R}^{n})}:=\left(\int_{\mathbb{R}^{n}}\int_{\mathbb{R}^{n}}
\frac{|u(x)-u(y)|^{p}}{|x-y|^{n+sp}}
dxdy\right)^{\frac{1}{p}}
\end{equation}
is the so-called Gagliardo (semi) norm of $u$. However, there is another case for $s\in (1,\infty)$ and $s$ is not an integer. In this case, we write $s=m+m'$, where $m$ is an integer and $m'\in (0,1)$. The space $W^{s,p}(\mathbb{R}^{n})$ consists of those equivalence classes of functions $u\in W^{m,p}(\mathbb{R}^{n})$ whose distributional derivatives $D^{\alpha}u$, with $|\alpha|=m$, belong to $W^{m',p}(\mathbb{R}^{n})$, namely
\begin{equation*}
W^{s,p}(\mathbb{R}^{n}):=\left\{u\in W^{m,p}(\mathbb{R}^{n}):~~D^{\alpha}u\in W^{m',p}(\mathbb{R}^{n})~~\mbox{for~~any}~~\alpha~s.t.~~|\alpha|=m\right\},
\end{equation*}
 and this is a Banach space with respect to the norm
 \begin{equation}\label{1.20}
\|u\|_{W^{s,p}(\mathbb{R}^{n})}:=\left( \left\|u\right\|^{p}_{W^{m,p}(\mathbb{R}^{n})}+\sum\limits_{|\alpha|=m}
\left\|D^{\alpha}u\right\|^{p}_{W^{m',p}(\mathbb{R}^{n})}\right)^{\frac{1}{p}}.
\end{equation}
 Clearly, if $s=m$ is an integer, the space $W^{s,p}(\mathbb{R}^{n})$ coincides with the Sobolev space $W^{m,p}(\mathbb{R}^{n})$.\\
 \indent Since for any $s>0$, the space $C^{\infty}_{0}(\mathbb{R}^{n})$ of smooth functions with compact support is dense in $W^{s,p}(\mathbb{R}^{n})$, we have $W^{s,p}_{0}(\mathbb{R}^{n})=W^{s,p}(\mathbb{R}^{n})$, where $W^{s,p}_{0}(\mathbb{R}^{n})$ denotes the closure of $C^{\infty}_{0}(\mathbb{R}^{n})$ in the space $W^{s,p}(\mathbb{R}^{n})$. However, for $s\in (0,1)$ and $p=2$, the fractional Sobolev spaces $W^{s,2}(\mathbb{R}^{n})$ and $W^{s,2}_{0}(\mathbb{R}^{n})$ turn out to be Hilbert spaces, and we usually label $W^{s,2}(\mathbb{R}^{n})=H^{s}(\mathbb{R}^{n})$ and $W^{s,2}_{0}(\mathbb{R}^{n})=H^{s}_{0}(\mathbb{R}^{n})$. That is,
\begin{equation}\label{1.21}
H^{s}(\mathbb{R}^{n}):=\left\{u\in L^{2}(\mathbb{R}^{n}):~~\frac{|u(x)-u(y)|}{|x-y|^{\frac{n}{2}+s}}\in L^{2}(\mathbb{R}^{n}\times\mathbb{R}^{n})\right\},
\end{equation}
i.e., an intermediary Hilbert space between $L^{2}(\mathbb{R}^{n})$ and $H^{1}(\mathbb{R}^{n})$, endowed with the natural norm
\begin{equation}\label{1.22}
\|u\|_{H^{s}(\mathbb{R}^{n})}:=\left(\int_{\mathbb{R}^{n}}|u|^{2}dx
+\int_{\mathbb{R}^{n}}\int_{\mathbb{R}^{n}}\frac{|u(x)-u(y)|^{2}}{|x-y|^{n+2s}}
dxdy\right)^{\frac{1}{2}},
\end{equation}
where the term
\begin{equation}\label{1.23}
[u]_{H^{s}(\mathbb{R}^{n})}:=\left(\int_{\mathbb{R}^{n}}\int_{\mathbb{R}^{n}}
\frac{|u(x)-u(y)|^{2}}{|x-y|^{n+2s}}
dxdy\right)^{\frac{1}{2}}
\end{equation}
is the so-called seminorm of $u$. \\
\indent There is alternative definition of the space $H^{s}(\mathbb{R}^{n})$ via the Fourier transform. For any real $s\geq 0$, we may define
\begin{equation}\label{1.24}
\widehat{H}^{s}(\mathbb{R}^{n}):=\left\{u\in L^{2}(\mathbb{R}^{n}):~~
\int_{\mathbb{R}^{n}}\left(1+|\xi|^{2s}\right)|\mathcal{F}u(\xi)|^{2}d\xi<\infty
\right\}.
\end{equation}
In the same manner, for $s<0$ there is an analogous definition for $H^{s}(\mathbb{R}^{n})$:
 \begin{equation}\label{1.25}
H^{s}(\mathbb{R}^{n}):=\left\{u\in \mathscr{S}'(\mathbb{R}^{n}):~~
\int_{\mathbb{R}^{n}}\left(1+|\xi|^{2 }\right)^{s}|\mathcal{F}u(\xi)|^{2}d\xi<\infty
\right\}.
\end{equation}
 On the other hand, let $s\in (0,1)$ and let $(-\Delta)^{s}:~~\mathscr{S}\rightarrow L^{2}(\mathbb{R}^{n})$ be the fractional Laplacian operator defined by \eqref{1.6}. Then\\
 \indent (1)  For any $u\in \mathscr{S}$,
  \begin{equation}\label{1.26}
  (-\Delta)^{s}u=\mathcal{F}^{-1}\left[|\xi|^{2s}\left(\mathcal{F}u\right)\right],~~\forall \xi \in \mathbb{R}^{n}.
\end{equation}
 \indent (2)  The fractional Sobolev space $H^{s}(\mathbb{R}^{n})$ defined in \eqref{1.21} coincides with $\widehat{H}^{s}(\mathbb{R}^{n})$ defined in \eqref{1.24}. In particular, for any $u\in H^{s}(\mathbb{R}^{n})$
 \begin{equation}\label{1.27}
 [u]^{2}_{H^{s}(\mathbb{R}^{n})}=2C(n,s)^{-1}
 \int_{\mathbb{R}^{n}}|\xi|^{2s}\left|\mathcal{F}u(\xi)\right|^{2}d\xi,
 \end{equation}
 where $C(n,s)$ is defined by \eqref{1.7}.\\
 \indent (3)  For $u\in H^{s}(\mathbb{R}^{n})$,
 \begin{equation}\label{1.28}
 [u]^{2}_{H^{s}(\mathbb{R}^{n})}=2C(n,s)^{-1}
 \left\|(-\Delta)^{\frac{s}{2}}u\right\|_{L^{2}(\mathbb{R}^{n})}^{2},
 \end{equation}
 where  $C(n,s)$ is defined by \eqref{1.7}.\hfill$\Box$
\end{defn}
\indent At the end of this section, we want to make some remarks on fractional Sobolev spaces and the fractional Laplacian in an open bounded set with smooth boundary.
\begin{rem}\label{r1.6} \rm
 Let $\Omega\subset \mathbb{R}^{n}$ be an open bounded set with smooth boundary. We reccall the following three conclusions.
 \\[0.3cm]
 (I)\quad We can identify $\mathcal{D}\left(\Lambda^{s} \right)(\Omega)$ with usual Sobolev spaces (\cite{Constantin3},\cite{Lions-J}):
 \begin{equation*}
\mathcal{D}\left(\Lambda^{s} \right)(\Omega)=\left\{
\begin{array}{ll}
\displaystyle H^{s}_{0}(\Omega)&\mbox{if}~s\in \left(\frac{1}{2},1\right],
\\[0.2cm]
\displaystyle H^{\frac{1}{2}}_{00}(\Omega)=\left\{u\in H^{\frac{1}{2}}_{0}(\Omega):~u/\sqrt{d(x)}\in L^{2}(\Omega)\right\}&\mbox{if}~s=\frac{1}{2},
\\[0.2cm]
 \displaystyle H^{s} (\Omega) &\mbox{if}~s\in \left[0,\frac{1}{2}\right).
\end{array}
\right.
\end{equation*}
Here, $d(x)=\mbox{dist}~(x,\partial \Omega)$.
\\[0.2cm]
 (II)\quad It follows from \cite{P-I} that the continuous embedding $\mathcal{D}\left(\Lambda^{s} \right)(\Omega)\subset H^{s}(\Omega)$ holds for all $s\geq0$.
 \\[0.2cm]
 (III)\quad We mention here some facts in \cite{P-I}. Let
  $$Q=\Omega\times \mathbb{R}_{+}=\left\{(x,z)~|~ x\in \Omega,~z>0\right\},$$
  $$H^{1}_{0L}(Q)=\left\{v\in H^{1} (Q)~ |~ v(x,z)=0,~ x\in \partial\Omega,~z>0\right\},$$
 $$V_{0}(\Omega)=\left\{f ~|~ \exists v\in H^{1}_{0L}(Q),~ f(x)=v(x,0),~ x\in \Omega \right\}.$$
Recall from \cite{Cabre-Tan} that, on one hand,
$$V_{0}(\Omega)=\left\{f\in  H^{\frac{1}{2}} (\Omega)~| ~
\int_{\Omega}\frac{f^{2}(x)}{d(x)}dx<\infty\right\}$$
with norm $\displaystyle\left\|f\right\|_{V_{0}}^{2}=\left\|f\right\|_{H^{\frac{1}{2}}(\Omega)}^{2}
+\int_{\Omega}\frac{f^{2}(x)}{d(x)}dx$. On the other hand, $V_{0}(\Omega)=\mathcal{D}\left(\Lambda^{\frac{1}{2}} \right)(\Omega)$, i.e., $$V_{0}(\Omega)=\left\{f\in   L^{2} (\Omega)| f=\sum\limits_{j}f_{j}\omega_{j},~\sum\limits_{j}
\lambda^{\frac{1}{2}}_{j}f^{2}_{j}<\infty\right\}$$
 with equivalent norm $\displaystyle\|f\|^{2}_{\mathcal{D}\left(\Lambda^{\frac{1}{2}} \right)}=\sum\limits_{j}
\lambda^{\frac{1}{2}}_{j}f^{2}_{j}=\left\|\Lambda ^{\frac{1}{2}}f\right\|^{2}_{L^{2}(\Omega)}$.\\
In particular, a smooth function with compact support satisfies the above three conclusions (I), (II) and (III) with $\Omega=\mathbb{R}^{n}$. \hfill$\Box$
 \end{rem}
      \indent The paper is organized as follows. In Section 2 we collect some preliminaries. In Section 3 we present the proof of the low-order regularity (Theorem \ref{t1.1}). We prove Theorem \ref{t1.2} in Section 4. Theorem \ref{t1.3} and Corollary \ref{c1.4} will be shown in the last Section (Section 5).
\section{Preliminaries}
In this section, we collect some preliminaries.\\
 \begin{lem} \label{l2.1}
 For $n=2,3$, let $\mathbf{u}$ and $\mathbf{v}$ be two smooth divergence free functions with compact support. Then
\begin{equation}\label{2.1}
\left\{
\begin{array}{ll}
\mathbf{u}\cdot\nabla \mathbf{v}+ \sum \limits^{n}_{j=1} v_{j}\nabla u_{j}=-\mathbf{u}\times(\nabla\times \mathbf{v})+\nabla(\mathbf{v}\cdot \mathbf{u}),
\\[0.2cm]
\left\langle \mathbf{u}\cdot\nabla \mathbf{v},\mathbf{u}\right\rangle+\left\langle \mathbf{v}\cdot\nabla \mathbf{u}^{T},\mathbf{u}\right\rangle=0,
\\[0.2cm]
\left\langle \mathbf{u}\times(\nabla\times \mathbf{v}),\mathbf{u}\right\rangle=0.
\end{array}
\right.
\end{equation}
\end{lem}
 {\bf Proof.} \quad Direct calculation leads to  \eqref{2.1}.\hfill$\Box$\\
  \begin{lem} \label{l2.2}
 For $n=2,3$, let $(\mathbf{u},\mathbf{v})$ be the solution of \eqref{1.1}-\eqref{1.2} with compact support. Then there holds
\begin{equation}\label{2.2}
\dfrac{1}{2}\dfrac{d}{dt}\left(\langle \mathbf{u},\mathbf{u}\rangle+\alpha^{2}\left\langle \nabla \mathbf{u},\nabla \mathbf{u}\right\rangle\right)  +\nu\left(\left\langle\Lambda^{s} \mathbf{u},\Lambda^{s} \mathbf{u}\right\rangle+\alpha^{2}\left\langle\nabla\Lambda^{s}  \mathbf{u},\nabla\Lambda^{s} \mathbf{u}\right\rangle\right) =0,\end{equation}
\begin{equation}\label{2.3}
\left.
\begin{array}{ll}
\displaystyle\left \|\mathbf{u}(\cdot,t)\right\|_{L^2(\mathbb{R}^{n})}^{2}+\alpha^{2}\left\|\nabla \mathbf{u}(\cdot,t)\right\|_{L^2(\mathbb{R}^{n})}^{2}
\\[0.3cm]
\qquad\qquad\quad\displaystyle+2\nu\int_0^t \left\|\Lambda^{s} \mathbf{u}(\cdot,t)\right\|_{L^2(\mathbb{R}^{n})}^{2}dt
+2\nu\alpha^{2}\int_0^t\left\|\nabla\Lambda^{s} \mathbf{u}(\cdot,t)\right\|_{L^2(\mathbb{R}^{n})}^{2}dt
\\[0.4cm]
\qquad\qquad\quad\displaystyle\leq\|\mathbf{u}_{0}\|_{L^2(\mathbb{R}^{n})}^{2}+\alpha^{2}\|\nabla \mathbf{u}_{0}\|_{L^2(\mathbb{R}^{n})}^{2}.
\end{array}
\right.
\end{equation}
\end{lem}
{\bf Proof.}\quad Note that $\nabla\cdot \mathbf{v}=\nabla\cdot \mathbf{u}=0$ with $n=2,3$, making inner product with $\mathbf{u}$ on the first equation in \eqref{1.1} gives rise to
$$\langle \mathbf{v}_{t},\mathbf{u}\rangle+\left\langle \mathbf{u}\cdot \nabla\mathbf{ v}+\mathbf{v}\cdot \nabla \mathbf{u}^{T},\mathbf{u}\right\rangle+\left\langle\nabla p,\mathbf{u}\right\rangle+\nu\left\langle(-\Delta)^{s} \mathbf{v},\mathbf{u}\right\rangle=0.$$
Integrating by parts yields
$$\langle \mathbf{v}_{t},\mathbf{u}\rangle+\nu\left\langle(-\Delta)^{s}\mathbf{v},  \mathbf{u}\right\rangle=\langle \mathbf{v}_{t},\mathbf{u}\rangle+\nu\left\langle \Lambda^{s}_{D}\mathbf{v},  \Lambda^{s}_{D}\mathbf{u}\right\rangle=0.$$
This leads to
\begin{eqnarray*}
&&
\left\langle \mathbf{u}_{t}-\alpha^{2} \Delta \mathbf{u}_{t},\mathbf{u}\right\rangle+\nu\left\langle(-\Delta)^{s}(\mathbf{u}-\alpha^{2}\Delta \mathbf{u}), \mathbf{u}\right\rangle\\
&&\qquad=\dfrac{1}{2}\dfrac{d}{dt} \left(\left\langle \mathbf{u},\mathbf{u}\right\rangle+\alpha^{2}\left\langle\nabla \mathbf{u},\nabla \mathbf{u}\right\rangle\right)
\\
&&\qquad\quad +\nu\left(\left\langle\Lambda^{s} \mathbf{u},\Lambda^{s}  \mathbf{u}\right\rangle+\alpha^{2}\left\langle\nabla\Lambda^{s}  \mathbf{u},\nabla\Lambda^{s}  \mathbf{u}\right\rangle\right)\\
&&\qquad=0.
\end{eqnarray*}
\qquad \hfill$\Box$\\

\indent We introduce the following notion of weak solutions to \eqref{1.1}-\eqref{1.2}. \\
\begin{defn}\label{d2.3}
For $\displaystyle\frac{n}{4}\leq s< 1$ and $n=2,3$, let $T>0$. A weak solution to \eqref{1.1}-\eqref{1.2} is a pair of functions $(\mathbf{v},\mathbf{u})$ such that
\begin{eqnarray*}
&&
\mathbf{v}\in \emph{L}^{\infty}\left([0,T];\emph{L}_{\sigma}^{2}(\mathbb{R}^{n})\right)\cap \emph{L}^{2}\left([0,T];\mathcal{D}_{\sigma}\left(\Lambda^{s}\right)(\mathbb{R}^{n})\right),
\\[0.2cm]
&&\partial_{t}\mathbf{v}\in\emph{L}^{2}\left([0,T];\mathcal{B}' \right),
\\[0.2cm]
&&\mathbf{u}\in\emph{L}^{\infty}\left([0,T]; \emph{H}_\sigma^{2} (\mathbb{R}^{n})\right)\cap \emph{L}^{2}\left([0,T];\mathcal{D}_\sigma\left(\Lambda^{2+s}\right)(\mathbb{R}^{n})\right),
\\[0.2cm]
 &&\mathbf{v}(0,x)=\mathbf{v}_{0}(x).
\end{eqnarray*}
Here, $\displaystyle \mathcal{B}=\left\{\begin{array}{ll}
\mathcal{D}_{\sigma}\left(\Lambda^{s}\right)(\mathbb{R}^{n})~~&\hbox{for}~~\frac{n}{4}< s < 1,
\\[0.3cm]
H_{\sigma}^{\frac{n}{4}}(\mathbb{R}^{n})~~&\hbox{for}~~s=\frac{n}{4}.\end{array}\right.$\\
\indent In addition, for every ${\phi}\in\emph{L}^{2}\left([0,T];\mathcal{E}\right)$ with ${\phi(T)}=0$, there holds
\begin{eqnarray*}
\vspace{5pt}
&&-\int_{0}^{T}\left\langle\mathbf{ v},\partial_{t}\phi\right\rangle ds+\int_{0}^{T}\left\langle \mathbf{u}\cdot\nabla \mathbf{v},\phi \right\rangle ds +\int_{0}^{T}\left\langle\phi\cdot\nabla \mathbf{u},\mathbf{v}\right\rangle ds+\nu\int_{0}^{T}\left\langle \Lambda ^{s}\mathbf{v},\Lambda ^{s} \phi\right\rangle ds\\
\\
&&\qquad\qquad\qquad\quad={\left\langle \mathbf{v}_{0},\phi(0)\right\rangle}.
\end{eqnarray*}
In particular, for $t\in [0,T]$ there holds
$$\left\langle\mathbf{u},\phi\right\rangle+\alpha^{2}\left\langle\nabla \mathbf{u}, \nabla\phi\right\rangle=\langle \mathbf{v},\phi\rangle. $$
Here, $\displaystyle \mathcal{E}=\left\{\begin{array}{ll}
\mathcal{D}_{\sigma}\left(\Lambda\right)(\mathbb{R}^{n})~~&\hbox{for}~~\frac{n}{4}< s < 1,
\\[0.3cm]
H_{\sigma}^{1}(\mathbb{R}^{n})~~&\hbox{for}~~s=\frac{n}{4}.\end{array}\right.$
$$\eqno \Box$$
\end{defn}
\indent It should be pointed out that if $\psi$ and $\phi$ belong to the Schwartz class $\mathscr{S}(\mathbb{R}^{n})$, Definition \eqref{1.4} of the fractional Laplacian together with Plancherel's theorem yields
$$\int_{\mathbb{R}^{n} }(-\Delta)^{s}\psi\phi dx=\int_{\mathbb{R}^{n} }|\xi|^{2s}\widehat{\psi}\widehat{\phi} d\xi=\int_{\mathbb{R}^{n} }(-\Delta)^{\frac{s}{2}}\psi(-\Delta)^{\frac{s}{2}}\phi dx.$$
\indent We then claim the following estimates.\\
\begin{lem}\label{l2.4}
For $\displaystyle\frac{n}{4}\leq s < 1$ and $n=2,3$, let $ \mathbf{v}\in \mathcal{D}_\sigma\left(\Lambda ^{s}\right)( \mathbb{R}^{n} )$ and $\mathbf{u}\in \mathcal{D}_\sigma\left(\Lambda ^{2+s}\right)(\mathbb{R}^{n})$ satisfy the Helmholtz equations
\begin{equation}\label{2.4}
\mathbf{u}-\alpha^2 \Delta \mathbf{u}=\mathbf{v}.
\end{equation}
This gives
$$\qquad\qquad\Lambda^{s}\mathbf{u}-\alpha^2 \Delta \Lambda^{s} \mathbf{u}=\Lambda^{s} \mathbf{v}.$$
In addition, there holds
\begin{equation}\label{2.5}
\left\{
\begin{array}{ll}
\left\|\mathbf{u}\right\|^{2}_{L^{2}(\mathbb{R}^{n})}+2\alpha^2\left\|\nabla \mathbf{u}\right\|^{2}_{L^{2}(\mathbb{R}^{n})}+\alpha^4 \left\|\Delta \mathbf{u}\right\|^{2}_{L^{2}(\mathbb{R}^{n})}=\left\|\mathbf{v}\right\|^{2}_{L^{2}(\mathbb{R}^{n})},
\\[0.3cm]
{ \left \|\Lambda^{s} \mathbf{u}\right\|^{2}_{L^{2}(\mathbb{R}^{n})}+2\alpha^2\left\|\nabla \Lambda^{s}  \mathbf{u}\right\|^{2}_{L^{2}(\mathbb{R}^{n})}+\alpha^4 \left\|\Delta \Lambda^{s} \mathbf{u}\right\|^{2}_{L^{2}(\mathbb{R}^{n})}=\left\|\Lambda^{s} \mathbf{v}\right\|^{2}_{L^{2}(\mathbb{R}^{n})},}
\\[0.3cm]
 \left\| \mathbf{u}\right\| _{L^{n}(\mathbb{R}^{n})}\leq C\left\|\mathbf{v}\right\|_{L^{2}(\mathbb{R}^{n})},~~
 \left\| \nabla \mathbf{u}\right\| _{L^{n}(\mathbb{R}^{n})}\leq C\left\|\mathbf{v}\right\|_{L^{2}(\mathbb{R}^{n})},
 \\[0.3cm]
\left\| \nabla \mathbf{u}\right\| _{L^{\frac{n}{s}}(\mathbb{R}^{n})}\leq C\left\|\Lambda^{s}\mathbf{v}\right\|_{L^{2}(\mathbb{R}^{n})},~~\left\| \Delta \mathbf{u}\right\| _{L^{\frac{n}{s}}(\mathbb{R}^{n})}\leq C\left\|\Lambda^{s}\mathbf{v}\right\|_{L^{2}(\mathbb{R}^{n})}. \end{array}
\right.
\end{equation}
\end{lem}
{\bf Proof.} Note that $\displaystyle\frac{s}{n}\geq\frac{1}{2}-\frac{s}{n} $ for $\displaystyle\frac{n}{4}\leq s<1$ with $n=2,3$, by virtue of the Gagliardo-Nirenberg-Sobolev inequality, we complete the proof of this lemma.\hfill$\Box$\\
\\
  \indent The following Lemma concerns the nonlocal version of the known inequalities established in these works of Ladyzhenskaya-Shkoller-Seregin \cite{Ladyzhenskaya-1,Ladyzhenskaya-2,Ladyzhenskaya-3,Ladyzhenskaya-4}.\\
 \begin{lem} \label{l2.5}
 For $n=2,3$ and $\mathbf{u}(x)\in H^{1}_{0}(\mathbb{R}^{n})$, $\forall$ $\varepsilon>0$, the following inequalities hold:
$$
\left\{
\begin{array}{lll}
\|\mathbf{u}\|^{2}_{L^{4}(\mathbb{R}^{n})}\leq \varepsilon \|\nabla \mathbf{u}\|^{2}_{L^{2}(\mathbb{R}^{n})}+\varepsilon^{-1}\| \mathbf{u}\|^{2}_{L^{2}(\mathbb{R}^{n})} &\mbox{for}~n=2,
\\
\\
\|\mathbf{u}\|^{2}_{L^{4}(\mathbb{R}^{n})}\leq 3^{-\frac{1}{4}}\sqrt{2\varepsilon }\|\nabla \mathbf{u}\|^{2}_{L^{2}(\mathbb{R}^{n})}+\sqrt{2}(3^{ \frac{5}{2}}\varepsilon)^{ -\frac{1}{6}}\| \mathbf{ u}\|^{2}_{L^{2}(\mathbb{R}^{n})}&\mbox{for}~n=3.
\end{array}
\right.\eqno(E-1)
$$
 The above inequalities (E-1) can be generalized to the following nonlocal version (fractional power Sobolev-type) estimates.\\
 \\
\indent $\heartsuit$ For $\displaystyle\frac{n}{4}<s< 1$ and $\mathbf{u}\in \mathcal{D}\left(\Lambda^{s}\right )(\mathbb{R}^{n})$, the following inequalities hold:
$$
\left\{
\begin{array}{lll}
\|\mathbf{u}\|^{2}_{L^{4}(\mathbb{R}^{n})}\leq \varepsilon \left\|\Lambda^{s} \mathbf{u}\right\|^{2}_{L^{2}(\mathbb{R}^{n})}+\varepsilon^{-1}\| \mathbf{u}\|^{2}_{L^{2}(\mathbb{R}^{n})} &\mbox{for}~n=2,\\
\\
\|\mathbf{u}\|^{2}_{L^{4}(\mathbb{R}^{n})}\leq  C(s)\varepsilon \left\|\Lambda^{s}  \mathbf{u}\right\|^{2}_{L^{2}(\mathbb{R}^{n})}+C(\varepsilon) \| \mathbf{u}\|^{2}_{L^{2}(\mathbb{R}^{n})} &\mbox{for}~n=3.
\end{array}
\right.\eqno(E-2)
$$
Here, $\varepsilon$, $C(s)$ and $C(\varepsilon)$ are constants; $C(s)$ depends only on spatial dimensions and $s$, and  $C(\varepsilon)=O(\varepsilon^{-\frac{1}{3}})$.
\\[0.3cm]
\indent $\heartsuit$ For the critical case $\displaystyle s=\frac{n}{4} $ and $\mathbf{u}\in \mathcal{D}\left(\Lambda^{\frac{n}{4}} \right)(\mathbb{R}^{n})$, the following inequality holds:
$$\|\mathbf{u}\|^{2}_{L^{4}(\mathbb{R}^{n} )}\leq C \left(\left\|\Lambda^{\frac{n}{4}}  \mathbf{u}\right\|^{2}_{L^{2}(\mathbb{R}^{n})}+\| \mathbf{u}\|^{2}_{L^{2}(\mathbb{R}^{n})}\right)~~\mbox{for}~n=2,3.\eqno(E-3)$$
Here, $C$ is a constant depending only on spatial dimension $n$.\hfill$\Box$
\end{lem}
\begin{rem}\label{r2.6}
For $\displaystyle\frac{n}{4}\leq s < 1$ and $n=2,3$, the fractional stationary stokes equation
$$\nu(-\Delta)^{s}\mathbf{u}+\nabla p=\mathbf{v}~~\mbox{in}~~\mathbb{R}^{n},\eqno(2.9)$$
has a solution $\mathbf{u}\in \mathcal{D}_{\sigma}\left(\Lambda^{s} \right)(\mathbb{R}^{n})$ for each $\mathbf{v}\in L^{2}_{\sigma}(\mathbb{R}^{n})$. Here, $\nu$ is a positive constant (the kinematic coefficient of viscosity), the "pressure" term $p$ is not known a priori but is determined by a posteriori from the solution itself, $(-\Delta)^{s}$ is the fractional Laplacian defined in Section 1. Solving (2.9) defines a continuous operator $L^{2}_{\sigma}(\mathbb{R}^{n})\rightarrow \mathcal{D}_{\sigma}\left(\Lambda^{s} \right)(\mathbb{R}^{n})$. Combining this with the compact inclusion $\mathcal{D}_{\sigma}\left(\Lambda^{s}\right )(\mathbb{R}^{n})\hookrightarrow L^{2}_{\sigma}(\mathbb{R}^{n})$  achieves a compact and self-adjoint operator $A:$ $L^{2}_{\sigma}(\mathbb{R}^{n})\rightarrow L^{2}_{\sigma}(\mathbb{R}^{n})$ ( the stokes operator) .\hfill$\Box$
\end{rem}

\indent We now collect some known estimates for the vector-valued fractional Leibniz rule.\\
\begin{lem} \label{l2.7}
Let $\Lambda^{s}=(-\Delta)^{\frac{s}{2}}$ be the standard Riesz potential of order $s\in \mathbb{R}$. Then the following two conclusions hold.\\
 \indent(I)  (see \cite{FGO})\quad Let $s_{1}, s_{2}\in [0,1]$, $s=s_{1}+ s_{2}$, and $p,~ p_{1},~p_{2}\in (1,\infty) $  such
  that $\displaystyle\frac{1}{p}=\frac{1}{p_{1}}+\frac{1}{p_{2}}$. Then the following bilinear estimate holds for all $f,g\in \mathscr{S}(\mathbb{R}^{n})$, $n\geq 1$:
  $$\left\|\Lambda^{s}(fg)-f\Lambda^{s}g-g\Lambda^{s}f\right\|_{L^{p}(\mathbb{R}^{n})}\leq C\left\|\Lambda^{s_{1}}f\right\|_{L^{p_{1}}(\mathbb{R}^{n})}
  \left\|\Lambda^{s_{2}}g\right\|_{L^{p_{2}}(\mathbb{R}^{n})}. \eqno (2.10)$$
  \\
 \indent (II) (see \cite{HTW})\quad Let $s>\max \left(0,\frac{n}{p }-n\right)$, or $s$ be a positive even integer, $\displaystyle\frac{1}{2 }<p<\infty$, $1<p_{1},p_{2}\leq\infty$, and $\displaystyle\frac{1}{p}=\frac{1}{p_{1}}+\frac{1}{p_{2}}$, then
 \\
  $$\left\|\Lambda^{s}\left(fg\right)\right\|_{L^{p}(\mathbb{R}^{n})}\leq C\left\|\Lambda^{s }f\right\|_{L^{p_{1}}(\mathbb{R}^{n})}
  \left\| g\right\|_{L^{p_{2}}(\mathbb{R}^{n})}+\left\|  f\right\|_{L^{p_{1}}(\mathbb{R}^{n})}
  \left\| \Lambda^{s}g\right\|_{L^{p_{2}}(\mathbb{R}^{n})}.\eqno (2.11) $$
  $$\quad\eqno \Box$$
\end{lem}
 \begin{rem}\label{r2.8}
For $0<s<1$ and $n=1$, Kenig, Ponce, and Vega in \cite{KPV} obtained the similar estimates for fractional derivatives as those in Lemma \ref{l2.7}:
  $$\left\|\Lambda^{s}(fg)-f\Lambda^{s}g-g\Lambda^{s}f\right\|_{L^{p}(\mathbb{R})}\leq C\left\|\Lambda^{s_{1}}f\right\|_{L^{p_{1}}(\mathbb{R})}
  \left\|\Lambda^{s_{2}}g\right\|_{L^{p_{2}}(\mathbb{R})}, $$
   $$\left\|\Lambda^{s}(fg)-f\Lambda^{s}g-g\Lambda^{s}f\right\|_{L^{p}(\mathbb{R})}\leq C\left\|g\right\|_{L^{\infty}(\mathbb{R})}
  \left\|\Lambda^{s }f\right\|_{L^{p }(\mathbb{R})},\qquad\qquad $$
  where $p, p_{1},p_{2}\in (1,\infty) $ and $\displaystyle\frac{1}{p}=\frac{1}{p_{1}}+\frac{1}{p_{2}}$, $0<s=s_{1}+s_{2}<1$, $s_{1}, s_{2}\geq 0$.\hfill$\Box$
\end{rem}
\begin{lem} \label{l2.9}
  For $\displaystyle\frac{n}{4}< s<1$ with $n=2,3$, let \vspace*{2 ex} $\displaystyle A=\frac{n}{2}+1-2s$. Direct calculation gives\\
 \vspace*{2 ex}
 (I) \quad $s\geq 1-s,~\dfrac{n-2}{2n} <\dfrac{2s-1}{n}<\dfrac{1}{n},  ~\dfrac{2s-1}{n}=\dfrac{1}{2}-\dfrac{A}{n},$\\
  \vspace*{2 ex}
  \indent\quad $\dfrac{n}{2}-1<A< 1<\dfrac{n}{2}<2~~\mbox{for}~~n=3$,\\
  \vspace*{2 ex}
  \indent\quad $\dfrac{n}{2}-1<A< 1=\dfrac{n}{2}<2~~\mbox{for}~~n=2$. \\
  \\
  (II) \quad In view of $\displaystyle\frac{n}{2s-1}=\frac{2n}{n-2A}$, for $\dfrac{n}{4}\leq s<1$ with $n=3$, and $\dfrac{n}{4}< s<1$ with $n=2$, we have the following fractional Sobolev-type continuous embedding between $\mathcal{D}\left(\Lambda^{A}\right)(\mathbb{R}^{n})$ and  $L^{\frac{n}{2s-1}}(\mathbb{R}^{n})$:
   $$\mathcal{D}\left(\Lambda^{A} \right)(\mathbb{R}^{n})\hookrightarrow L^{ n/(2s-1)}(\mathbb{R}^{n}).\eqno \Box$$
\end{lem}
\begin{lem} \label{l2.10}
For $s\in (0,1)$ and $n=2,3$, the inclusion $\mathcal{D}\left(\Lambda^{s}\right)(\mathbb{R}^{n})\rightarrow L^{2}_{0}(\mathbb{R}^{n})$ is compact and the embedding $\mathcal{D}\left(\Lambda^{\gamma}\right)(\mathbb{R}^{n})\hookrightarrow H^{\gamma}(\mathbb{R}^{n})$ is continuous for all $\gamma\geq 0$.\hfill$\Box$
\end{lem}
 \section{Low-order regularity}
\indent In this section, we prove the low-order regularity result (Theorem \ref{t1.1}).\\
{\bf Proof of Theorem \ref{t1.1}}.\\
\indent From Theorem 3.1 in Gan-Lin-Tong \cite{GLT}, it follows that there exists a weak solution to the Cauchy problem \eqref{1.1}-\eqref{1.2} in the sense of definition \ref{d2.3}. In particular, it satisfies the bound \eqref{1.11}.
\\[0.3cm]
\indent We then prove \eqref{1.12}.
\\[0.3cm]
\indent Let $\phi\in L^{2}\left([0,T],\mathcal{D}_{\sigma}\left(\Lambda^{s}\right)(\mathbb{R}^{n})\right)$ for $\dfrac{n}{4}< s<1$, and  $\phi\in L^{2}\left([0,T],H_{\sigma}^{\frac{n}{4}}(\mathbb{R}^{n})\right)$ for $s=\dfrac{n}{4}$. Making inner product for the first equation in \eqref{1.1} with $\phi$ yields
  \begin{equation}\label{3.1}
\left\langle\partial_{t}\mathbf{v} ,\phi \right\rangle+\left\langle \mathbf{u} \cdot\nabla \mathbf{v} ,\phi\right\rangle+\left\langle \mathbf{v} \cdot \nabla \mathbf{u}^{T},\phi\right\rangle =-\nu \left\langle (-\Delta)^{s} \mathbf{v} ,\phi\right\rangle.\end{equation}
This leads to
 \begin{equation}\label{3.2}
\left.
\begin{array}{ll}
  \left|\left\langle\partial_{t}\mathbf{v} ,\phi\right\rangle\right|  &\leq \left|\left\langle \mathbf{u} \cdot\nabla \mathbf{v} ,\phi\right\rangle\right|+\left|\left\langle \mathbf{v} \cdot \nabla \mathbf{u} ^{T},\phi\right\rangle\right|+\nu\left|\left\langle (-\Delta)^{s} \mathbf{v} ,\phi\right\rangle\right|
  \\[0.3cm]
& := A_{1}+A_{2}+A_{3}.
\end{array}
\right.
\end{equation}
 We then estimate $A_{1},A_{2},A_{3}$  one by one through considering two cases:\\
\\
 \indent \quad {\bf Case 1}~~ $\dfrac{n}{4}< s<1$ for $n= 2,3$;\\
 \indent \quad {\bf Case 2}~~ $s= \dfrac{n}{4}$ for $n=2,3$.\\
 $\displaystyle\Large{\diamondsuit}$ We first consider  {\bf Case 1}~~  $\dfrac{n}{4}< s<1$ for $n=2, 3$.\\
 \indent In this case, a straightforward computation shows that
 \begin{equation}\label{3.3}
\left\{
\begin{array}{ll}
\dfrac{n}{s}=\dfrac{2n}{n-2(n-2s)/2},~~\dfrac{1}{2}=\dfrac{n-2}{2}<\dfrac{n-2s}{2}
  <\dfrac{3}{4}~~
  &\mbox{for}~~n=3,
  \\[0.4cm]
  0=\dfrac{n-2}{2}<\dfrac{n-2s}{2}< \dfrac{1}{2}~~
  &\mbox{for}~~n=2,
  \\[0.3cm]
 B=\dfrac{n-2s}{2}+1-s=\dfrac{n}{2}+1-2s, ~~\dfrac{n}{2}-1<B< 1 ~~  &\mbox{for}~~n=2,3.
  \end{array}
\right.
\end{equation}
 Thanks to Lemma \ref{l2.4}, Lemma \ref{l2.7} and Lemma \ref{l2.9}, there holds
 \begin{equation}\label{3.4}
 A_{1} = \left|\left\langle \mathbf{u} \cdot\nabla \mathbf{v} ,\phi\right\rangle\right|  \lesssim\left\|\Lambda^{s} \mathbf{v} \right\|_{L^{2} (\mathbb{R}^{n})}\left\| \Lambda^{1-s} (\mathbf{u} \phi)\right\|_{L^{2} (\mathbb{R}^{n})}.
 \end{equation}
 The term $\left\| \Lambda^{1-s}(\mathbf{u}\phi)\right\|_{L^{2} (\mathbb{R}^{n})}$ in \eqref{3.4} can be estimated as
\begin{equation}\label{3.5}
\left.
\begin{array}{ll}
\left\| \Lambda^{1-s}(\mathbf{u}\phi)\right\|_{L^{2} (\mathbb{R}^{n})}&\lesssim \left\| \Lambda^{1-s} (\mathbf{u} \phi)-\mathbf{u} \Lambda^{1-s} \phi-\phi\Lambda^{1-s}  \mathbf{u} \right\|_{L^{2} ( \mathbb{R}^{n})}
  \\[0.3cm]
 &\quad+\left\| \mathbf{u} \Lambda^{1-s}  \phi \right\|_{L^{2} (\mathbb{R}^{n})}+\left\| \phi\Lambda^{1-s}  \mathbf{u} \right\|_{L^{2} (\mathbb{R}^{n} )}.
\end{array}
\right.
\end{equation}
By interpolation inequality, note that the first estimate in \eqref{2.5},
by (I) of Lemma \ref{l2.7}, Lemma \ref{l2.9} and \eqref{3.3}, for $\displaystyle\frac{1}{2}=\frac{1}{n/s}+\frac{1}{2n/(n-2s)}$ and $0<1-s<s<1$, the three terms on the right hand side of \eqref{3.5} can be bounded as follows.
\\[0.2cm]
\begin{equation}\label{3.6}
\left.
\begin{array}{ll}
 \displaystyle\left\| \Lambda^{1-s}(\mathbf{u}\phi)-\mathbf{u}\Lambda^{1-s} \phi-\phi\Lambda^{1-s} \mathbf{u}\right \|_{L^{2} (\mathbb{R}^{n})}
 \\[0.3cm]
 \displaystyle\qquad\qquad\lesssim\left\| \Lambda^{1-s}\mathbf{u}  \right \|_{L^{\frac{n}{s}} (\mathbb{R}^{n})} \left\| \phi  \right\|_{L^{\frac{2n}{n-2s}} (\mathbb{R}^{n})}
 \\[0.3cm]
\displaystyle\qquad\qquad \lesssim\left\| \Lambda^{\frac{n}{2}+1-2s} \mathbf{u}   \right\|_{L^{2} (\mathbb{R}^{n})} \left\| \Lambda^{s} \phi  \right\|_{L^{2} (\mathbb{R}^{n})}
\\[0.3cm]
\displaystyle\qquad\qquad \lesssim\left\|   \mathbf{u}   \right\|^{\frac{1}{2}}_{L^{2} (\mathbb{R}^{n})}\left\|  \nabla \mathbf{u}   \right\|^{\frac{1}{2}}_{L^{2} (\mathbb{R}^{n})}\left \| \Lambda^{s} \phi  \right\|_{L^{2} (\mathbb{R}^{n})},
 \end{array}
\right.
\end{equation}
\\[0.1cm]
\begin{equation}\label{3.7}
\left.
\begin{array}{ll}
\displaystyle\left\| \mathbf{u} \Lambda^{1-s}  \phi  \right\|_{L^{2} (\mathbb{R}^{n})}
 \lesssim\left\|  \mathbf{u}   \right\|_{L^{\frac{n}{2s-1}} (\mathbb{R}^{n})} \left\| \Lambda^{1-s} \phi  \right\|_{L^{\frac{2n}{n-2(2s-1)}} (\mathbb{R}^{n})}
  \\[0.3cm]
\displaystyle \qquad\qquad\qquad \lesssim\left\|  \mathbf{u}    \right\|_{L^{\frac{2n}{n-2\left(\frac{n}{2}+1-2s\right)}} (\mathbb{R}^{n})} \left\| \Lambda^{s} \phi  \right\|_{L^{2} (\mathbb{R}^{n})}
\\[0.4cm]
\displaystyle \qquad\qquad\qquad \lesssim\left\|  \Lambda^{\frac{n}{2}+1-2s} \mathbf{u}   \right\|_{L^{2} (\mathbb{R}^{n})} \left\| \Lambda^{s} \phi  \right\|_{L^{2} (\mathbb{R}^{n})}
\\[0.3cm]
\displaystyle\qquad\qquad\qquad\lesssim\left\|   \mathbf{u}   \right\|^{\frac{1}{2}}_{L^{2} (\mathbb{R}^{n})}\left\|  \nabla \mathbf{u}   \right\|^{\frac{1}{2}}_{L^{2} (\mathbb{R}^{n})}\left\| \Lambda^{s} \phi  \right\|_{L^{2} (\mathbb{R}^{n})},
 \end{array}
\right.
\end{equation}
and
\begin{equation}\label{3.8}
\left.
\begin{array}{ll}
\left\|  \phi \Lambda^{1-s} \mathbf{u}  \right\|_{L^{2} (\mathbb{R}^{n})}
 \lesssim\left\| \phi  \right\|_{L^{\frac{2n}{n-2s}} (\mathbb{R}^{n})}\left\| \Lambda^{1-s}  \mathbf{u}   \right\|_{L^{\frac{n}{s}} (\mathbb{R}^{n})}
 \\[0.3cm]
 \displaystyle\qquad\qquad\qquad \lesssim\left\| \nabla \mathbf{u } \right \|_{L^{2} (\mathbb{R}^{n})} \left\| \Lambda^{s} \phi  \right\|_{L^{2} (\mathbb{R}^{n})}.
 \end{array}
\right.
\end{equation}
\\[0.2cm]
Combining \eqref{3.4} with \eqref{3.5}, \eqref{3.6}, \eqref{3.7} and \eqref{3.8} gives rise to
\begin{equation}\label{3.9}
\left.
\begin{array}{ll}
A_{1}&\lesssim \left(\left\|   \mathbf{u}   \right\|^{\frac{1}{2}}_{L^{2} (\mathbb{R}^{n})}\left\|  \nabla \mathbf{u}   \right\|^{\frac{1}{2}}_{L^{2} (\mathbb{R}^{n})}+\left\| \nabla \mathbf{u} \right\|_{L^{2} (\mathbb{R}^{n})}\right) \left\| \Lambda^{s}\mathbf{ v}  \right\|_{L^{2} (\mathbb{R}^{n})}\left\| \Lambda^{s} \phi \right\|_{L^{2} (\mathbb{R}^{n})}
 \\[0.3cm]
 &\lesssim\left\|\mathbf{ v} \right \|_{L^{2} (\mathbb{R}^{n})} \left\| \Lambda^{s} \mathbf{v} \right \|_{L^{2} (\mathbb{R}^{n})}\left\| \Lambda^{s} \phi \right\|_{L^{2} (\mathbb{R}^{n})}.
\end{array}
\right.
\end{equation}
 We next bound $A_{2}$ and $A_{3}$. Thanks to \eqref{2.5} and \eqref{3.3}, we have
\begin{equation}\label{3.10}
\left.
\begin{array}{ll}
\displaystyle A_{2}=\left|\left\langle \mathbf{v} \cdot \nabla \mathbf{u} ^{T},\phi\right\rangle\right|
 \\[0.3cm]
\displaystyle\quad\lesssim \left\| \nabla  \mathbf{u }\right\|_{L^{2} (\mathbb{R}^{n})}\left\| \mathbf{v } \phi\right\|_{L^{2} (\mathbb{R}^{n})}
 \\[0.3cm]
\displaystyle\quad\lesssim\left\| \nabla  \mathbf{u} \right\|_{L^{2} (\mathbb{R}^{n})}\left\| \mathbf{v} \right\|_{L^{\frac{2n}{n-2s}} (\mathbb{R}^{n})}\left\| \phi \right\|_{L^{\frac{n}{s}} (\mathbb{R}^{n})}
 \\[0.3cm]
\displaystyle\quad \lesssim\left\| \mathbf{v}\right \|_{L^{2} (\mathbb{R}^{n})}\left\| \Lambda^{s} \mathbf{v } \right\|_{L^{2} (\mathbb{R}^{n})}\left\| \Lambda^{s} \phi\right\|_{L^{2} (\mathbb{R}^{n})}.
\end{array}
\right.
\end{equation}
\begin{equation}\label{3.11}
A_{3}= \nu\left|\left\langle (-\Delta)^{s} \mathbf{v} ,\phi\right\rangle\right| \lesssim  \left\| \Lambda^{s} \mathbf{v} \right\|_{L^{2} (\mathbb{R}^{n})}\left\| \Lambda^{s} \phi\right\|_{L^{2} (\mathbb{R}^{n})}.\qquad\qquad
\end{equation}
Collecting \eqref{3.9}, \eqref{3.10} and \eqref{3.11} together yields that for  $\displaystyle\frac{n}{4}< s<1$ with $n=2, 3$,
\begin{equation}\label{3.12}
\left|\left\langle \partial_{t}\mathbf{v} ,\phi\right\rangle\right|\lesssim\left\| \mathbf{v} \right\|_{L^{2} (\mathbb{R}^{n})} \left\| \Lambda^{s} \mathbf{v}\right \|_{L^{2} (\mathbb{R}^{n})} \left\| \Lambda^{s} \phi\right\|_{L^{2} (\mathbb{R}^{n})}.
\end{equation}
\vspace{0.3cm}
As  $\phi\in L^{2}\left([0,T],\mathcal{D}_{\sigma}\left(\Lambda^{s}\right)(\mathbb{R}^{n})\right)$  can be chosen arbitrarily, using H\"{o}lder's inequality,
\eqref{3.12} together with \eqref{1.8}, \eqref{1.11} and \eqref{2.3} concludes that for  $\displaystyle\frac{n}{4}< s<1$ with $n= 2,3$,
\begin{equation}\label{3.13}
 \left\|\partial_{t}\mathbf{v} \right\|_{L^{2}\left([0,T],\left[\mathcal{D}_{\sigma}\left(\Lambda^{s}\right)(\mathbb{R}^{n})\right]'\right)}\leq C\left(n,s,\alpha,\nu,\|\mathbf{v}_{0}\|_{L^{2}(\mathbb{R}^{n})}\right).\end{equation}\\
$\displaystyle\diamondsuit$ We next consider  {\bf Case 2}~~ $ \displaystyle s= \frac{n}{4}$ for $n=2,3$.
\\[0.3cm]
\indent  In this case, direct calculation yields
 \begin{equation}\label{3.14}
\left.
\begin{array}{ll}
A_{1} = \left|\left\langle \mathbf{u} \cdot\nabla\mathbf{ v} ,\phi\right\rangle\right|  \lesssim\left\| \Lambda^{\frac{n}{4}} \mathbf{v} \right\|_{L^{2} (\mathbb{R}^{n})}\left\|\Lambda^{1-n/4} (\mathbf{u} \phi) \right\|_{L^{2} (\mathbb{R}^{n})}.\end{array}
\right.
\end{equation}
In particular, $\displaystyle\left\|\Lambda^{1-n/4}  (\mathbf{u} \phi) \right\|_{L^{2} (\mathbb{R}^{n})}$ can be bounded as
\begin{equation}\label{3.15}
\left.
\begin{array}{ll}
\displaystyle\left\|\Lambda^{1-n/4}  (\mathbf{u} \phi) \right\|_{L^{2} (\mathbb{R}^{n})}
  \lesssim\left\|\Lambda^{1-n/4}  (\mathbf{u} \phi)-u \Lambda^{1-n/4}  \phi  -\phi\Lambda^{1-n/4} \mathbf{u}   \right\|_{L^{2} (\mathbb{R}^{n})}
  \\[0.3cm]
 \qquad \qquad \qquad\qquad\displaystyle+\left\| \mathbf{u} \Lambda^{1-n/4} \phi\right\|_{L^{2} (\mathbb{R}^{n})}+\left\|\phi\Lambda^{1-n/4}  \mathbf{u} \right\|_{L^{2} (\mathbb{R}^{n})}.
\end{array}
\right.
\end{equation}
 We first estimate the first term on the right hand side of inequality \eqref{3.15}.  Recall (I) of Lemma \ref{l2.7}, one achieves for $\displaystyle 0<s_{1}<1-\frac{n}{4}$,
\begin{equation}\label{3.16}
\left.
\begin{array}{ll}
\displaystyle\left\|\Lambda^{1-n/4} (\mathbf{u }\phi)-\mathbf{u} \Lambda^{1-n/4}   \phi  -\phi\Lambda^{1-n/4} \mathbf{u}    \right\|_{L^{2} (\mathbb{R}^{n})}
\\[0.3cm]
 \displaystyle\qquad\qquad\lesssim\left\| \Lambda^{1-n/4-s_{1}}  \mathbf{u}  \right\|_{L^{\frac{2n}{n-2\left(n/4+s_{1}\right)}} (\mathbb{R}^{n})}
\left\|\Lambda^{ s_{1}} \phi\right\|_{L^{\frac{2n}{n-2\left(n/4-s_{1}\right)}} (\mathbb{R}^{n})}
\\[0.3cm]
\displaystyle\qquad\qquad\lesssim\left\| \nabla \mathbf{u}  \right\|_{L^{2}  (\mathbb{R}^{n})}
\left\|\Lambda^{n/4} \phi \right\|_{L^{2} (  \mathbb{R}^{n})}
\\[0.3cm]
\displaystyle\qquad\qquad \lesssim\| \mathbf{v}  \|_{L^{2} (\mathbb{R}^{n})}
\left\|\Lambda^{n/4} \phi \right\|_{L^{2} (  \mathbb{R}^{n})},
\end{array}
\right.
\end{equation}
where $\displaystyle s_{1}\in \left(0,\frac{1}{2}\right), 1-\frac{n}{4}-s_{1}\in \left(0,\frac{1}{2}\right)$, $\displaystyle\frac{1}{2}=\frac{1}{p_{1}}+\frac{1}{p_{2}}$ with $\displaystyle p_{1}, p_{2}\in (1,\infty)$, $\displaystyle p_{1}=\frac{2n}{n-2(n/4+s_{1})}$, $\displaystyle p_{2}=\frac{2n}{n-2(n/4-s_{1})}$. In the same manner, by virtue of \eqref{2.3}, Lemma \ref{l2.10}, Agmon's inequality and interpolation inequality, the second and the third terms on the right hand side of inequality \eqref{3.15} can be bounded as follows:
\begin{equation}\label{3.17}
\left.
\begin{array}{ll}
\displaystyle\left\| \mathbf{u} \Lambda^{1-n/4} \phi\right\|_{L^{2} (\mathbb{R}^{n})}+\left\|\phi\Lambda^{1-n/4}  \mathbf{u}\right\|_{L^{2} (\mathbb{R}^{n})}
\\[0.3cm]
\displaystyle
 \qquad\qquad\lesssim\|   \mathbf{u}  \|_{L^{\infty} (\mathbb{R}^{n})}
\left\|\Lambda^{ 1-n/4} \phi\right\|_{L^{2} (\mathbb{R}^{n})}+  \| \phi \|_{L^{4} (\mathbb{R}^{n})}
\left\|\Lambda^{ 1-n/4}  \mathbf{u} \right\|_{L^{4} (\mathbb{R}^{n})}
\\[0.3cm]
\displaystyle
\qquad\qquad \lesssim  \|   \mathbf{u} \|^{\frac{1}{2}}_{H^{1} (\mathbb{R}^{n})}
\|   \mathbf{u}  \|^{\frac{1}{2}}_{H^{2} (\mathbb{R}^{n})}\left\| \phi\right\|^{\frac{1}{2}}_{L^{2} (\mathbb{R}^{n})}
\left\|\Lambda^{  n/4} \phi\right\|^{\frac{1}{2}}_{L^{2} (\mathbb{R}^{n})}
\\[0.3cm]
\displaystyle
 \qquad\qquad\quad +  \left(\left\|\Lambda^{ 1-n/4}  \mathbf{u}\right \|^{2} _{L^{2} (\mathbb{R}^{n})}+\left\|\Lambda^{1-n/4+n/4}    \mathbf{u} \right\|^{2}  _{L^{2} (\mathbb{R}^{n})}\right)^{\frac{1}{2}}
 \\[0.3cm]
\displaystyle
\qquad\qquad\qquad \times  \left( \|\phi\|^{2} _{L^{2}(\mathbb{R}^{n})}+\left\|\Lambda^{n/4}  \phi\right\|^{2}_{L^{2} (\mathbb{R}^{n})}\right)^{\frac{1}{2}}
\\[0.3cm]
\displaystyle \qquad\qquad \lesssim  \| \mathbf{v}  \|^{2}_{L^{2} (\mathbb{R}^{n})}\left( \|\phi\| _{L^{2}(\mathbb{R}^{n})}+\left\|\Lambda^{n/4}  \phi\right\|^{2}_{L^{2} (\mathbb{R}^{n})}\right)^{\frac{1}{2}}.
\end{array}
\right.
\end{equation}
From \eqref{3.14},\eqref{3.15},\eqref{3.16} and \eqref{3.17} it follows that
\begin{equation}\label{3.18}
\displaystyle A_{1}\lesssim  \left\| \mathbf{v} \right\|_{L^{2} (\mathbb{R}^{n})}\left\| \Lambda^{\frac{n}{4}} \mathbf{v} \right\|_{L^{2} (\mathbb{R}^{n})}\left( \|\phi\|^{2}  _{L^{2}(\mathbb{R}^{n})}+\left\|\Lambda^{n/4}  \phi\right\|^{2} _{L^{2} (\mathbb{R}^{n})}\right)^{\frac{1}{2}}.\end{equation}
Note that the assumption of this theorem (the smallness of initial data), calculations like those employed in $A_{1}$ above imply
\begin{equation}\label{3.19}
\left.
\begin{array}{ll}
A_{2}=\left|\left\langle \mathbf{v} \cdot \nabla \mathbf{u} ^{T},\phi\right\rangle\right|
\\[0.3cm]
\displaystyle
~\quad \lesssim \left\| \nabla  \mathbf{u }\right\|_{L^{2} (\mathbb{R}^{n})}\left\| \mathbf{v } \phi\right\|_{L^{2} ( \mathbb{R}^{n})}
 \\[0.3cm]
 ~\quad \lesssim \left\| \nabla  \mathbf{u} \right\|_{L^{2} (\mathbb{R}^{n})}\| \mathbf{v} \|_{L^{4} (\mathbb{R}^{n})}\left\|  \phi\right\|_{L^{4} (\mathbb{R}^{n})}
 \\[0.3cm]
 ~\quad \lesssim  \left\| \nabla  \mathbf{u} \right\|_{L^{2} (\mathbb{R}^{n})}
\left(\| \mathbf{v}  \|^{2}_{L^{2} (\mathbb{R}^{n})}+\left\| \Lambda^{n/4} \mathbf{v} \right\|^{2}_{L^{2} (\mathbb{R}^{n})}\right)^{\frac{1}{2}}
 \\[0.3cm]
 ~\qquad \qquad  \times\left(\left\| \phi \right\|^{2}_{L^{2} (\mathbb{R}^{n})}+\left\| \Lambda^{n/4} \phi\right\|^{2}_{L^{2} (\mathbb{R}^{n})}\right)^{\frac{1}{2}}
  \\[0.3cm]
 ~\quad \lesssim \| \mathbf{v} \|_{L^{2} (\mathbb{R}^{n})}
\left\| \Lambda^{n/4} \mathbf{v}  \right\|_{L^{2} (\mathbb{R}^{n})}
\left( \|\phi\|^{2} _{L^{2}(\mathbb{R}^{n})}+\left\|\Lambda^{n/4}  \phi\right\|^{2}_{L^{2} (\mathbb{R}^{n})}\right)^{\frac{1}{2}}, \end{array}
\right.
\end{equation}
and
\begin{equation}\label{3.20}
A_{3}=\nu\left|\left\langle (-\Delta)^{n/4} \mathbf{v} ,\phi\right\rangle\right|\lesssim\left\| \Lambda^{n/4} \mathbf{v} \right\|_{L^{2} (\mathbb{R}^{n})}\left\| \Lambda^{n/4} \phi\right\|_{L^{2} (\mathbb{R}^{n})}.\end{equation}
Recall \eqref{1.8} and \eqref{1.10}, combining \eqref{3.2} with \eqref{3.18}, \eqref{3.19} and \eqref{3.20} yields
\begin{equation}\label{3.21}
\left|\left\langle \partial_{t}\mathbf{v} ,\phi\right\rangle\right|\lesssim  \| \mathbf{v} \|_{L^{2} (\mathbb{R}^{n})} \left\| \Lambda^{n/4} \mathbf{v} \right\|_{L^{2} (\mathbb{R}^{n})}\left( \|\phi\| ^{2}_{L^{2}(\mathbb{R}^{n})}+\left\|\Lambda^{n/4}  \phi\right\|^{2}_{L^{2} (\mathbb{R}^{n})}\right)^{\frac{1}{2}}.
\end{equation}
As we can choose arbitrarily for $\phi\in L^{2}\left([0,T],H^{ \frac{n}{4} }_{\sigma}(\mathbb{R}^{n})\right)$, recall \eqref{1.10}, \eqref{1.11} and the assumption of Theorem \ref{t1.1}, we deduce that for $\displaystyle s=\frac{n}{4}$ with $n=2,3$,
\begin{equation}\label{3.22}
\left\| \partial_{t}\mathbf{v} \right\|_{L^{2} \left([0,T],\left[H^{\frac{n}{4}}_{0}(\mathbb{R}^{n})\right]'  \right)}\leq  C\left(n,\alpha,\nu,\left\| \mathbf{v}_{0} \right\|_{L^{2} (\mathbb{R}^{n})}\right).
\end{equation}
Hence \eqref{3.13} together with \eqref{3.22} concludes the estimate \eqref{1.12} for $\displaystyle \frac{n}{4}\leq s<1$ with $n=2,3$.\\
\indent This completes the proof of Theorem \ref{t1.1}.\hfill$\Box$\\

\section{Higher-order regularity}
 In this section we prove Theorem \ref{t1.2}.\\
\indent {\bf Proof of Theorem 1.2.}\\
\indent In Theorem \ref{t1.1} we have established the existence of the weak solutions to the Cauchy problem \eqref{1.1}-\eqref{1.2} in the sense of Definition \ref{d2.3} . In particular, the bound in \eqref{1.13} can be applied to the solution established in Theorem \ref{t1.1}. With this reason, we shall derive formally the high-order regularity properties of the solutions constructed in  Theorem \ref{t1.1}. The proof of Theorem \ref{t1.2} can be devided into two aspects. We first prove the high-order regularity with respect to space. We then verify the high-order regularity with respect to space-time. We shall finish the proof of Theorem \ref{t1.2} after proving the following two theorems.\\ \indent Firstly, we claim the following result concerning the high-order regularity with respect to space.\\
\indent
\begin{thm}[High-order regularity w. s. t. space]\label{t4.1}\end{thm}

 For $n=2,3$, assume that
\\[0.3cm]
 \indent (1) for $\displaystyle\frac{n}{4}< s < 1$, the initial data $\mathbf{v}_{0}\in \mathcal{D}_{\sigma}(\Lambda^{K} )(\mathbb{R}^{n})$,\\
  and\\
 \indent (2) for $\displaystyle s=\frac{n}{4}$, the initial data $\mathbf{v}_{0}\in H_{\sigma}^{K}(\mathbb{R}^{n})$, and in addition, there esists $\varepsilon^{*}=\varepsilon^{*}(\alpha,\nu,n)$ sufficiently small such that $\left\| \mathbf{v}_{0}\right\|_{H_{0}^{K}(\mathbb{R}^{n})}\leq \varepsilon^{*}$.
 \\[0.3cm]
  Then the solutions to the Cauchy problem \eqref{1.1}-\eqref{1.2} constructed in Theorem \ref{t1.1} satisfy that for all $M\leq K$,
\begin{equation}\label{4.1}
\left\|\nabla^{M}\mathbf{v}\right\|^{2}_{L^{2}(\mathbb{R}^{n})}+\nu\int_{0}^{T}
\left\|\nabla^{M}\Lambda^{s} \mathbf{v}\right\|^{2}_{L^{2}(\mathbb{R}^{n})}dt\leq C\left(n,\alpha,\nu,\left\| \mathbf{v}_{0}\right\|_{\mathcal{A}_{1}}\right).
\end{equation}
Here, $\displaystyle \mathcal{A}_{1}=\left\{\begin{array}{ll}
\mathcal{D}\left(\Lambda^{K}\right)(\mathbb{R}^{n})~~&\displaystyle \hbox{for}~~\frac{n}{4}< s < 1
\\[0.3cm]
H_{0}^{K}(\mathbb{R}^{n})~~&\displaystyle \hbox{for}~~s=\frac{n}{4}\end{array}\right.$, \quad $M$ and $K$ are both integers.
\\[0.3cm]
{\bf proof.}\quad We prove it by induction with three steps.
\\[0.3cm]
 {\bf Step 1}\quad We first give the inductive assumption. Let $\displaystyle\frac{n}{4}\leq s < 1$. Assume that for all $m<M$, the following inductive bound holds:
\begin{equation}\label{4.2}
\left.
\begin{array}{ll}
\displaystyle\left\| \nabla^{m}\mathbf{v}\right\|^{2}_{L^{2}(\mathbb{R}^{n})}+\nu\int^{T}_{0}\left\| \nabla^{m}\Lambda^{s} \mathbf{v}\right\|^{2}_{L^{2}(\mathbb{R}^{n})}dt\leq C\left(n,\alpha,\nu,\left\| \mathbf{v}_{0}\right\|_{\mathcal{A}_{1}}\right).
\end{array}
\right.
\end{equation}
\\[0.3cm]
{\bf Step 2}\quad By Theorem \ref{t1.1}, it is easy to verify that the inductive assumption \eqref{4.2} holds for the base case $m=0$.
\\[0.3cm]
 {\bf Step 3}. We will show that the inductive assumption \eqref{4.2} holds for $m=M$.
  \\[0.3cm]
 \indent Multiplying the first equation in \eqref{1.1} by $\Delta^{M}\mathbf{v}$ and integrating in space yields, after some integration by parts,
\begin{equation}\label{4.3}
\left.
\begin{array}{ll}
\dfrac{1}{2}\dfrac{d}{dt}\left\| \nabla^{M}\mathbf{v}\right\|^{2}_{L^{2}(\mathbb{R}^{n})}+\nu\left\| \nabla^{M}\Lambda^{s} \mathbf{v}\right\|^{2}_{L^{2}(  \mathbb{R}^{n})}
\\[0.3cm]
\qquad \leq \left|\left\langle \mathbf{u}\cdot \nabla \mathbf{v}, \Delta^{M}\mathbf{v}\right\rangle\right|+\left|\left\langle \mathbf{v}\cdot \nabla \mathbf{u}^{T}, \Delta^{M}\mathbf{v}\right\rangle\right|
\\[0.3cm]
\qquad  :=I_{M}+J_{M}.
\end{array}
\right.
\end{equation}
We shall estimate $I_{M}$ and $J_{M}$ in \eqref{4.3} through two cases:
\\[0.3cm]
\indent {\bf Case (1)}\quad $\dfrac{n}{4}<s<1$,  $n=2,3$;
\\[0.2cm]
 \indent {\bf Case (2)}\quad $s=\dfrac{n}{4} $,  $n=2,3$.
 \\[0.3cm]
We first consider  {\bf Case (1)}\quad $\dfrac{n}{4}<s<1$,  $n=2,3$.
\\[0.3cm]
 \indent In this case, recall that
 $\left<\mathbf{u}\cdot \nabla \mathbf{v}, \mathbf{v}\right>=0$, thanks to Cauchy's inequality, H\"{o}lder's inequality and the fractional Gagliardo-Nirenberg-Sobolev inequality, one deduces that
\begin{equation}\label{4.4}
\left.
\begin{array}{ll}
\vspace{5pt}
\displaystyle I_{M}=\left|\left\langle \mathbf{u}\cdot \nabla \mathbf{v}, \Delta^{M}\mathbf{v}\right\rangle\right|\\
\vspace{5pt}
\quad\displaystyle=\sum\limits^{M}_{m=1}\left(
\begin{array}{l}M\\
m\end{array}\right)\left\langle\nabla^{m}\mathbf{u}\cdot \nabla\nabla^{M-m} \mathbf{v}, \nabla^{M}\mathbf{v}\right\rangle\\
\vspace{5pt}
\quad\displaystyle\lesssim\sum\limits^{M}_{m=1} \left\|\nabla^{M+1-m} \mathbf{v}\right\|_{L^{2}(\mathbb{R}^{n})}\left\|\nabla^{ m} \mathbf{u}\nabla^{ M} \mathbf{v}\right\|_{L^{2}(\mathbb{R}^{n})}\\
\vspace{5pt}
\quad\displaystyle\lesssim\sum\limits^{M}_{m=1} \left\|\nabla^{M+1-m} \mathbf{v}\right\|_{L^{2}(\mathbb{R}^{n})}\left\|\nabla^{ m} \mathbf{u} \right\|_{L^{\frac{n}{s}}(\mathbb{R}^{n})}\left\| \nabla^{M}\mathbf{v}\right\|_{L^{\frac{2n}{n-2s}}(\mathbb{R}^{n})}\\
\vspace{5pt}
\quad\displaystyle\lesssim\sum\limits^{M}_{m=1} \left\|\nabla^{M+1-m} \mathbf{v}\right\|_{L^{2}(\mathbb{R}^{n})}\left\| \Lambda^{s} \nabla^{ m} \mathbf{u} \right\|_{L^{2}(\mathbb{R}^{n})}\left\|  \Lambda^{s} \nabla^{M}\mathbf{v}\right\|_{L^{2}(\mathbb{R}^{n})}\\
\vspace{5pt}
\quad \displaystyle\lesssim\sum\limits^{M}_{m=1}\left \|\nabla^{M+1-m} \mathbf{v}\right\|^{2}_{L^{2}(\mathbb{R}^{n})}\left\| \Lambda^{s} \nabla^{ m-1} \mathbf{v } \right\|^{2}_{L^{2}(\mathbb{R}^{n})}+\frac{\nu}{4}\left\|  \Lambda^{s} \nabla^{M}\mathbf{v}\right\|^{2}_{L^{2}(\mathbb{R}^{n})}.
\end{array}
\right.
\end{equation}
and
\begin{equation}\label{4.5}
\left.
\begin{array}{ll}
\vspace{5pt}
J_{M}=\left|\left\langle \mathbf{v}\cdot \nabla \mathbf{u}^{T}, \Delta^{M}\mathbf{v}\right\rangle\right|\\
\vspace{5pt}
\quad=\sum\limits^{M}_{m=0}\left(
\begin{array}{l}M\\
m\end{array}\right)\left|\left\langle\nabla^{M}\mathbf{v}\cdot \nabla\nabla^{m} \mathbf{u}, \nabla^{M-m}\mathbf{v}\right\rangle\right|\\
\vspace{5pt}
\quad\lesssim\left|\left\langle\nabla^{M}\mathbf{v}\cdot \nabla \mathbf{u}, \nabla^{M}\mathbf{v}\right\rangle\right|
 +\sum\limits^{M}_{m=1}\left(
\begin{array}{l}M\\
m\end{array}\right)\left|\left\langle\nabla^{M}\mathbf{v}\cdot \nabla^{m+1} \mathbf{u}, \nabla^{M-m}\mathbf{v}\right\rangle\right|.
  \end{array}
\right.
\end{equation}
On the other hand, by Lemma \ref{l2.5} we have
\begin{equation}\label{4.6}
\left.
\begin{array}{ll}
\vspace{5pt}
\left|\left\langle\nabla^{M}\mathbf{v}\cdot \nabla \mathbf{u}, \nabla^{M}\mathbf{v}\right\rangle\right|\\
 \vspace{5pt}
 \qquad \lesssim\left\| \nabla^{M} \mathbf{v}\right\|^{2}_{L^{4}(\mathbb{R}^{n})}
\left\| \nabla \mathbf{u}\right\|_{L^{2}(\mathbb{R}^{n})}\\
\vspace{5pt}
\qquad\lesssim\left(C(\varepsilon)\left\| \nabla^{M} \mathbf{v}\right\|^{2}_{L^{2}(\mathbb{R}^{n})}+\varepsilon\left\| \Lambda^{s} \nabla^{M} \mathbf{v}\right\|^{2}_{L^{2}(\mathbb{R}^{n})}\right)\left\| \nabla \mathbf{u}\right\|_{L^{2}(\mathbb{R}^{n})}.
  \end{array}
\right.
\end{equation}
and
\begin{equation}\label{4.7}
\left.
\begin{array}{ll}
\vspace{5pt}
\sum\limits^{M}_{m=1}\left(
\begin{array}{l}M\\
m\end{array}\right)\left|\left\langle\nabla^{M}\mathbf{v}\cdot \nabla^{m+1} \mathbf{u}, \nabla^{M-m}\mathbf{v}\right\rangle\right|\\
\vspace{5pt}
\qquad\lesssim\sum\limits^{M}_{m=1} \left\| \nabla^{M-m} \mathbf{v}\right\|_{L^{2}(\mathbb{R}^{n})}\left\| \nabla^{m+1} \mathbf{u}\cdot\nabla^{M}\mathbf{v} \right\|_{L^{2}(\mathbb{R}^{n})}\\
\vspace{5pt}
\qquad\lesssim\sum\limits^{M}_{m=1} \left\| \nabla^{M-m} \mathbf{v}\right\|_{L^{2}(\mathbb{R}^{n})}\left\| \nabla^{m+1} \mathbf{u } \right\|_{L^{\frac{n}{s}}(\mathbb{R}^{n})}\left\| \nabla^{M} \mathbf{v} \right\|_{L^{\frac{2n}{n-2s}}(\mathbb{R}^{n})}\\
\vspace{5pt}
\qquad\lesssim\sum\limits^{M}_{m=1} \left\| \nabla^{M-m} \mathbf{v}\right\|_{L^{2}(\mathbb{R}^{n})}\left\| \Lambda^{s} \nabla^{m-1} \mathbf{v} \right\|_{L^{2}(\mathbb{R}^{n})}\left\| \Lambda^{s} \nabla^{M} \mathbf{v} \right\|_{L^{2}(\mathbb{R}^{n})}\\
\vspace{5pt}
\qquad\lesssim\displaystyle\sum\limits^{M}_{m=1} \left\| \nabla^{M-m} \mathbf{v}\right\|^{2}_{L^{2}(\mathbb{R}^{n})}\left\| \Lambda^{s} \nabla^{m-1} \mathbf{v} \right\|^{2}_{L^{2}(\mathbb{R}^{n})}+\frac{\nu}{4}\left\| \Lambda^{s} \nabla^{M} \mathbf{v} \right\|^{2}_{L^{2}(\mathbb{R}^{n})}.\end{array}
\right.
\end{equation}
Note that \eqref{2.3}, under the assumptions of Theorem \ref{t4.1}, choosing $\varepsilon\leq \dfrac{\nu}{2 \left\| \nabla \mathbf{u}_{0} \right\|^{2}_{L^{2}(\mathbb{R}^{n})}}$, by \eqref{4.3}, \eqref{4.4}, \eqref{4.5}, \eqref{4.6} and \eqref{4.7}, one has
\begin{equation}\label{4.8}
\left.
\begin{array}{ll}
\vspace{5pt}
\dfrac{d}{dt}\left\| \nabla^{M}\mathbf{v}\right\|^{2}_{L^{2}(\mathbb{R}^{n})}\lesssim\sum\limits^{M}_{m=1}  \left\|\nabla^{M+1-m} \mathbf{v}\right\|^{2}_{L^{2}(\mathbb{R}^{n})}\left\| \Lambda^{s} \nabla^{ m-1} \mathbf{v} \right\|^{2}_{L^{2}(\mathbb{R}^{n})}\\
\vspace{5pt}
\qquad \qquad\qquad\quad+\sum\limits^{M}_{m=1} \left\| \nabla^{M-m} \mathbf{v}\right\|^{2}_{L^{2}(\mathbb{R}^{n})}\left\| \Lambda^{s} \nabla^{m-1} \mathbf{v } \right\|^{2}_{L^{2}(\mathbb{R}^{n})}\\
\vspace{5pt}
\qquad \qquad\qquad\quad+ \left\| \mathbf{v}_{0}\right\|_{L^{2}(\mathbb{R}^{n})}\left\| \nabla^{M} \mathbf{v}\right\|^{2}_{L^{2}(\mathbb{R}^{n})}.
\end{array}
\right.
\end{equation}
This together with Theorem \ref{t1.1}, \eqref{4.2} and Gronwall's inequality leads to \eqref{4.1} under the case $\dfrac{n}{4}<s<1$ with $n=2,3$.
\\[0.3cm]
We then consider {\bf Case (2)}\quad $s=\dfrac{n}{4}$ for $n=2,3$.
\\[0.3cm]
\indent In this case, recall again that
 $\left\langle \mathbf{u}\cdot \nabla \mathbf{v}, \mathbf{v}\right\rangle=0$,  Cauchy's inequality,  H\"{o}lder's inequality and the Gagliardo-Nirenberg-Sobolev inequality give rise to the following two estimates:
\begin{equation}\label{4.9}
\left.
\begin{array}{ll}
\vspace{5pt}
I_{M}=\left|\left\langle \mathbf{u}\cdot \nabla \mathbf{v}, \Delta^{M}\mathbf{v}\right\rangle\right|\\
\vspace{5pt}
\quad=\sum\limits^{M}_{m=1}\left(
\begin{array}{l}M\\
m\end{array}\right)\left\langle\nabla^{m}\mathbf{u}\cdot \nabla\nabla^{M-m} \mathbf{v}, \nabla^{M}\mathbf{v}\right\rangle\\
\vspace{5pt}
\quad\lesssim\sum\limits^{M}_{m=1} \left\|\nabla^{M+1-m} \mathbf{v}\right\|_{L^{2}(\mathbb{R}^{n})}\left\|\nabla^{ m} u\nabla^{ M} \mathbf{v}\right\|_{L^{2}(\mathbb{R}^{n})}\\
\vspace{5pt}
\quad\lesssim\sum\limits^{M}_{m=1} \left\|\nabla^{M+1-m} \mathbf{v}\right\|_{L^{2}(\mathbb{R}^{n})}\left\|\nabla^{ m} \mathbf{u} \right\|_{L^{4}(\mathbb{R}^{n} )}\left\| \nabla^{M}\mathbf{v}\right\|_{L^{4}(\mathbb{R}^{n})}\\
\vspace{5pt}
\quad\lesssim\sum\limits^{M}_{m=1}\left \|\nabla^{M+1-m} \mathbf{v}\right\|_{L^{2}(\mathbb{R}^{n})} \left(\left\| \nabla^{ m} \mathbf{u} \right\|^{2}_{L^{2}(\mathbb{R}^{n})}+\left\|\Lambda^{\frac{n}{4}} \nabla^{ m} \mathbf{u} \right\|^{2}_{L^{2}(\mathbb{R}^{n})}\right)^{\frac{1}{2}}
\\[0.3cm]
\qquad\qquad\qquad\qquad\times\left(\left\| \nabla^{M}\mathbf{v}\right\|^{2}_{L^{2}(\mathbb{R}^{n})}+
\left\| \Lambda^{\frac{n}{4}} \nabla^{M}\mathbf{v}\right\|^{2}_{L^{2}(\mathbb{R}^{n})}\right)^{\frac{1}{2}}\\
\vspace{5pt}
\quad\lesssim\sum\limits^{M}_{m=1}\left \|\nabla^{M+1-m} \mathbf{v}\right\|_{L^{2}(\mathbb{R}^{n})} \left(\left\| \nabla^{ m-1} \mathbf{v} \right\|^{2}_{L^{2}(\mathbb{R}^{n})}+\left\|\Lambda^{\frac{n}{4}} \nabla^{ m-1} \mathbf{v} \right\|^{2}_{L^{2}(\mathbb{R}^{n})}\right)^{\frac{1}{2}}
\\[0.3cm]
\qquad\qquad\qquad\qquad\times
\left(\left\| \nabla^{M}\mathbf{v}\right\|^{2}_{L^{2}(\mathbb{R}^{n})}+
\left\| \Lambda^{\frac{n}{4}} \nabla^{M}\mathbf{v}\right\|^{2}_{L^{2}(\mathbb{R}^{n})}\right)^{\frac{1}{2}},\end{array}
\right.
\end{equation}
and
\begin{equation}\label{4.10}
\left.
\begin{array}{ll}
\vspace{5pt}
J_{M}=\left|\left\langle \mathbf{v}\cdot \nabla \mathbf{u}^{T}, \Delta^{M}\mathbf{v}\right\rangle\right|\\
\vspace{5pt}
\quad=\sum\limits^{M}_{m=0}\left(
\begin{array}{l}M\\
m\end{array}\right)\left|\left\langle\nabla^{M}\mathbf{v}\cdot  \nabla^{m+1} \mathbf{u}, \nabla^{M-m}\mathbf{v}\right\rangle\right|\\
\vspace{5pt}
\quad\lesssim\sum\limits^{M}_{m=0}\left(
\begin{array}{l}M\\
m\end{array}\right)\left\| \nabla^{M-m} \mathbf{v}\right\|_{L^{2}(\mathbb{R}^{n})}\left\| \nabla^{m+1} \mathbf{u}\cdot\nabla^{M}\mathbf{v}\right\|_{L^{2}(\mathbb{R}^{n})}\\
\vspace{5pt}
\quad\lesssim\sum\limits^{M}_{m=0} \left\| \nabla^{M-m} \mathbf{v}\right\|_{L^{2}(\mathbb{R}^{n})}\left\| \nabla^{m+1} \mathbf{u} \right\|_{L^{4}(\mathbb{R}^{n})}\left\| \nabla^{M}\mathbf{v}\right\|_{L^{4}(\mathbb{R}^{n})}\\
\vspace{5pt}
\quad\lesssim\sum\limits^{M}_{m=0}\left \|\nabla^{M-m} \mathbf{v}\right\|_{L^{2}(\mathbb{R}^{n})} \left(\left\| \nabla^{ m+1} \mathbf{u} \right\|^{2}_{L^{2}(\mathbb{R}^{n})}+\left\|\Lambda^{\frac{n}{4}} \nabla^{ m+1} \mathbf{u} \right\|^{2}_{L^{2}(\mathbb{R}^{n})}\right)^{\frac{1}{2}}
\\[0.2cm]
\qquad\qquad\qquad\qquad\times
\left(\left\| \nabla^{M}\mathbf{v}\right\|^{2}_{L^{2}(\mathbb{R}^{n})}+
\left\| \Lambda^{\frac{n}{4}} \nabla^{M}\mathbf{v}\right\|^{2}_{L^{2}(\mathbb{R}^{n})}\right)^{\frac{1}{2}}\\
\vspace{5pt}
\quad\lesssim\sum\limits^{M}_{m=1}\left \|\nabla^{M+1-m} \mathbf{v}\right\|_{L^{2}(\mathbb{R}^{n})} \left(\left\| \nabla^{ m-1} \mathbf{v} \right\|^{2}_{L^{2}(\mathbb{R}^{n})}+\left\|\Lambda^{\frac{n}{4}} \nabla^{ m-1} \mathbf{v} \right\|^{2}_{L^{2}(\mathbb{R}^{n})}\right)^{\frac{1}{2}}
\\[0.2cm]
\qquad\qquad\qquad\qquad\times
\left(\left\| \nabla^{M}\mathbf{v}\right\|^{2}_{L^{2}(\mathbb{R}^{n})}+
\left\| \Lambda^{\frac{n}{4}} \nabla^{M}\mathbf{v}\right\|^{2}_{L^{2}(\mathbb{R}^{n})}\right)^{\frac{1}{2}}.
\end{array}
\right.
\end{equation}
 Combining \eqref{4.9} with \eqref{4.10} yields
 \begin{equation}\label{4.11}
\left.
\begin{array}{ll}
I_{M}+J_{M}&\lesssim\sum\limits^{M}_{m=1}\left \|\nabla^{M+1-m} \mathbf{v}\right\|_{L^{2}(\mathbb{R}^{n})}\left(\left\| \nabla^{ m-1} \mathbf{v} \right\|^{2}_{L^{2}(\mathbb{R}^{n})}+\left\|\Lambda^{\frac{n}{4}} \nabla^{ m-1} \mathbf{v} \right\|^{2}_{L^{2}(\mathbb{R}^{n})}\right)^{\frac{1}{2}}
\\[0.2cm]
 &\qquad\qquad\qquad\qquad\times\left(\left\| \nabla^{ M} \mathbf{v} \right\|^{2}_{L^{2}(\mathbb{R}^{n})}+\left\|\Lambda^{\frac{n}{4}} \nabla^{ M} \mathbf{v} \right\|^{2}_{L^{2}(\mathbb{R}^{n})}\right)^{\frac{1}{2}}
 \\[0.2cm]
 &\lesssim\left\| \nabla^{ M} \mathbf{v} \right\|_{L^{2}(\mathbb{R}^{n})}\left(\left\| \mathbf{v} \right\|^{2}_{L^{2}(\mathbb{R}^{n})}+\left\|\Lambda^{\frac{n}{4}} \mathbf{v} \right\|^{2}_{L^{2}(\mathbb{R}^{n})}\right)^{\frac{1}{2}}
 \\[0.2cm]
 &\qquad\qquad\qquad\times\left(\left\| \nabla^{ M} \mathbf{v} \right\|^{2}_{L^{2}(\mathbb{R}^{n})}+\left\|\Lambda^{\frac{n}{4}} \nabla^{ M} \mathbf{v} \right\|^{2}_{L^{2}(\mathbb{R}^{n})}\right)^{\frac{1}{2}}
 \\[0.2cm]
 &\quad+\sum\limits^{M-1}_{m=2}\left \|\nabla^{M+1-m} \mathbf{v}\right\|_{L^{2}(\mathbb{R}^{n})}\left(\left\| \nabla^{ m-1} \mathbf{v} \right\|^{2}_{L^{2}(\mathbb{R}^{n})}+\left\|\Lambda^{\frac{n}{4}} \nabla^{ m-1} \mathbf{v} \right\|^{2}_{L^{2}(\mathbb{R}^{n})}\right)^{\frac{1}{2}}
\\[0.2cm]
 &\qquad\qquad\qquad\qquad\times\left(\left\| \nabla^{ M} \mathbf{v} \right\|^{2}_{L^{2}(\mathbb{R}^{n})}+\left\|\Lambda^{\frac{n}{4}} \nabla^{ M} \mathbf{v} \right\|^{2}_{L^{2}(\mathbb{R}^{n})}\right)^{\frac{1}{2}}
 \\[0.3cm]
  &\quad+ \left \|\nabla \mathbf{v}\right\|_{L^{2}(\mathbb{R}^{n})}\left(\left\| \nabla^{ M-1} \mathbf{v} \right\|^{2}_{L^{2}(\mathbb{R}^{n})}+\left\|\Lambda^{\frac{n}{4}} \nabla^{ M-1} \mathbf{v} \right\|^{2}_{L^{2}(\mathbb{R}^{n})}\right)^{\frac{1}{2}}
\\[0.3cm]
 &\qquad\qquad\qquad \times\left(\left\| \nabla^{ M} \mathbf{v} \right\|^{2}_{L^{2}(\mathbb{R}^{n})}+\left\|\Lambda^{\frac{n}{4}} \nabla^{ M} \mathbf{v} \right\|^{2}_{L^{2}(\mathbb{R}^{n})}\right)^{\frac{1}{2}}
 \\[0.3cm]
 &\lesssim \left\| \nabla^{ M} \mathbf{v} \right\|^{2}_{L^{2}(\mathbb{R}^{n})}+\left(\left\| \mathbf{v} \right\|^{2}_{L^{2}(\mathbb{R}^{n})}+\left\|\Lambda^{\frac{n}{4}} \mathbf{v} \right\|^{2}_{L^{2}(\mathbb{R}^{n})}\right)
  \\[0.3cm]
 &\qquad\qquad\qquad\times\left(\left\| \nabla^{ M} \mathbf{v} \right\|^{2}_{L^{2}(\mathbb{R}^{n})}+\left\|\Lambda^{\frac{n}{4}} \nabla^{ M} \mathbf{v} \right\|^{2}_{L^{2}(\mathbb{R}^{n})}\right)
 \\[0.3cm]
 & \quad+\sum\limits^{M-1}_{m=2}\left(\left \|\nabla^{M+1-m} \mathbf{v}\right\|^{2}_{L^{2}(\mathbb{R}^{n})}+\left(\left\| \nabla^{ m-1} \mathbf{v} \right\|^{2}_{L^{2}(\mathbb{R}^{n})}+\left\|\Lambda^{\frac{n}{4}} \nabla^{ m-1} \mathbf{v} \right\|^{2}_{L^{2}(\mathbb{R}^{n})}\right)\right.
\\[0.3cm]
 &\qquad\qquad\qquad\qquad\qquad\times\left.\left(\left\| \nabla^{ M} \mathbf{v} \right\|^{2}_{L^{2}(\mathbb{R}^{n})}+\left\|\Lambda^{\frac{n}{4}} \nabla^{ M} \mathbf{v} \right\|^{2}_{L^{2}(\mathbb{R}^{n})}\right)\right)
 \\[0.3cm]
 & \quad+ \left \|\nabla \mathbf{v}\right\|^{2}_{L^{2}(\mathbb{R}^{n})} \left(\left\| \nabla^{ M-1} \mathbf{v} \right\|^{2}_{L^{2}(\mathbb{R}^{n})}+\left\|\Lambda^{\frac{n}{4}} \nabla^{ M-1} \mathbf{v} \right\|^{2}_{L^{2}(\mathbb{R}^{n})}\right)
\\[0.3cm]
 &\qquad\qquad\qquad\times \left(\left\| \nabla^{ M} \mathbf{v} \right\|^{2}_{L^{2}(\mathbb{R}^{n})}+\left\|\Lambda^{\frac{n}{4}} \nabla^{ M} \mathbf{v} \right\|^{2}_{L^{2}(\mathbb{R}^{n})}\right).
  \end{array}
\right.
\end{equation}
Using interpolation inequality, we obtain for $M\leq K$,
 \begin{equation*}
\left.
\begin{array}{ll}
\vspace{5pt}
 \left(\left\| \mathbf{v}\right\|^{2}_{L^{2}(\mathbb{R}^{n})} +
\left\| \Lambda^{\frac{n}{4}} \mathbf{v}\right\|^{2}_{L^{2}(\mathbb{R}^{n})}\right)
   +\sum\limits^{M-1}_{m=2} \left(\left\| \nabla^{ m-1} \mathbf{v} \right\|^{2}_{L^{2}(\mathbb{R}^{n})}+\left\|\Lambda^{\frac{n}{4}} \nabla^{ m-1} \mathbf{v} \right\|^{2}_{L^{2}(\mathbb{R}^{n})}\right)
  \\[0.2cm]
 \qquad+ \left(\left\| \nabla^{ M-1} \mathbf{v} \right\|^{2}_{L^{2}(\mathbb{R}^{n})}+\left\|\Lambda^{\frac{n}{4}} \nabla^{ M-1} \mathbf{v} \right\|^{2}_{L^{2}(\mathbb{R}^{n})}\right)
\\[0.2cm]
\qquad\quad=\sum\limits^{M}_{m=1} \left(\left\| \nabla^{ m-1} \mathbf{v} \right\|^{2}_{L^{2}(\mathbb{R}^{n})}+\left\|\Lambda^{\frac{n}{4}} \nabla^{ m-1} \mathbf{v} \right\|^{2}_{L^{2}(\mathbb{R}^{n})}\right)
\\[0.2cm]
\qquad\quad \lesssim \left\|   \mathbf{v} \right\|^{2}_{L^{2}(\mathbb{R}^{n})}+\left\|\Lambda^{\frac{n}{4}} \nabla^{ M-1} \mathbf{v} \right\|^{2}_{L^{2}(\mathbb{R}^{n})}
\\[0.2cm]
\qquad\quad\lesssim \left\|   \mathbf{v} \right\|^{2}_{H_{0}^{\frac{n}{4}+M-1}(\mathbb{R}^{n})}\lesssim\left\|   \mathbf{v} _{0} \right\|^{2}_{H_{0}^{K}(\mathbb{R}^{n})}.
 \end{array}
\right.
\end{equation*}
Due to the assumption of this theorem, choosing $\varepsilon^{*}$ sufficiently small such that
$$\left\|   \mathbf{v} _{0} \right\|^{2}_{H_{0}^{K}(\mathbb{R}^{n})}\leq \varepsilon^{*}<\frac{\nu}{2},$$
which together with the inductive assumption \eqref{4.2}, \eqref{1.11}, interpolation inequality and Gronwall's inequality leads to \eqref{4.1} under the case  $\displaystyle s=\frac{n}{4}$ with $n=2,3$. \\
  \indent This finishes the proof of Theorem \ref{t4.1}.\hfill$\Box$\\
  \begin{thm}[High-order regularity w. r. t. space-time]\label{t4.2}\rm
   For $n=2,3$, assume that
 \\[0.3cm]
 \indent (1)\quad for $\displaystyle\frac{n}{4}< s < 1$, the initial data $\mathbf{v}_{0}\in \mathcal{D}_{\sigma}\left(\Lambda^{K} \right)(\mathbb{R}^{n})$,\\
  and\\
   \indent (2) for $\displaystyle s=\frac{n}{4}$, the initial data $\mathbf{v}_{0}\in H^{K}_{\sigma} (\mathbb{R}^{n})$, and in addition, there exists an $\varepsilon^{**}=\varepsilon^{**}(\alpha,\nu,n)>0$ such that $\left\| \mathbf{v}_{0}\right\|_{H^{K}_{0} (\mathbb{R}^{n})}\leq \varepsilon^{**}$.
   \\[0.3cm]
   Then for all $M+2Ps\leq K$, the solutions to the Cauchy problem \eqref{1.1}-\eqref{1.2} constructed in Theorem \ref{t1.1} admits the bound
\begin{equation}\label{4.12}
\left\|\partial^{P}_{t}\nabla^{M}\mathbf{v}\right\|^{2}_{L^{2}(\mathbb{R}^{n})}
+\nu\int_{0}^{T}
\left\|\partial^{P}_{t}\nabla^{M}\Lambda^{s} \mathbf{v}\right\|^{2}_{L^{2}(\mathbb{R}^{n})}dt\leq C\left(n,\alpha,\nu,\| \mathbf{v}_{0}\|_{\mathcal{A}_{2}}\right).
\end{equation}
Here, $\displaystyle \mathcal{A}_{2}=\left\{\begin{array}{ll}
\mathcal{D}\left(\Lambda^{K}\right)(\mathbb{R}^{n})~~&\displaystyle\hbox{for}~~\frac{n}{4}< s < 1
\\[0.3cm]
H_{0}^{K}(\mathbb{R}^{n})~~&\displaystyle\hbox{for}~~s=\frac{n}{4} \end{array}\right.$, ~~$M$, $K$ and $P$ are all integers.

\end{thm}
  {\bf Proof.}\quad Applying $\partial^{P}_{t}\nabla^{M}$ to the solutions of \eqref{1.1}-\eqref{1.2}, we have
  \begin{equation}\label{4.13}
 \partial^{P+1}_{t}\nabla^{M}\mathbf{v}+\partial^{P}_{t}\nabla^{M}\left(\mathbf{u}\cdot \nabla \mathbf{v}\right)+\partial^{P}_{t}\nabla^{M}\left(\mathbf{v}\cdot \nabla \mathbf{u}^{T}\right)+\partial^{P}_{t}\nabla^{M}\nabla p=-\nu\partial^{P}_{t}\nabla^{M}(-\Delta)^{s}\mathbf{v}.
  \end{equation}
 Direct calculation gives
  \begin{equation*}
 \partial^{P+1}_{t}\nabla^{M}\mathbf{v}+\partial^{P}_{t}\nabla^{M}\nabla p=-\partial^{P}_{t}\nabla^{M}\left(\mathbf{u}\cdot \nabla \mathbf{v}\right)-\partial^{P}_{t}\nabla^{M}\left(\mathbf{v}\cdot \nabla \mathbf{u}^{T}\right)-\nu\partial^{P}_{t}\nabla^{M}(-\Delta)^{s}\mathbf{v},
  \end{equation*}
 this yields that
 \begin{equation*}
 \left.
\begin{array}{ll}
\left\| \partial^{P+1}_{t}\nabla^{M}\mathbf{v}+\partial^{P}_{t}\nabla^{M}\nabla p\right\|^{2}_{L^{2}(\mathbb{R}^{n})}
\\[0.3cm]
\qquad\qquad\lesssim\left\| -\partial^{P}_{t}\nabla^{M}\left(\mathbf{u}\cdot \nabla \mathbf{v}\right)-\partial^{P}_{t}\nabla^{M}\left(\mathbf{v}\cdot \nabla \mathbf{u}^{T}\right)
-\nu\partial^{P}_{t}\nabla^{M}(-\Delta)^{s}\mathbf{v}\right\|^{2}_{L^{2}(\mathbb{R}^{n})}
\\[0.3cm]
\qquad\qquad\lesssim  \left\| \partial^{P}_{t} \nabla^{M}\Lambda^{2s} \mathbf{v}\right\|^{2}_{L^{2}(\mathbb{R}^{n})}+
\left\| \partial^{P}_{t} \nabla^{M}\left(\mathbf{u}\cdot \nabla \mathbf{v}\right)\right\|^{2}_{L^{2}(\mathbb{R}^{n})}+ \left\| \partial^{P}_{t} \nabla^{M}\left(\mathbf{v}\cdot \nabla \mathbf{u}^{T}\right)\right\|^{2}_{L^{2}(\mathbb{R}^{n})}.
 \end{array}
\right.
\end{equation*}
Note that
 \begin{equation*}
 \left.
\begin{array}{ll}
\displaystyle\int_{\mathbb{R}^{n}}\left | \partial^{P+1}_{t}\nabla^{M}\mathbf{v}+\partial^{P}_{t}\nabla^{M}\nabla p\right |^{2} dx
\\[0.3cm]
\qquad\qquad=\displaystyle\left\|\partial^{P+1}_{t}\nabla^{M}\mathbf{v}\right\|^{2}_{L^{2}(\mathbb{R}^{n})}
+\left\|\partial^{P}_{t}\nabla^{M}\nabla p\right\|^{2}_{L^{2}(\mathbb{R}^{n})}+
2\int_{\mathbb{R}^{n}} \partial^{P+1}_{t}\nabla^{M}\mathbf{v}\cdot \partial^{P}_{t}\nabla^{M}\nabla p dx,
 \end{array}
\right.
\end{equation*}
 and
 $$\int_{\mathbb{R}^{n}} \partial^{P+1}_{t}\nabla^{M}\mathbf{v}\cdot \partial^{P}_{t}\nabla^{M}\nabla p dx=-\int_{\mathbb{R}^{n}} \partial^{P+1}_{t}\nabla^{M}\nabla\cdot\mathbf{v}\cdot \partial^{P}_{t}\nabla^{M} p dx=0$$
 from $\mbox{div}~\mathbf{v}=0$ , we get
\begin{equation}\label{4.14}
\left\| \partial^{P+1}_{t}\nabla^{M}\mathbf{v}\right\|^{2}_{L^{2}(\mathbb{R}^{n})}\lesssim  \left\| \partial^{P}_{t} \nabla^{M}\Lambda^{2s} \mathbf{v}\right\|^{2}_{L^{2}(\mathbb{R}^{n})}+
\left\| \partial^{P}_{t} \nabla^{M}\left(\mathbf{u}\cdot \nabla \mathbf{v}\right)\right\|^{2}_{L^{2}(\mathbb{R}^{n})}+ \left\| \partial^{P}_{t} \nabla^{M}\left(\mathbf{v}\cdot \nabla \mathbf{u}^{T}\right)\right\|^{2}_{L^{2}(\mathbb{R}^{n})}.
\end{equation}
We shall estimate the last two terms on the right hand side of \eqref{4.14} through considering two cases:
\\[0.3cm]
\indent {\bf Case (1)} \quad $\displaystyle\frac{n}{4}<s<1$,  $n=2,3$;
\\[0.15cm]
 \indent {\bf Case (2)} \quad $\displaystyle s=\frac{n}{4} $,  $n=2,3$.
 \\[0.3cm]
 We first consider {\bf Case (1)}\quad $\displaystyle\frac{n}{4}<s<1$,  $n=2,3$.\\
  \indent In this case, direct calculation gives $\displaystyle 0< \frac{n}{2}-2s+1<1 $. Thanks to H\"{o}lder's inequality and the Gagliardo-Nirenberg-Sobolev inequality, one obtains
\begin{equation}\label{4.15}
\left.
\begin{array}{ll}
\vspace{5pt}
\left\| \partial^{P}_{t} \nabla^{M}\left(\mathbf{u}\cdot \nabla \mathbf{v}\right)\right\|^{2}_{L^{2}(\mathbb{R}^{n})}\\
\vspace{5pt}
\qquad \lesssim \sum\limits^{P}_{p=0}\sum\limits^{M}_{m=0}\left(
\begin{array}{l}P\\
p\end{array}\right)\left(
\begin{array}{l}M\\
m\end{array}\right) \left\| \partial^{p}_{t}\nabla^{m}\mathbf{u}\cdot\partial^{P-p}_{t}\nabla^{M-m+1}\mathbf{v}\right\|^{2}_{L^{2}(\mathbb{R}^{n})}\\
\vspace{5pt}
\qquad\lesssim \sum\limits^{P}_{p=0}\sum\limits^{M}_{m=0}\left(
\begin{array}{l}P\\
p\end{array}\right)\left(
\begin{array}{l}M\\
m\end{array}\right) \left\| \partial^{p}_{t}\nabla^{m}\mathbf{u}\right\|^{2} _{L^{\frac{n}{2s-1}}(\mathbb{R}^{n})}
\left\|\partial^{P-p}_{t}\nabla^{M-m+1}\mathbf{v}\right\|^{2}_{L^{\frac{2n}{n-2(2s-1)}}(\mathbb{R}^{n})}\\
\vspace{5pt}
\qquad\lesssim\sum\limits^{P}_{p=0}\sum\limits^{M}_{m=0}
\left  \| \partial^{p}_{t}\nabla^{m}\Lambda^{\frac{n}{2}-2s+1} \mathbf{u}\right\|^{2}_{L^{2}(\mathbb{R}^{n})}
\left\| \partial^{P-p}_{t}\nabla^{M-m }\Lambda^{2s} \mathbf{v}\right\|^{2}_{L^{2}(\mathbb{R}^{n})}\\
\vspace{5pt}
\qquad\lesssim\sum\limits^{P}_{p=0}\sum\limits^{M}_{m=0}
 \left\| \partial^{p}_{t}\nabla^{m}\mathbf{v}\right\|^{2}_{L^{2}(\mathbb{R}^{n})}
\left\| \partial^{P-p}_{t}\nabla^{M-m }\Lambda^{2s} \mathbf{v}\right\|^{2}_{L^{2}(\mathbb{R}^{n})}.
 \end{array}
\right.
\end{equation}
In the same manner, a straightforward computation shows that
\begin{equation}\label{4.16}
\left.
\begin{array}{ll}
\vspace{5pt}
\left\| \partial^{P}_{t} \nabla^{M}\left(\mathbf{v}\cdot \nabla \mathbf{u}^{T}\right)\right\|^{2}_{L^{2}(\mathbb{R}^{n})}\\
\vspace{5pt}
 \qquad\lesssim \sum\limits^{P}_{p=0}\sum\limits^{M}_{m=0}\left(
\begin{array}{l}P\\
p\end{array}\right)\left(
\begin{array}{l}M\\
m\end{array}\right)\left \| \partial^{p}_{t}\nabla^{m+1}\mathbf{u}\right\|^{2}_{L^{\frac{n}{s}}(\mathbb{R}^{n})}
\left\| \partial^{P-p}_{t}\nabla^{M-m}\mathbf{v}\right\|^{2}_{L^{\frac{2n}{n-2s}}(\mathbb{R}^{n})}\\
\vspace{5pt}
\qquad\lesssim\sum\limits^{P}_{p=0}\sum\limits^{M}_{m=0}\left(
\begin{array}{l}P\\
p\end{array}\right)\left(
\begin{array}{l}M\\
m\end{array}\right)
\left \| \partial^{p}_{t}\nabla^{m+1}\Lambda^{s} \mathbf{u}\right\|^{2}_{L^{2}(\mathbb{R}^{n})}
\left\| \partial^{P-p}_{t}\nabla^{M-m}\Lambda^{s} \mathbf{v}\right\|^{2}_{L^{2}(\mathbb{R}^{n})}\\
\vspace{5pt}
\qquad
\lesssim\sum\limits^{P}_{p=0}\sum\limits^{M}_{m=0}\left(
\begin{array}{l}P\\
p\end{array}\right)\left(
\begin{array}{l}M\\
m\end{array}\right)
 \left\| \partial^{p}_{t}\nabla^{m}\Lambda^{s} \mathbf{v}\right\|^{2}_{L^{2}(\mathbb{R}^{n})}
\left\| \partial^{P-p}_{t}\nabla^{M-m}\Lambda^{s} \mathbf{v}\right\|^{2}_{L^{2}(\mathbb{R}^{n})}.
 \end{array}
\right.
\end{equation}
By the aid of interpolation inequality, recall that \eqref{1.11},
combining \eqref{4.14} with \eqref{4.15}-\eqref{4.16} gives rise to
$$\left\| \partial^{P+1}_{t}\nabla^{M}\mathbf{v}\right\|^{2}_{L^{2}(\mathbb{R}^{n})}\lesssim\left\| \partial^{P}_{t}\mathbf{v}\right\|^{2}_{\mathcal{D}\left(\Lambda^{M+2s}\right)(\mathbb{R}^{n})}.$$
This implies, for all $M,P$ such that $M+2Ps\leq K$,
 \begin{equation}\label{4.17}
\left \| \partial^{P}_{t}\nabla^{M}\mathbf{v}\right\|^{2}_{L^{2}(\mathbb{R}^{n})}\lesssim\left\| \mathbf{v}\right\|^{2}_{\mathcal{D}\left(\Lambda^{K}\right )(\mathbb{R}^{n})}.
 \end{equation}
 \\[0.2cm]
 We next deal with {\bf Case (2)}\quad $\displaystyle s=\frac{n}{4} $,  $n=2,3$.\\
   \indent In this case, using Lemma \ref{l2.2} and Lemma \ref{l2.5}, we obtain
\begin{equation}\label{4.18}
\left.
\begin{array}{ll}
\left\| \partial^{P}_{t} \nabla^{M}\left(\mathbf{u}\cdot \nabla \mathbf{v}\right)\right\|^{2}_{L^{2}(\mathbb{R}^{n})}
\\[0.3cm]
\qquad \lesssim \sum\limits^{P}_{p=0}\sum\limits^{M}_{m=0}\left(
\begin{array}{l}P\\
p\end{array}\right)\left(
\begin{array}{l}M\\
m\end{array}\right) \left\| \partial^{p}_{t}\nabla^{m}\mathbf{u}\cdot\partial^{P-p}_{t}\nabla^{M-m+1}\mathbf{v}\right\|^{2}_{L^{2}(\mathbb{R}^{n})}
 \\[0.4cm]
\qquad\lesssim\sum\limits^{P}_{p=0}\sum\limits^{M}_{m=0}\left(
\begin{array}{l}P\\
p\end{array}\right)\left(
\begin{array}{l}M\\
m\end{array}\right) \left\| \partial^{p}_{t}\nabla^{m}\mathbf{u}\right\|^{2} _{L^{4}(\mathbb{R}^{n})} \left\|\partial^{P-p}_{t}\nabla^{M-m+1}\mathbf{v}\right\|^{2}_{L^{4}(\mathbb{R}^{n})}
\\[0.4cm]
\qquad\lesssim\sum\limits^{P}_{p=0}\sum\limits^{M}_{m=0}
\left(\left \| \partial^{p}_{t}\nabla^{m}  \mathbf{u}\right\|^{2}_{L^{2}(\mathbb{R}^{n})}+\left \| \partial^{p}_{t}\nabla^{m}\Lambda^{n/4}  \mathbf{u}\right\|^{2}_{L^{2}(\mathbb{R}^{n})}\right)
\\[0.3cm]
\qquad\qquad\qquad\times
\left(\left\| \partial^{P-p}_{t}\nabla^{M+1-m } \mathbf{v}\right\|^{2}_{L^{2}(\mathbb{R}^{n})}+\left\| \partial^{P-p}_{t}\nabla^{M+1-m }\Lambda^{n/4} \mathbf{v}\right\|^{2}_{L^{2}(\mathbb{R}^{n})}\right)
\\[0.3cm]
\qquad\lesssim\sum\limits^{P}_{p=0}\sum\limits^{M}_{m=0}
\left(\left \| \partial^{p}_{t}\nabla^{m}  \mathbf{v}\right\|^{2}_{L^{2}(\mathbb{R}^{n})}+\left \| \partial^{p}_{t}\nabla^{m}\Lambda^{n/4}  \mathbf{v}\right\|^{2}_{L^{2}(\mathbb{R}^{n})}\right)
\\[0.3cm]
\qquad\qquad\qquad\times
\left(\left\| \partial^{P-p}_{t}\nabla^{M+1-m } \mathbf{v}\right\|^{2}_{L^{2}(\mathbb{R}^{n})}+\left\| \partial^{P-p}_{t}\nabla^{M+1-m }\Lambda^{n/4} \mathbf{v}\right\|^{2}_{L^{2}(\mathbb{R}^{n})}\right).
\end{array}
\right.
\end{equation}
The same argument leads to
\begin{equation}\label{4.19}
\left.
\begin{array}{ll}
\left\| \partial^{P}_{t} \nabla^{M}\left(\mathbf{v}\cdot \nabla \mathbf{u}^{T}\right)\right\|^{2}_{L^{2}(\mathbb{R}^{n})}
\\[0.3cm]
\qquad \lesssim \sum\limits^{P}_{p=0}\sum\limits^{M}_{m=0}\left(
\begin{array}{l}P\\
p\end{array}\right)\left(
\begin{array}{l}M\\
m\end{array}\right) \left\| \partial^{p}_{t}\nabla^{m+1}\mathbf{u}\cdot\partial^{P-p}_{t}\nabla^{M-m}\mathbf{v}
\right\|^{2}_{L^{2}(\mathbb{R}^{n})}
\\[0.4cm]
\qquad\lesssim \sum\limits^{P}_{p=0}\sum\limits^{M}_{m=0}\left(
\begin{array}{l}P\\
p\end{array}\right)\left(
\begin{array}{l}M\\
m\end{array}\right) \left\| \partial^{p}_{t}\nabla^{m+1}\mathbf{u}\right\|^{2}_{L^{4}(\mathbb{R}^{n})}
\left\| \partial^{P-p}_{t}\nabla^{M-m}\mathbf{v}\right\|^{2}_{L^{4}(\mathbb{R}^{n})}
\\[0.4cm]
\qquad\lesssim\sum\limits^{P}_{p=0}\sum\limits^{M}_{m=0}
\left(\left \| \partial^{p}_{t}\nabla^{m+1}  \mathbf{u}\right\|^{2}_{L^{2}(\mathbb{R}^{n})}+\left \| \partial^{p}_{t}\nabla^{m+1}\Lambda^{n/4}  \mathbf{u}\right\|^{2}_{L^{2}(\mathbb{R}^{n})}\right)
\\[0.3cm]
\qquad\qquad\qquad\times
\left(\left\| \partial^{P-p}_{t}\nabla^{M-m } \mathbf{v}\right\|^{2}_{L^{2}(\mathbb{R}^{n})}+\left\| \partial^{P-p}_{t}\nabla^{M-m }\Lambda^{n/4} \mathbf{v}\right\|^{2}_{L^{2}(\mathbb{R}^{n})}\right)
\\[0.3cm]
\qquad\lesssim\sum\limits^{P}_{p=0}\sum\limits^{M}_{m=0}
\left(\left \| \partial^{p}_{t}\nabla^{m}  \mathbf{v}\right\|^{2}_{L^{2}(\mathbb{R}^{n})}+\left \| \partial^{p}_{t}\nabla^{m}\Lambda^{n/4}  \mathbf{v}\right\|^{2}_{L^{2}(\mathbb{R}^{n})}\right)
\\[0.3cm]
\qquad\qquad\qquad\times
\left(\left\| \partial^{P-p}_{t}\nabla^{M-m } \mathbf{v}\right\|^{2}_{L^{2}(\mathbb{R}^{n})}+\left\| \partial^{P-p}_{t}\nabla^{M-m }\Lambda^{n/4} \mathbf{v}\right\|^{2}_{L^{2}(\mathbb{R}^{n})}\right).
 \end{array}
\right.
\end{equation}
Combining \eqref{4.14} with \eqref{4.18} and \eqref{4.19} gives
\begin{equation}\label{4.20}
\left\| \partial^{P+1}_{t}\nabla^{M}\mathbf{v}\right\|^{2}_{L^{2}(\mathbb{R}^{n})}\lesssim\left\| \partial^{P}_{t}\mathbf{v}\right\|^{2}_{H^{M+\frac{n}{2}}_{0} (\mathbb{R}^{n})}.
\end{equation}
By another inductive discussion, for all $M,P$ such that $\displaystyle M+\frac{nP}{2}\leq K$ we get
 \begin{equation}\label{4.21}
 \left\| \partial^{P}_{t}\nabla^{M}\mathbf{v}\right\|^{2}_{L^{2}(\mathbb{R}^{n})}\lesssim\left\| \mathbf{v}\right\|^{2}_{H^{K}_{0} (\mathbb{R}^{n})}.
 \end{equation}
 Multiplying \eqref{4.13} by $\partial^{P}_{t}\nabla^{M}\mathbf{v}$, then integrating by parts with respect to space yield
 \begin{equation}\label{4.22}
\left.
\begin{array}{ll}
\displaystyle \frac{1}{2}\frac{d}{dt}\left\| \partial^{P}_{t} \nabla^{M} \mathbf{v} \right\|^{2}_{L^{2}(\mathbb{R}^{n})}+\nu\left\| \partial^{P}_{t} \nabla^{M}\Lambda^{s} \mathbf{v} \right\|^{2}_{L^{2}(\mathbb{R}^{n})}
\\[0.3cm]
\qquad \leq \left| \left\langle \partial^{p}_{t}\nabla^{M}\left(\mathbf{u}\cdot \nabla \mathbf{v}\right), \partial^{p}_{t}\nabla^{M}\mathbf{v}
\right\rangle\right|+ \left| \left\langle \partial^{p}_{t}\nabla^{M}\left(\mathbf{v}\cdot \nabla \mathbf{u}^{T}\right), \partial^{p}_{t}\nabla^{M}\mathbf{v}
\right\rangle\right|.
 \end{array}
\right.
\end{equation}
 Making the similar argument employed in the proof of Theorem \ref{t4.1}, in particular, for $\displaystyle s=\frac{n}{4}$,
by the assumptions of this theorem:\quad $\left\| \mathbf{v}_{0}\right\| _{H_{0}^{K}(\mathbb{R}^{n})}\leq \varepsilon^{**} $ for $\varepsilon^{**}$ sufficiently small, combining \eqref{4.17} with \eqref{4.21}, for $\displaystyle\frac{n}{4}\leq s<1$ with $n=2,3$, we finish the proof of Theorem \ref{t4.2}.\hfill$\Box$
\\[0.3cm]
\indent  By Theorem \ref{t4.1} and Theorem \ref{t4.2}, we complete the proof of Theorem \ref{t1.2}.\hfill$\Box$

\section{Uniqueness}
In this section we shall show Theorem \ref{t1.3}. Indeed, we will show the continuous dependence of weak solutions to the Cauchy problem \eqref{1.1}-\eqref{1.2} constructed in Theorem \ref{t1.1} on the initial data and, in particular, we show the uniqueness of the weak solutions.\\
\indent  Let $(\mathbf{v},\mathbf{u})$ and $(\mathbf{w},\mathbf{q})$ be two weak solutions of the Cauchy problem \eqref{1.1}-\eqref{1.2} with the same initial data on the interval $[0,T]$. Then the two solutions $(\mathbf{v},\mathbf{u})$ and $(\mathbf{w},\mathbf{q})$ satisfy the following equations in $\mathbb{R}^{n}$ ($n=2,3$):
\begin{equation}\label{5.1}
\left\{
\begin{array}{ll}
\partial_{t}(\mathbf{v}-\mathbf{w})+\nu(-\Delta)^{s}(\mathbf{v}-\mathbf{w})+\nabla \tau
+\mathbf{u}\cdot\nabla \mathbf{v}-\mathbf{q}\cdot\nabla \mathbf{w}+\mathbf{v}\cdot\nabla \mathbf{u}^{T}-\mathbf{w}\cdot\nabla \mathbf{q}^{T}=0,
\\[0.3cm]
(\mathbf{u}-\mathbf{q})-\alpha^{2}\Delta(\mathbf{u}-\mathbf{q})=\mathbf{v}-\mathbf{w},
\\[0.3cm]
\nabla\cdot (\mathbf{v}-\mathbf{w})=0,\nabla\cdot \mathbf{v}=\nabla\cdot \mathbf{u}=\nabla\cdot \mathbf{w}=\nabla\cdot \mathbf{q}=0,
\\[0.3cm]
(\mathbf{v}-\mathbf{w})(0,x)=(\mathbf{v}_{0}(x)-\mathbf{w}_{0}(x))=0.
\end{array}
\right.
\end{equation}
Here, $\nabla \tau$ denotes the difference of the pressures corresponding to $\mathbf{v}$ and $\mathbf{w}$, respectively, and
\begin{equation}\label{5.2}
\mathbf{u}-\alpha^{2}\Delta \mathbf{u}=\mathbf{v},~~\mathbf{q}-\alpha^{2}\Delta \mathbf{q}=\mathbf{w}.
\end{equation}
  Note that
  \begin{equation}\label{5.3}
\left\| \mathbf{v}-\mathbf{w}\right\|^{2}_{L^{2}(\mathbb{R}^{n})}=\left\| \mathbf{u}-\mathbf{q}\right\|^{2}_{L^{2}(\mathbb{R}^{n})}+2\alpha^{2}\left\| \nabla\left(\mathbf{u}-\mathbf{q}\right)\right\|^{2}_{L^{2}(\mathbb{R}^{n})}+\alpha^{4}\left\| \Delta\left(\mathbf{u}-\mathbf{q}\right)\right\|^{2}_{L^{2}(\mathbb{R}^{n})},
\end{equation}
 by virtue of \eqref{1.11}, we will show the uniqueness through two steps for $\displaystyle\frac{n}{4}\leq s< 1$ with $n=2,3$. We first show for any $t\in [0,T]$, $\displaystyle\left\| \mathbf{u}-\mathbf{q}\right\|^{2}_{L^{2}(\mathbb{R}^{n})}=0$ and
 $\displaystyle\left\| \nabla\left(\mathbf{u}-\mathbf{q}\right)\right\|^{2}_{L^{2}(\mathbb{R}^{n})}=0$, then verify $\displaystyle\left\| \Delta\left(\mathbf{u}-\mathbf{q}\right)\right\|^{2}_{L^{2}(\mathbb{R}^{n})}=0$.
 \\[0.3cm]
 {\bf Step 1}\quad We prove that for any  $t\in [0,T]$, $\displaystyle\left\| \mathbf{u}-\mathbf{q}\right\|^{2}_{L^{2}(\mathbb{R}^{n})}=0$ and
 $\displaystyle\left\| \nabla\left(\mathbf{u}-\mathbf{q}\right)\right\|^{2}_{L^{2}(\mathbb{R}^{n})}=0$.
 \\[0.3cm]
 \indent Multiplying the first equation in \eqref{5.1} by $\mathbf{u}-\mathbf{q}$ and integrating in space yields, after some integration by parts,
\begin{equation}\label{5.4}
\left.
\begin{array}{ll}
\displaystyle\frac{1}{2}\frac{d}{dt} \left(\left\| \mathbf{u}-\mathbf{q}\right\|^{2}_{L^{2}(\mathbb{R}^{n})}+\alpha^{2}\left\| \nabla\left(\mathbf{u}-\mathbf{q}\right)\right\|^{2}_{L^{2}(\mathbb{R}^{n})}\right)
\\[0.3cm]
\qquad\quad+\nu\left(\left\| \Lambda^{s} \left(\mathbf{u}-\mathbf{q}\right)\right\|^{2}_{L^{2}(\mathbb{R}^{n})}+\alpha^{2}\left\| \nabla\Lambda^{s} \left(\mathbf{u}-\mathbf{q}\right)\right\|^{2}_{L^{2}(\mathbb{R}^{n})}\right)
\\[0.3cm]
\qquad=\displaystyle-\int_{\mathbb{R}^{n}}\left(\mathbf{u}\cdot \nabla\mathbf{ v}-\mathbf{q}\cdot \nabla \mathbf{w}+\mathbf{v}\cdot \nabla \mathbf{u}^{T}-\mathbf{w}\cdot \nabla \mathbf{q}^{T}\right)\left(\mathbf{u}-\mathbf{q}\right)dx
\\[0.4cm]
\qquad\leq\displaystyle\left|\int_{\mathbb{R}^{n}}\left(\mathbf{u}-\mathbf{q}\right) \nabla\left(\mathbf{u}-\alpha^{2}\Delta \mathbf{u}\right) \left(\mathbf{u}-\mathbf{q}\right)dx\right|
\\[0.4cm]
\qquad\quad+\displaystyle\left|\int_{\mathbb{R}^{n}}\mathbf{q} \nabla\left[\left(\mathbf{u}-\mathbf{q}\right)-\alpha^{2}\Delta \left(\mathbf{u}-\mathbf{q}\right) \right]\left(\mathbf{u}-\mathbf{q}\right)dx\right|
\\[0.4cm]
\qquad\quad+\displaystyle\left|\int_{ \mathbb{R}^{n}}\left(\mathbf{u}-\alpha^{2}\Delta \mathbf{u}\right) \nabla \left(\mathbf{u}^{T}-\mathbf{q}^{T}\right) \left(\mathbf{u}-\mathbf{q}\right)dx\right|
\\[0.4cm]
  \qquad\quad+\displaystyle\left|\int_{\mathbb{R}^{n}}\left[\left(\mathbf{u}-\mathbf{q}\right)
  -\alpha^{2}\Delta \left(\mathbf{u}-\mathbf{q}\right) \right]\nabla \mathbf{q}^{T}\left(\mathbf{u}-\mathbf{q}\right)dx\right|
  \\[0.4cm]
\qquad:=I_{1}+I_{2}+I_{3}+I_{4}.
  \end{array}
\right.
\end{equation}
There are two cases to consider for estimating \eqref{5.4}:
\\[0.3cm]
 \indent {\bf Case (I)}~~ $\displaystyle\frac{n}{4}<s< 1$ for $n=2,3$;\\
 \indent {\bf Case (II)}~~ $s=\displaystyle\frac{n}{4}$ for $n=2,3$.
 \\[0.3cm]
 We first deal with {\bf Case (I)}~~ $\displaystyle\frac{n}{4}<s< 1$  for $n=2,3$.
 \\[0.3cm]
  \indent Thanks to  H\"{o}lder's inequality, Gagliardo-Nirenberg-Sobolev inequality, interpolation inequality, Lemma \ref{l2.5} and the third equation in \eqref{5.1},  noting that the $i$th component of $\mathbf{v}\cdot \nabla \mathbf{u}^{T}$ is $\left(\mathbf{v}\cdot \nabla \mathbf{u}^{T}\right)_{i}=\sum\limits^{n}_{j=1} v_{j}\partial_{i}u_{j}$,  one deduces the following a priori estimates for $I_{i}$, $i=1,2,3,4$ :
\begin{equation}\label{5.5}
\left.
\begin{array}{ll}
I_{1}=\left|\displaystyle\int_{\mathbb{R}^{n}}\left(\mathbf{u}-\mathbf{q}\right) \nabla\left(\mathbf{u}-\alpha^{2}\Delta \mathbf{u}\right) \left(\mathbf{u}-\mathbf{q}\right)dx\right|
\\[0.4cm]
 \qquad=\left|\displaystyle\int_{\mathbb{R}^{n}}\left(\mathbf{u}-\mathbf{q}\right) \nabla \mathbf{ u}\left(\mathbf{u}-\mathbf{q}\right)dx-
 \alpha^{2}\displaystyle\int_{\mathbb{R}^{n}}\left(\mathbf{u}-\mathbf{q}\right) \nabla \Delta \mathbf{u}  \left(\mathbf{u}-\mathbf{q}\right)dx\right|
 \\[0.4cm]
\qquad\lesssim\left|\displaystyle\int_{\mathbb{R}^{n}}\left(\mathbf{u}-\mathbf{q}\right) \nabla \mathbf{u}\left(\mathbf{u}-\mathbf{q}\right)dx\right|
+\alpha^{2}\left|\displaystyle\int_{\mathbb{R}^{n}}\left(\mathbf{u}-\mathbf{q}\right) \nabla \Delta \mathbf{u}  \left(\mathbf{u}-\mathbf{q}\right)dx\right|
\\[0.4cm]
\qquad\lesssim \left\|\nabla \mathbf{u} \right\|_{L^{2}(\mathbb{R}^{n})}
\left\|\mathbf{u}-\mathbf{q}\right\|^{2}_{L^{4}(\mathbb{R}^{n})}+
\left\|\Delta \mathbf{u}\right\|_{L^{2}(\mathbb{R}^{n})}\left\| \left(\mathbf{u}-\mathbf{q}\right)\nabla\left(\mathbf{u}-\mathbf{q}\right)\right\|_{L^{2}(\mathbb{R}^{n})}
\\[0.3cm]
\qquad\lesssim\left\|\nabla \mathbf{u}\right\|_{L^{2}(\mathbb{R}^{n})}\left(\varepsilon \left\| \Lambda^{s} \left(\mathbf{u}-\mathbf{q}\right)\right\|^{2}_{L^{2}(\mathbb{R}^{n})}
+C(\varepsilon)\left\| \mathbf{u}-\mathbf{q} \right\|^{2}_{L^{2}(\mathbb{R}^{n})}\right)
\\[0.3cm]
 \qquad\quad +\left\|\Delta \mathbf{u}\right\|_{L^{2}(\mathbb{R}^{n})}\left\|\mathbf{u}-\mathbf{q}\right\|_{L^{4}(\mathbb{R}^{n})}
 \left\|\nabla\left(\mathbf{u}-\mathbf{q}\right)\right\|_{L^{4}(\mathbb{R}^{n})}
 \\[0.3cm]
  \qquad\lesssim\left\|\nabla \mathbf{u}\right\|_{L^{2}(\mathbb{R}^{n})}\left(\varepsilon \left\| \Lambda^{s} \left(\mathbf{u}-\mathbf{q}\right)\right\|^{2}_{L^{2}(\mathbb{R}^{n})}
  +C(\varepsilon)\left\| \mathbf{u}-\mathbf{q} \right\|^{2}_{L^{2}(\mathbb{R}^{n})}\right)
  \\[0.3cm]
  \qquad\quad + \left\|\Delta \mathbf{u}\right\|_{L^{2}(\mathbb{R}^{n})}\left(\varepsilon\left\|\Lambda^{s} \left(\mathbf{u}-\mathbf{q}\right)\right\|^{2}_{L^{2}(\mathbb{R}^{n})}+C(\varepsilon)
 \left\|\mathbf{u}-\mathbf{q}\right\|^{2}_{L^{2}(\mathbb{R}^{n})}\right)
 \\[0.3cm]
 \qquad\quad +\left\|\Delta \mathbf{u}\right\|_{L^{2}(\mathbb{R}^{n})}\left(\varepsilon\left\|\Lambda^{s}  \nabla \left(\mathbf{u}-\mathbf{q}\right)\right\|^{2}_{L^{2}(\mathbb{R}^{n})}+C(\varepsilon)
 \left\|\nabla\left(\mathbf{u}-\mathbf{q}\right)\right\|^{2}_{L^{2}(\mathbb{R}^{n})}\right)
 \\[0.3cm]
  \qquad\lesssim C(\varepsilon)\left( \left\|\nabla \mathbf{u}\right\|_{L^{2}(\mathbb{R}^{n})}+\left\|\Delta \mathbf{u}\right\|_{L^{2}(\mathbb{R}^{n})}\right)   \cdot\left(\left\|\mathbf{u}-\mathbf{q}\right\|^{2}_{L^{2}(\mathbb{R}^{n})}
  +\left\|\nabla\left(\mathbf{u}-\mathbf{q}\right)\right\|^{2}_{L^{2}(\mathbb{R}^{n})}\right)
  \\[0.3cm]
 \qquad\quad + \varepsilon\left(\left\|\nabla \mathbf{u} \right\|_{L^{2}(\mathbb{R}^{n})}+\left\| \Delta \mathbf{ u}\right\|_{L^{2}(\mathbb{R}^{n})}\right)\left( \left\|\Lambda^{s} \left(\mathbf{u}-\mathbf{q}\right)\right\|^{2}_{L^{2}(\mathbb{R}^{n})}
  +\alpha^{2}\left\|\nabla\Lambda^{s} \left(\mathbf{u}-\mathbf{q}\right)\right\|^{2}_{L^{2}(\mathbb{R}^{n})}\right);
    \end{array}
\right.
\end{equation}
 \begin{equation}\label{5.6}
\left.
\begin{array}{ll}
I_{2}=\left|\displaystyle\int_{\mathbb{R}^{n}}\mathbf{q}\nabla
\left[\left(\mathbf{u}-\mathbf{q}\right)-\alpha^{2}\Delta\left(u-q\right)\right]
   \left(\mathbf{u}-\mathbf{q}\right)dx\right|
    \\[0.4cm]
 \qquad=\left|\alpha^{2}\displaystyle\int_{\mathbb{R}^{n}}\mathbf{q}\nabla\Delta\left(\mathbf{u}-\mathbf{q}\right)
 \left(\mathbf{u}-\mathbf{q}\right)dx\right|
  \\[0.4cm]
\qquad\lesssim\left| \displaystyle\int_{\mathbb{R}^{n}}\mathbf{q}\nabla\left(\mathbf{u}-\mathbf{q}\right)
\Delta\left(\mathbf{u}-\mathbf{q}\right)dx\right|
 \\[0.4cm]
 \qquad\lesssim\left| \displaystyle\int_{\mathbb{R}^{n}}\nabla \mathbf{q}\nabla\left(\mathbf{u}-\mathbf{q}\right)\nabla\left(\mathbf{u}-\mathbf{q}\right) dx+\displaystyle\int_{\mathbb{R}^{n}} \mathbf{q}\nabla^{2}\left(\mathbf{u}-\mathbf{q}\right)\nabla\left(\mathbf{u}-\mathbf{q}\right) dx\right|
  \\[0.4cm]
 \qquad\lesssim\left| \displaystyle\int_{\mathbb{R}^{n}}\nabla \mathbf{q}\left|\nabla\left(\mathbf{u}-\mathbf{q}\right)\right|^{2} dx\right|
  \\[0.4cm]
 \qquad\lesssim  \left\|\nabla \mathbf{ q}\right\|_{L^{2}(\mathbb{R}^{n})}
 \left\|\nabla\left(\mathbf{u}-\mathbf{q}\right)\right\|^{2}_{L^{4}(\mathbb{R}^{n})}
  \\[0.3cm]
\qquad\lesssim\left\|\nabla \mathbf{q}\right\|_{L^{2}(\mathbb{R}^{n})}\left(\varepsilon \left\| \nabla\Lambda^{s} \left(\mathbf{u}-\mathbf{q}\right)\right\|^{2}_{L^{2}(\mathbb{R}^{n})}
+C(\varepsilon)\left\|\nabla\left( \mathbf{u}-\mathbf{q} \right) \right\|^{2}_{L^{2}(\mathbb{R}^{n})}\right).
    \end{array}
\right.
\end{equation}
\begin{equation}\label{5.7}
\left.
\begin{array}{ll}
\vspace{5pt}
  I_{3}=\left|\displaystyle\int_{\mathbb{R}^{n}}\left(\mathbf{u}-\alpha^{2}\Delta \mathbf{u}\right) \nabla \left(\mathbf{u}^{T}-\mathbf{q}^{T}\right)\left(\mathbf{u}-\mathbf{q}\right)dx\right|\\
 \vspace{5pt}
 \qquad\lesssim  \| \mathbf{u}\|_{L^{2}(\mathbb{R}^{n})}\left\| \left(\mathbf{u}-\mathbf{q}\right)\nabla\left(\mathbf{u}-\mathbf{q}\right)\right\|_{L^{2}(\mathbb{R}^{n})} +   \left\|\Delta \mathbf{u}\right\|_{L^{2}(\mathbb{R}^{n})}\left\|\nabla\left(\mathbf{u}-\mathbf{q}\right) \left(\mathbf{u}-\mathbf{q}\right)\right\|_{L^{2}(\mathbb{R}^{n})}\\
  \vspace{5pt}
  \qquad\lesssim  \|\mathbf{u}\|_{L^{2}(\mathbb{R}^{n})}\left(\left\| \left(\mathbf{u}-\mathbf{q}\right)\right\|^{2}_{L^{4}(\mathbb{R}^{n})}+\left\| \nabla\left(\mathbf{u}-\mathbf{q}\right)\right\|^{2}_{L^{4}(\mathbb{R}^{n})}\right)\\
  \vspace{5pt}
  \qquad\quad+   \left\|\Delta \mathbf{u}\right\|_{L^{2}(\mathbb{R}^{n})}\left(\left\| \left(\mathbf{u}-\mathbf{q}\right)\right\|^{2}_{L^{4}(\mathbb{R}^{n})}+\left\| \nabla\left(\mathbf{u}-\mathbf{q}\right)\right\|^{2}_{L^{4}(\mathbb{R}^{n})}\right)\\
  \vspace{5pt}
  \qquad\lesssim
  \left(\|\mathbf{u}\|_{L^{2}(\mathbb{R}^{n})}+\left\|\Delta \mathbf{u}\right\|_{L^{2}(\mathbb{R}^{n})}\right)\left(\varepsilon \left\| \Lambda^{s}\left(\mathbf{u}-\mathbf{q}\right) \right\|^{2}_{L^{2}(\mathbb{R}^{n})}+C(\varepsilon)\left\| \mathbf{u}-\mathbf{q}\right\|^{2}_{L^{2}(\mathbb{R}^{n})}\right)\\
     \vspace{5pt}
     \qquad\quad+ \left(\|\mathbf{u}\|_{L^{2}(\mathbb{R}^{n})}+\left\|\Delta \mathbf{u}\right\|_{L^{2}(\mathbb{R}^{n})}\right)\left(\varepsilon \left\| \nabla\Lambda^{s} \left(\mathbf{u}-\mathbf{q}\right)\right\|^{2}_{L^{2}(\mathbb{R}^{n})}+C(\varepsilon)\left\|\nabla\left( \mathbf{u}-\mathbf{q} \right) \right\|^{2}_{L^{2}(\mathbb{R}^{n})}\right)\\
 \vspace{5pt}
 \qquad\lesssim C(\varepsilon)\left(\|\mathbf{u}\|_{L^{2}(\mathbb{R}^{n})}+\left\|\Delta \mathbf{u}\right\|_{L^{2}(\mathbb{R}^{n})}\right)\cdot\left(  \left\| \mathbf{u}-\mathbf{q} \right\|^{2}_{L^{2}(\mathbb{R}^{n})}+ \left\|\nabla\left( \mathbf{u}-\mathbf{q} \right) \right\|^{2}_{L^{2}(\mathbb{R}^{n})}\right)\\
\vspace{5pt}
\qquad\quad+\varepsilon\left(\|\mathbf{u}\|_{L^{2}(\mathbb{R}^{n})}+\left\|\Delta \mathbf{u}\right\|_{L^{2}(\mathbb{R}^{n})}\right)\left(\left\| \Lambda^{s} \left(\mathbf{u}-\mathbf{q}\right)\right\|^{2}_{L^{2}(\mathbb{R}^{n})}+\left\|\nabla \Lambda^{s} \left(\mathbf{u}-\mathbf{q}\right)\right\|^{2}_{L^{2}(\mathbb{R}^{n})}\right).
      \end{array}
\right.
\end{equation}
\begin{equation}\label{5.8}
\left.
\begin{array}{ll}
\vspace{5pt}
|I_{4}|=\left|\displaystyle\int_{\mathbb{R}{n}}\left[\left(\mathbf{u}-\mathbf{q}\right)-\alpha^{2}\Delta \left(\mathbf{u}-\mathbf{q}\right)\right] \nabla \mathbf{q}^{T} \left(\mathbf{u}-\mathbf{q}\right)dx\right|\\
  \vspace{5pt}
  \qquad\lesssim  \left\|\mathbf{q}\right\|_{L^{2}(\mathbb{R}^{n})}\left\| \nabla\left(\mathbf{u}-\mathbf{q}\right)\left(\mathbf{u}-\mathbf{q}\right)\right\|_{L^{2}(\mathbb{R}^{n})} + \left\|\nabla\mathbf{ q}\right\|_{L^{2}(\mathbb{R}^{n})}\left\| \nabla\left(\mathbf{u}-\mathbf{q}\right)\right\|^{2}_{L^{4}(\mathbb{R}^{n})}\\
 \vspace{5pt}
 \qquad\quad+   \left\|\Delta \mathbf{q}\right\|_{L^{2}(\mathbb{R}^{n})}\left\|\nabla\left(\mathbf{u}-\mathbf{q}\right) \left(\mathbf{u}-\mathbf{q}\right)\right\|_{L^{2}(\mathbb{R}^{n})}  \\
    \vspace{5pt}
    \qquad\lesssim  \left\| \mathbf{q}\right\|_{L^{2}(\mathbb{R}^{n})}\left\| \nabla\left(\mathbf{u}-\mathbf{q}\right)\right\|_{L^{4}(\mathbb{R}^{n})}\left\| \mathbf{u}-\mathbf{q} \right \|_{L^{4}(\mathbb{R}^{n})}\\
   \vspace{5pt}
    \qquad\quad+\left \|\nabla \mathbf{q}\right\|_{L^{2}(\mathbb{R}^{n})}\left(\varepsilon \left\|\nabla\Lambda^{s} \left( \mathbf{u}-\mathbf{q}\right) \right\|^{2}_{L^{2}(\mathbb{R}^{n})}+C(\varepsilon)\left\|\nabla\left( \mathbf{u}-\mathbf{q}\right ) \right\|^{2}_{L^{2}(\mathbb{R}^{n})}\right)\\
  \vspace{5pt}
  \qquad\quad+   \left\|\Delta \mathbf{q}\right\|_{L^{2}(\mathbb{R}^{n})}\left\|\nabla\left(\mathbf{u}-\mathbf{q}\right) \right\|_{L^{4}(\mathbb{R}^{n})}\left\| \mathbf{u}-\mathbf{q} \right \|_{L^{4}(\mathbb{R}^{n})}\\
    \vspace{5pt}
    \qquad\lesssim  \left\| \mathbf{q}\right\|_{L^{2}(\mathbb{R}^{n})}\left(\left\| \nabla\left(\mathbf{u}-\mathbf{q}\right)\right\|^{2}_{L^{4}(\mathbb{R}^{n})}+\left\| \mathbf{u}-\mathbf{q} \right\|^{2}_{L^{4}(\mathbb{R}^{n})}\right)\\
    \vspace{5pt}
    \qquad\quad +\left\|\nabla \mathbf{q}\right\|_{L^{2}(\mathbb{R}^{n})}\left(\varepsilon \left\|\nabla\Lambda^{s} \left( \mathbf{u}-\mathbf{q}\right) \right\|^{2}_{L^{2}(\mathbb{R}^{n})}+C(\varepsilon)\left\|\nabla\left( \mathbf{u}-\mathbf{q} \right) \right\|^{2}_{L^{2}(\mathbb{R}^{n})}\right)\\
  \vspace{5pt}
  \qquad\quad+  \left \|\Delta \mathbf{q}\right\|_{L^{2}(\mathbb{R}^{n})}\left( \left\| \mathbf{u}-\mathbf{q}  \right\|^{2}_{L^{4}(\mathbb{R}^{n})}+\left\|\nabla\left(\mathbf{u}-\mathbf{q}\right) \right\|^{2}_{L^{4}(\mathbb{R}^{n})}\right)\\
  \vspace{5pt}
  \qquad\lesssim \left(\left\|\mathbf{ q}\right\|_{L^{2}(\mathbb{R}^{n})}+\left\|\nabla \mathbf{q}\right\|_{L^{2}(\mathbb{R}^{n})}+\left\| \Delta \mathbf{q}\right\|_{L^{2}(\mathbb{R}^{n})}\right)\\
  \vspace{5pt}
  \qquad\qquad\cdot\left(\varepsilon \left\|\nabla\Lambda^{s}\left( \mathbf{u}-\mathbf{q}\right)\right \|^{2}_{L^{2}(\mathbb{R}^{n})}+C(\varepsilon)\left\|\nabla\left( \mathbf{u}-\mathbf{q}\right ) \right\|^{2}_{L^{2}(\mathbb{R}^{n})}\right)\\
  \vspace{5pt}
  \qquad\quad+ \left(\left\|\mathbf{ q}\right\|_{L^{2}(\mathbb{R}^{n})}+\left\|\nabla \mathbf{q}\right\|_{L^{2}(\mathbb{R}^{n})}+\left\| \Delta \mathbf{q}\right\|_{L^{2}(\mathbb{R}^{n})}\right)  \\
  \vspace{5pt}
  \qquad\qquad\cdot  \left(\varepsilon \left\| \Lambda^{s} \left( \mathbf{u}-\mathbf{q}\right)\right \|^{2}_{L^{2}(\mathbb{R}^{n})}+C(\varepsilon)\left\| \mathbf{u}-\mathbf{q}   \right\|^{2}_{L^{2}(\mathbb{R}^{n})}\right)\\
  \vspace{5pt}
   \qquad\lesssim C(\varepsilon) \left( \left\| \mathbf{q}\right\|_{L^{2}(\mathbb{R}^{n})}+\left\|\nabla  \mathbf{q}\right\|_{L^{2}(\mathbb{R}^{n})}+\left\| \Delta \mathbf{q}\right\|_{L^{2}(\mathbb{R}^{n})}\right)\\
  \vspace{5pt}
  \qquad\qquad\cdot   \left(\left\|  \mathbf{u}-\mathbf{q}   \right\|^{2}_{L^{2}(\mathbb{R}^{n})}+\alpha^{2}\left\|\nabla\left( \mathbf{u}-\mathbf{q} \right) \right\|^{2}_{L^{2}(\mathbb{R}^{n})}\right)\\
   \vspace{5pt}
   \qquad\quad+ \varepsilon \left( \left\|\mathbf{q}\right\|_{L^{2}(\mathbb{R}^{n})}+\left\| \nabla \mathbf{q}\right\|_{L^{2}(\mathbb{R}^{n})}+\left \|\Delta \mathbf{q}\right\|_{L^{2}(\mathbb{R}^{n})}\right)\\
  \vspace{5pt}
  \qquad\qquad\cdot\left(\left\| \Lambda^{s} \left( \mathbf{u}-\mathbf{q}\right) \right\|^{2}_{L^{2}(\mathbb{R}^{n})}+\alpha^{2}\left\|\nabla\Lambda^{s} \left( \mathbf{ u}-\mathbf{q}\right) \right\|^{2}_{L^{2}(\mathbb{R}^{n})}\right).
   \end{array}
\right.
\end{equation}
Thanks to  \eqref{1.11}, \eqref{2.2} and \eqref{2.3}, choosing $\varepsilon$ small sufficiently such that
$$ \varepsilon \left( \left\|\mathbf{q}\right\|^{2}_{L^{2}(\mathbb{R}^{n})}+\left\|\nabla \mathbf{q}\right\|^{2}_{L^{2}(\mathbb{R}^{n})}+\left\| \Delta \mathbf{q}\right\|_{L^{2}(\mathbb{R}^{n})} +\left\|\mathbf{ u}\right\|_{L^{2}(\mathbb{R}^{n})}+ \left\| \nabla \mathbf{u}\right\|_{L^{2}(\mathbb{R}^{n})}+\left\| \Delta \mathbf{u}\right\|_{L^{2}(\mathbb{R}^{n})} \right)<\frac{\nu}{2},$$ combining \eqref{5.4} with \eqref{5.5}, \eqref{5.6}, \eqref{5.7} and \eqref{5.8} yields that
\begin{equation}\label{5.9}
\left.
\begin{array}{ll}
\vspace{5pt}
\displaystyle\frac{1}{2}\frac{d}{dt} \left(\left\| \mathbf{u}-\mathbf{q}\right\|^{2}_{L^{2}(\mathbb{R}^{n})}+\alpha^{2}\left\| \nabla\left(\mathbf{u}-\mathbf{q}\right)\right\|^{2}_{L^{2}(\mathbb{R}^{n})}\right)\\
  \vspace{5pt}
  \qquad \quad \displaystyle+\frac{\nu}{2}\left(\left\| \Lambda^{s} \left(\mathbf{u}-\mathbf{q}\right)\right\|^{2}_{L^{2}(\mathbb{R}^{n})}+\alpha^{2}\left\| \nabla\Lambda^{s} \left (\mathbf{u}-\mathbf{q}\right)\right\|^{2}_{L^{2}(\mathbb{R}^{n})}\right)\\
 \vspace{5pt}\qquad\lesssim\left\|  \mathbf{u}-\mathbf{q}   \right\|^{2}_{L^{2}(\mathbb{R}^{n})}+\alpha^{2}\left\|\nabla\left( \mathbf{u}-\mathbf{q}\right ) \right \|^{2}_{L^{2}(\mathbb{R}^{n})}.
 \end{array}
\right.
\end{equation}
We then deal with  {\bf Case (II)} ~~ $\displaystyle s=\frac{n}{4}$ for $n=2,3$.
 \\[0.3cm]
\indent In this case, note that \eqref{5.4}, making the similar a priori estimates to those employed in \eqref{5.5}-\eqref{5.9} imply
\begin{equation}\label{5.10}
\left.
\begin{array}{ll}
\vspace{5pt}
I_{1}=\left|\displaystyle\int_{\mathbb{R}^{n}}\left(\mathbf{u}-\mathbf{q}\right) \nabla\left(\mathbf{u}-\alpha^{2}\Delta \mathbf{u}\right) \left(\mathbf{u}-\mathbf{q}\right)dx\right|\\
 \vspace{5pt}
 \qquad\lesssim  \left\|\nabla \mathbf{u} \right\|_{L^{2}(\mathbb{R}^{n})}\left\|\mathbf{u}-\mathbf{q}\right\|^{2}_{L^{4}(\mathbb{R}^{n})}+ \left\|\Delta \mathbf{u}  \right\|_{L^{2}(\mathbb{R}^{n})}\left\|\left(\mathbf{u}-\mathbf{q}\right)\nabla\left(u-q\right)\right\|
 _{L^{2}(\mathbb{R}^{n})}\\
 \vspace{5pt}
 \qquad\lesssim\left\|\nabla \mathbf{u}\right\|_{L^{2}(\mathbb{R}^{n})}\left(  \left\| \Lambda^{n/4} \left(\mathbf{u}-\mathbf{q}\right)\right\|^{2}_{L^{2}(\mathbb{R}^{n})}+ \left\| \mathbf{u}-\mathbf{q} \right \|^{2}_{L^{2}(\mathbb{R}^{n})}\right)\\
  \vspace{5pt}
  \qquad\quad + \left\|\Delta \mathbf{u}\right\|_{L^{2}(\mathbb{R}^{n})}\left\|\nabla\left(\mathbf{u}-\mathbf{q}\right)\right\|_{L^{4}(\mathbb{R}^{n})}
 \left\| \mathbf{u}-\mathbf{q} \right\|_{L^{4}(\mathbb{R}^{n})}\\
 \vspace{5pt}
 \qquad\lesssim\left\|\nabla \mathbf{u}\right\|_{L^{2}(\mathbb{R}^{n})}\left(  \left\| \Lambda^{n/4} \left(\mathbf{u}-\mathbf{q}\right)\right\|^{2}_{L^{2}(\mathbb{R}^{n})}+\left \| \mathbf{u}-\mathbf{q} \right\|^{2}_{L^{2}(\mathbb{R}^{n})}\right)\\
 \vspace{5pt}
  \qquad\quad + \left\|\Delta \mathbf{u}\right\|_{L^{2}(\mathbb{R}^{n})}\left(  \left\| \Lambda^{n/4} \left(\mathbf{u}-\mathbf{q}\right)\right\|^{2}_{L^{2}(\mathbb{R}^{n})}+ \left\| \mathbf{u}-\mathbf{q }\right\|^{2}_{L^{2}(\mathbb{R}^{n})}\right)\\
 \vspace{5pt}
 \qquad\quad+ \left\|\Delta \mathbf{u}\right\|_{L^{2}(\mathbb{R}^{n})}\left(  \left\|\nabla \Lambda^{n/4} \left(\mathbf{u}-\mathbf{q}\right)\right\|^{2}_{L^{2}(\mathbb{R}^{n})}+ \left\| \nabla\left(\mathbf{u}-\mathbf{q}\right) \right\|^{2}_{L^{2}(\mathbb{R}^{n})}\right)\\
  \vspace{5pt}
  \qquad\lesssim \left(\left\|\nabla \mathbf{u}\right\|_{L^{2}(\mathbb{R}^{n})}+\left\|\Delta \mathbf{u}\right\|_{L^{2}(\mathbb{R}^{n})} \right)\cdot \left( \left\|  \left(\mathbf{u}-\mathbf{q}\right)\right\|^{2}_{L^{2}(\mathbb{R}^{n})}+\alpha^{2}
 \left\|\nabla\left(\mathbf{u}-\mathbf{q}\right)\right\|^{2}_{L^{2}(\mathbb{R}^{n})}\right) \\
\vspace{5pt}
\qquad\quad + \left(\left\|\nabla \mathbf{u}\right\|_{L^{2}(\mathbb{R}^{n})}+\left\|\Delta \mathbf{u}\right\|_{L^{2}(\mathbb{R}^{n})} \right)\left(\left\|\nabla\Lambda^{n/4} \left(\mathbf{u}-\mathbf{q}\right)\right\|^{2}_{L^{2}(\mathbb{R}^{n})}+\left\|\Lambda^{n/4} \left(\mathbf{u}-\mathbf{q}\right)\right\|^{2}_{L^{2}(\mathbb{R}^{n})}\right).
   \end{array}
\right.
\end{equation}
 \begin{equation}\label{5.11}
\left.
\begin{array}{ll}
\vspace{5pt}
  |I_{2}|=\left|\displaystyle\int_{\mathbb{R}^{n}}\mathbf{q}\nabla\left[\left(\mathbf{u}-\mathbf{q}\right)
  -\alpha^{2}\Delta
  \left(\mathbf{u}-\mathbf{q}\right)\right]
   \left(\mathbf{u}-\mathbf{q}\right)dx\right|\\
 \vspace{5pt}
 \qquad=\left|\alpha^{2}\displaystyle\int_{\mathbb{R}^{n}}q\nabla\Delta\left(\mathbf{u}-\mathbf{q}\right)
 \left(\mathbf{u}-\mathbf{q}\right)dx\right|\\
 \vspace{5pt}
 \qquad\lesssim\left| \displaystyle\int_{\mathbb{R}^{n}} \mathbf{q} \nabla\left(\mathbf{u}-\mathbf{q}\right)\Delta \left(\mathbf{u}-\mathbf{q}\right)  dx\right|\\
 \vspace{5pt}
 \qquad\lesssim\left| \displaystyle\int_{\mathbb{R}^{n}} \nabla \mathbf{q} \nabla\left(\mathbf{u}-\mathbf{q}\right)\nabla \left(\mathbf{u}-\mathbf{q}\right)  dx+\displaystyle\int_{\mathbb{R}^{n}}  \mathbf{q} \nabla^{2}\left(\mathbf{u}-\mathbf{q}\right)\nabla \left(\mathbf{u}-\mathbf{q}\right)  dx\right|\\
\vspace{5pt}
\qquad\lesssim \left| \displaystyle\int_{\mathbb{R}^{n}} \nabla \mathbf{q } \left|\nabla\left(\mathbf{u}-\mathbf{q}\right)\right|^{2}dx\right|\\
\vspace{5pt}
\qquad\lesssim  \left\|\nabla \mathbf{q}\right\|_{L^{2}(\mathbb{R}^{n})}\left\|\nabla\left(\mathbf{u}-\mathbf{q}\right)\right\|^{2}
_{L^{4}(\mathbb{R}^{n})}\\
 \vspace{5pt}
 \qquad\lesssim \left\|\nabla \mathbf{q}\right\|_{L^{2}(\mathbb{R}^{n})}\left(  \left\| \nabla\Lambda^{n/4} \left(\mathbf{u}-\mathbf{q}\right)\right\|^{2}_{L^{2}(\mathbb{R}^{n})}+ \left\|\nabla\left( \mathbf{u}-\mathbf{q }\right) \right\|^{2}_{L^{2}(\mathbb{R}^{n})}\right).
 \end{array}
\right.
\end{equation}
\begin{equation}\label{5.12}
\left.
\begin{array}{ll}
\vspace{5pt}
   I_{3}=\left|\displaystyle\int_{\mathbb{R}^{n}}\left(\mathbf{u}-\alpha^{2}\Delta \mathbf{u}\right) \nabla \left(\mathbf{u}^{T}-\mathbf{q}^{T}\right)\left(\mathbf{u}-\mathbf{q}\right)dx\right|\\
  \vspace{5pt}
  \qquad\lesssim  \left(\| \mathbf{u}\|_{L^{2}(\mathbb{R}^{n})}+\left\| \Delta \mathbf{u}\right\|_{L^{2}(\mathbb{R}^{n})}\right)\left\| \left(\mathbf{u}-\mathbf{q}\right)\nabla\left(\mathbf{u}-\mathbf{q}\right)\right\|_{L^{2}(\mathbb{R}^{n})}\\
 \vspace{5pt}
  \qquad\lesssim  \left(\left\| \mathbf{u}\right\|_{L^{2}(\mathbb{R}^{n})}+\left\| \Delta \mathbf{u}\right\|_{L^{2}(\mathbb{R}^{n})}\right)\left\|  \mathbf{u}-\mathbf{q} \right \|_{L^{4}(\mathbb{R}^{n})}\left\|  \nabla\left(\mathbf{u}-\mathbf{q}\right)  \right\|_{L^{4}(\mathbb{R}^{n})}\\
 \vspace{5pt}
 \qquad\lesssim \left(\| \mathbf{u}\|_{L^{2}(\mathbb{R}^{n})}+\left\| \Delta \mathbf{u}\right\|_{L^{2}(\mathbb{R}^{n})}\right)\left(\left\|  \mathbf{u}-\mathbf{q}  \right\|^{2}_{L^{4}(\mathbb{R}^{n})}+\left\|  \nabla\left(\mathbf{u}-\mathbf{q}\right)  \right\|^{2}_{L^{4}(\mathbb{R}^{n})}\right)\\
 \vspace{5pt}
 \qquad\lesssim \left(\| \mathbf{u}\|_{L^{2}(\mathbb{R}^{n})}+\left\| \Delta \mathbf{u}\right\|_{L^{2}(\mathbb{R}^{n})}\right)\left(  \left\| \Lambda^{n/4} \left(\mathbf{u}-\mathbf{q}\right)\right\|^{2}_{L^{2}(\mathbb{R}^{n})}+ \left\| \mathbf{u}-\mathbf{q} \right\|^{2}_{L^{2}(\mathbb{R}^{n})}\right)
 \\
 \vspace{5pt}
 \qquad\quad+\left(\| \mathbf{u}\|_{L^{2}(\mathbb{R}^{n})}+\left\| \Delta \mathbf{u}\right\|_{L^{2}(\mathbb{R}^{n})}\right)\left(  \left\| \nabla\Lambda^{n/4} \left(\mathbf{u}-\mathbf{q}\right)\right\|^{2}_{L^{2}(\mathbb{R}^{n})}+ \left\| \nabla \left(\mathbf{u}-\mathbf{q}\right) \right\|^{2}_{L^{2}(\mathbb{R}^{n})}\right)\\
 \vspace{5pt}
 \qquad\lesssim\left(\| \mathbf{u}\|_{L^{2}(\mathbb{R}^{n})}+\left\| \Delta \mathbf{u}\right\|_{L^{2}(\mathbb{R}^{n})}\right)\left( \left\| \mathbf{u}-\mathbf{q} \right\|^{2}_{L^{2}(\mathbb{R}^{n})}+ \left\| \nabla \left(\mathbf{u}-\mathbf{q}\right) \right\|^{2}_{L^{2}(\mathbb{R}^{n})}\right)\\
 \vspace{5pt}
 \qquad\quad+\left(\| \mathbf{u}\|_{L^{2}(\mathbb{R}^{n})}+\left\| \Delta \mathbf{u}\right\|_{L^{2}(\mathbb{R}^{n})}\right)\left( \left\| \Lambda^{n/4} \left(\mathbf{u}-\mathbf{q}\right) \right\|^{2}_{L^{2}(\mathbb{R}^{n})}+ \left\| \nabla \Lambda^{n/4} \left(\mathbf{u}-\mathbf{q}\right)  \right\|^{2}_{L^{2}(\mathbb{R}^{n})}\right).
 \end{array}
\right.
\end{equation}

\begin{equation}\label{5.13}
\left.
\begin{array}{ll}
\vspace{5pt}
I_{4}=\left|\displaystyle\int_{\mathbb{R}^{n}}\left[\left(\mathbf{u}-\mathbf{q}\right)-\alpha^{2}\Delta \left(\mathbf{u}-\mathbf{q}\right)\right] \nabla \mathbf{q}^{T} \left(\mathbf{u}-\mathbf{q}\right)dx\right|\\
 \vspace{5pt}
 \qquad\lesssim \left\| \nabla \mathbf{q}\right\|_{L^{2}(\mathbb{R}^{n})}\left(\left\| \mathbf{u}-\mathbf{q}\right\|^{2}_{L^{4}(\mathbb{R}^{n})}+\left\| \nabla\left(\mathbf{u}-\mathbf{q}\right)\right\|^{2}_{L^{4}(\mathbb{R}^{n})}\right)\\
 \vspace{5pt}
 \qquad\quad+\left\| \Delta \mathbf{q}\right\|_{L^{2}(\mathbb{R}^{n})}\left\| \left(\mathbf{u}-\mathbf{q}\right)\nabla\left(u-q\right)\right\|_{L^{2}(\mathbb{R}^{n})} \\
  \vspace{5pt}
  \qquad\lesssim \left(\left\| \nabla \mathbf{q}\right\|_{L^{2}(\mathbb{R}^{n})}+\left\| \Delta \mathbf{q}\right\|_{L^{2}(\mathbb{R}^{n})}\right)\left(\left\| \mathbf{u}-\mathbf{q}\right\|^{2}_{L^{4}(\mathbb{R}^{n})}+\left\| \nabla\left(\mathbf{u}-\mathbf{q}\right)\right\|^{2}_{L^{4}(\mathbb{R}^{n})}\right)\\
      \vspace{5pt}
    \qquad\lesssim \left(\left\| \nabla \mathbf{q}\right\|_{L^{2}(\mathbb{R}^{n})}+\left\| \Delta \mathbf{q}\right\|_{L^{2}(\mathbb{R}^{n})}\right)\left(\left\| \Lambda^{n/4} \left( \mathbf{u}-\mathbf{q}\right) \right\|^{2}_{L^{2}(\mathbb{R}^{n})}+\left\| \mathbf{u}-\mathbf{q} \right\|^{2}_{L^{2}(\mathbb{R}^{n})}\right)\\
   \vspace{5pt}
 \qquad\quad+\left(\left\| \nabla \mathbf{q}\right\|_{L^{2}(\mathbb{R}^{n})}+\left\| \Delta \mathbf{q}\right\|_{L^{2}(\mathbb{R}^{n})}\right)\left(\left\| \nabla\Lambda^{n/4} \left( \mathbf{u}-\mathbf{q}\right) \right\|^{2}_{L^{2}(\mathbb{R}^{n})}+\left\| \nabla\left(\mathbf{u}-\mathbf{q }\right)\right\|^{2}_{L^{2}(\mathbb{R}^{n})}\right)\\
    \vspace{5pt}
 \qquad\lesssim \left(\left\| \nabla \mathbf{q}\right\|_{L^{2}(\mathbb{R}^{n})}+\left\| \Delta \mathbf{q}\right\|_{L^{2}(\mathbb{R}^{n})}\right)\left(\left\| \mathbf{ u}-\mathbf{q} \right\|^{2}_{L^{2}(\mathbb{R}^{n})}+\left\| \nabla\left(\mathbf{u}-\mathbf{q}\right) \right\|^{2}_{L^{2}(\mathbb{R}^{n})}\right)\\
   \vspace{5pt}
 \qquad\quad+\left(\left\| \nabla \mathbf{q}\right\|_{L^{2}(\mathbb{R}^{n})}+\left\| \Delta \mathbf{q}\right\|_{L^{2}(\mathbb{R}^{n})}\right)\left(\left\|  \Lambda^{n/4} \left( \mathbf{u}-\mathbf{q}\right) \right\|^{2}_{L^{2}(\mathbb{R}^{n})}+\left\| \nabla\Lambda^{n/4} \left( \mathbf{u}-\mathbf{q}\right) \right\|^{2}_{L^{2}(\mathbb{R}^{n})}\right).
 \end{array}
\right.
\end{equation}
Recall \eqref{1.11} and \eqref{2.5}, due to
$$ \left( \left\|\nabla \mathbf{u}\right\|_{L^{2}(\mathbb{R}^{n})}+\left\| \nabla \mathbf{q}\right\|_{L^{2}(\mathbb{R}^{n})} +\left\|\Delta \mathbf{u}\right\|_{L^{2}(\mathbb{R}^{n})}+\left\| \Delta \mathbf{q}\right\|_{L^{2}(\mathbb{R}^{n})}\right)\lesssim \left\|\mathbf{v} \right\|^{2}_{L^{2}(\mathbb{R}^{n})}+\left\|\mathbf{w} \right\|^{2}_{L^{2}(\mathbb{R}^{n})}\lesssim \left\|\mathbf{v}_{0} \right\|^{2}_{L^{2}(\mathbb{R}^{n})},$$
note that  \eqref{5.4}, \eqref{5.10}-\eqref{5.13} and the assumption of Theorem \ref{t1.1} that for $\displaystyle s=\frac{n}{4}$,   $ \displaystyle \left\|\mathbf{v}_{0}\right\|^{2}_{L^{2}(\mathbb{R}^{n})}\lesssim \varepsilon^{*}$ for $\varepsilon^{*}$ sufficiently small,  in particular, choosing $ \displaystyle \varepsilon^{*}\leq \frac{\nu}{8}$, we infer that
\begin{equation}\label{5.14}
\left.
\begin{array}{ll}
\vspace{5pt}
\displaystyle\frac{1}{2}\frac{d}{dt} \left(\left\| \mathbf{u}-\mathbf{q}\right\|^{2}_{L^{2}(\mathbb{R}^{n})}+\alpha^{2}\left\| \nabla\left(\mathbf{u}-\mathbf{q}\right)\right\|^{2}_{L^{2}(\mathbb{R}^{n})}\right)\\
 \vspace{5pt}
 \qquad\quad \displaystyle+\frac{\nu}{2}\left(\left\| \Lambda^{n/4} \left(\mathbf{u}-\mathbf{q}\right)\right\|^{2}_{L^{2}(\mathbb{R}^{n})}+\alpha^{2}\left\| \nabla\Lambda^{n/4} \left(\mathbf{u}-\mathbf{q}\right)\right\|^{2}_{L^{2}(\mathbb{R}^{n})}\right)\\
 \vspace{5pt}
 \qquad\lesssim \displaystyle\left(\left\|  \mathbf{u}-\mathbf{q}   \right\|^{2}_{L^{2}(\mathbb{R}^{n})}+\alpha^{2}\left\|\nabla\left( \mathbf{u}-\mathbf{q}\right ) \right\|^{2}_{L^{2}(\mathbb{R}^{n})}\right).\end{array}
\right.
\end{equation}
Note that \eqref{5.3}, combining \eqref{5.9} with \eqref{5.14}, note that $\left\|  \mathbf{u}_{0}-\mathbf{q}_{0}   \right\|^{2}_{L^{2}(\mathbb{R}^{n})}=0$ and  $\left\| \nabla\left( \mathbf{u}_{0}-\mathbf{q}_{0}\right )  \right\|^{2}_{L^{2}(\mathbb{R}^{n})}=0$ from $\mathbf{v}_{0}-\mathbf{w}_{0}=0$,
Gronwall's inequality then yields that for $\displaystyle\frac{n}{4}\leq s< 1$ and for any $t\in [0,T]$,
\begin{equation}\label{5.15}
 \left\|\mathbf{u} -\mathbf{q}\right \|^{2}_{L^{2}(\mathbb{R}^{n})}+\alpha^{2}\left\|\nabla\left(\mathbf{u} -\mathbf{q}\right) \right\|^{2}_{L^{2}(\mathbb{R}^{n})}=0.
\end{equation}
In the following, we need to show $\left\|\Delta\left(\mathbf{u} -\mathbf{q}\right) \right\|^{2}_{L^{2}(\mathbb{R}^{n})}=0$.
\\[0.3cm]
 {\bf Step 2}\quad We prove that for any  $t\in [0,T]$, ~$\displaystyle\left\| \Delta\left(\mathbf{u}-\mathbf{q}\right)\right\|^{2}_{L^{2}(\mathbb{R}^{n})}=0$.
 \\[0.3cm]
 \indent Multiplying the first equation in \eqref{5.1} by $\Delta\left(\mathbf{u}-\mathbf{q}\right)$ and integrating in space yields, after some integration by parts,
\begin{equation}\label{5.16}
\left.
\begin{array}{ll}
\displaystyle\frac{1}{2}\frac{d}{dt} \left(\left\|\nabla\left( \mathbf{u}-\mathbf{q}\right)\right\|^{2}_{L^{2}(\mathbb{R}^{n})}+\alpha^{2}\left\| \Delta\left(\mathbf{u}-\mathbf{q}\right)\right\|^{2}_{L^{2}(\mathbb{R}^{n})}\right)
\\[0.3cm]
\qquad\quad+\nu\left(\left\| \Lambda^{s} \nabla \left(\mathbf{u}-\mathbf{q}\right)\right\|^{2}_{L^{2}(\mathbb{R}^{n})}+\alpha^{2}\left\| \Lambda^{s}\Delta \left(\mathbf{u}-\mathbf{q}\right)\right\|^{2}_{L^{2}(\mathbb{R}^{n})}\right)
\\[0.3cm]
\qquad=\displaystyle\int_{\mathbb{R}^{n}}\left(\mathbf{u}\cdot \nabla\mathbf{ v}-\mathbf{q}\cdot \nabla \mathbf{w}+\mathbf{v}\cdot \nabla \mathbf{u}^{T}-\mathbf{w}\cdot \nabla \mathbf{q}^{T}\right)\Delta\left(\mathbf{u}-\mathbf{q}\right)dx
\\[0.4cm]
\qquad\leq\displaystyle\left|\int_{\mathbb{R}^{n}}\left(\mathbf{u}-\mathbf{q}\right) \nabla\left(\mathbf{u}-\alpha^{2}\Delta \mathbf{u}\right) \Delta\left(\mathbf{u}-\mathbf{q}\right)dx\right|
\\[0.4cm]
\qquad\quad+\displaystyle\left|\int_{\mathbb{R}^{n}}\mathbf{q} \nabla\left[\left(\mathbf{u}-\mathbf{q}\right)-\alpha^{2}\Delta \left(\mathbf{u}-\mathbf{q}\right) \right]\Delta\left(\mathbf{u}-\mathbf{q}\right)dx\right|
\\[0.4cm]
\qquad\quad+\displaystyle\left|\int_{ \mathbb{R}^{n}}\left(\mathbf{u}-\alpha^{2}\Delta \mathbf{u}\right) \nabla \left(\mathbf{u}^{T}-\mathbf{q}^{T}\right) \Delta\left(\mathbf{u}-\mathbf{q}\right)dx\right|
\\[0.4cm]
  \qquad\quad+\displaystyle\left|\int_{\mathbb{R}^{n}}\left[\left(\mathbf{u}-\mathbf{q}\right)
  -\alpha^{2}\Delta \left(\mathbf{u}-\mathbf{q}\right) \right]\nabla \mathbf{q}^{T}\Delta\left(\mathbf{u}-\mathbf{q}\right)dx\right|
  \\[0.4cm]
\qquad:=A_{1}+A_{2}+A_{3}+A_{4}.
  \end{array}
\right.
\end{equation}
We will consider two cases to bound \eqref{5.16}:
\\[0.3cm]
 \indent {\bf Case 1}~~ $\displaystyle\frac{n}{4}<s< 1$ for $n=2,3$;\\
 \indent {\bf Case 2}~~ $s=\displaystyle\frac{n}{4}$ for $n=2,3$.
 \\[0.3cm]
 We first tackle {\bf Case 1}~~ $\displaystyle\frac{n}{4}<s< 1$  for $n=2,3$.
 \\[0.3cm]
 \indent In this case, by virtue of H\"{o}lder's inequality, Gagliardo-Nirenberg-Sobolev inequality and interpolation inequality, $A_{1},~A_{2},~A_{3},~A_{4}$ can be bounded as follows:
\begin{equation}\label{5.17}
\left.
\begin{array}{ll}
A_{1}=\left|\displaystyle\int_{\mathbb{R}^{n}}\left(\mathbf{u}-\mathbf{q}\right) \nabla\left(\mathbf{u}-\alpha^{2}\Delta \mathbf{u}\right) \Delta\left(\mathbf{u}-\mathbf{q}\right)dx\right|
\\[0.4cm]
 \qquad=\left|\displaystyle\int_{\mathbb{R}^{n}}\left(\mathbf{u}-\mathbf{q}\right) \nabla \mathbf{ u}\Delta\left(\mathbf{u}-\mathbf{q}\right)dx-
 \alpha^{2}\displaystyle\int_{\mathbb{R}^{n}}\left(\mathbf{u}-\mathbf{q}\right) \nabla \Delta \mathbf{u}  \Delta\left(\mathbf{u}-\mathbf{q}\right)dx\right|
 \\[0.4cm]
\qquad\lesssim\left|\displaystyle\int_{\mathbb{R}^{n}}\left(\mathbf{u}-\mathbf{q}\right) \nabla \mathbf{u}\Delta\left(\mathbf{u}-\mathbf{q}\right)dx\right|
+\alpha^{2}\left|\displaystyle\int_{\mathbb{R}^{n}}\left(\mathbf{u}-\mathbf{q}\right) \nabla \Delta \mathbf{u} \Delta \left(\mathbf{u}-\mathbf{q}\right)dx\right|.
   \end{array}
\right.
\end{equation}
Direct calculation gives
\begin{equation}\label{5.18}
\left.
\begin{array}{ll}
 \left|\displaystyle\int_{\mathbb{R}^{n}}\left(\mathbf{u}-\mathbf{q}\right) \nabla \mathbf{u}\Delta\left(\mathbf{u}-\mathbf{q}\right)dx\right| &\lesssim \left\|\nabla \mathbf{u} \right\|_{L^{2}(\mathbb{R}^{n})}
 \left\|\left(\mathbf{u}-\mathbf{q}\right)\Delta
 \left(\mathbf{u}-\mathbf{q}\right)\right\|_{L^{2}(\mathbb{R}^{n})}
 \\[0.3cm]
  &\lesssim \displaystyle\left\|\nabla \mathbf{u} \right\|_{L^{2}(\mathbb{R}^{n})}
 \left\| \mathbf{u}-\mathbf{q} \right\|_{L^{\frac{n}{s}}(\mathbb{R}^{n})}
 \left\| \Delta\left(\mathbf{u}-\mathbf{q}\right) \right\|_{L^{\frac{2n}{n-2s}}(\mathbb{R}^{n})}
 \\[0.3cm]
  &\lesssim \displaystyle\left\|\nabla \mathbf{u} \right\|^{2}_{L^{2}(\mathbb{R}^{n})}
 \left\| \nabla\left(\mathbf{u}-\mathbf{q}\right) \right\|^{2}_{L^{2}(\mathbb{R}^{n})}
 +\frac{\nu}{16}\left\|\Lambda^{s} \Delta\left(\mathbf{u}-\mathbf{q}\right) \right\|^{2}_{L^{2}(\mathbb{R}^{n})}.
   \end{array}
\right.
\end{equation}
On the other hand, by Lemma \ref{l2.4},  Lemma \ref{l2.7} and  Lemma \ref{l2.9}, we have
\begin{equation}\label{5.19}
\left.
\begin{array}{ll}
 &\left|\displaystyle\int_{\mathbb{R}^{n}}\left(\mathbf{u}-\mathbf{q}\right) \nabla\Delta \mathbf{u}\Delta\left(\mathbf{u}-\mathbf{q}\right)dx\right|
 \\[0.3cm]
  &\qquad \lesssim \displaystyle\left\|\Lambda^{s}\Delta \mathbf{u} \right\|_{L^{2}(\mathbb{R}^{n})}
  \left\|\Lambda^{1-s}\left[\left(\mathbf{u}-\mathbf{q}\right)
  \Delta\left(\mathbf{u}-\mathbf{q}\right)\right]\right\|_{L^{2}(\mathbb{R}^{n})}
 \\[0.3cm]
  &\qquad \lesssim\displaystyle \left\|\Lambda^{s}\Delta \mathbf{u} \right\|_{L^{2}(\mathbb{R}^{n})}
 \left\{\left\|\Lambda^{1-s}\left[\left(\mathbf{u}-\mathbf{q}\right)
  \Delta\left(\mathbf{u}-\mathbf{q}\right)\right]\right.\right.
   \\[0.4cm]
  &\qquad\qquad\qquad\qquad\qquad\left.-\left(\mathbf{u}-\mathbf{q}\right)
  \Lambda^{1-s}\Delta\left(\mathbf{u}-\mathbf{q}\right)
  -\Lambda^{1-s}\left(\mathbf{u}-\mathbf{q}\right)
  \Delta\left(\mathbf{u}-\mathbf{q}\right)\right\|_{L^{2}(\mathbb{R}^{n})}
  \\[0.3cm]
  &\qquad\qquad\qquad + \displaystyle
 \left.\left\| \left(\mathbf{u}-\mathbf{q}\right)
  \Lambda^{1-s}\Delta\left(\mathbf{u}-\mathbf{q}\right)
   \right\|_{L^{2}(\mathbb{R}^{n})}+\left\|
  \Lambda^{1-s}\left(\mathbf{u}-\mathbf{q}\right)
  \Delta\left(\mathbf{u}-\mathbf{q}\right)\right\|_{L^{2}(\mathbb{R}^{n})}\right\}.
    \end{array}
\right.
\end{equation}
By the aid of (I) of Lemma \ref{l2.7},  Lemma \ref{l2.9} and  \eqref{3.3}, note that $\displaystyle \frac{1}{2}=\frac{1}{n/s}+\frac{1}{2n/(n-2s)}$ and $0<1-s<s<1$, we can make the following estimates for the right hand side of \eqref{5.19}:
\begin{equation}\label{5.20}
\left.
\begin{array}{ll}
    & \displaystyle  \left\|\Lambda^{1-s}\left[\left(\mathbf{u}-\mathbf{q}\right)
  \Delta\left(\mathbf{u}-\mathbf{q}\right)\right]-\left(\mathbf{u}-\mathbf{q}\right)
  \Lambda^{1-s}\Delta\left(\mathbf{u}-\mathbf{q}\right)
  -\Lambda^{1-s}\left(\mathbf{u}-\mathbf{q}\right)
  \Delta\left(\mathbf{u}-\mathbf{q}\right)\right\|_{L^{2}(\mathbb{R}^{n})}
  \\[0.3cm]
  &\qquad\lesssim \left\| \Lambda^{1-s}\left(\mathbf{u}-\mathbf{q}\right) \right\|_{L^{\frac{n}{s}}(\mathbb{R}^{n})}
 \left\| \Delta\left(\mathbf{u}-\mathbf{q}\right) \right\|_{L^{\frac{2n}{n-2s}}(\mathbb{R}^{n})}
   \\[0.3cm]
  &\qquad\lesssim \left\| \Lambda^{\frac{n}{2}+1-2s}\left(\mathbf{u}-\mathbf{q}\right) \right\|_{L^{2}(\mathbb{R}^{n})}
 \left\| \Lambda^{s}\Delta\left(\mathbf{u}-\mathbf{q}\right) \right\|_{L^{2}(\mathbb{R}^{n})}
 \\[0.3cm]
  &\qquad\lesssim \left\| \nabla\left(\mathbf{u}-\mathbf{q}\right) \right\|_{L^{2}(\mathbb{R}^{n})}
 \left\| \Lambda^{s}\Delta\left(\mathbf{u}-\mathbf{q}\right) \right\|_{L^{2}(\mathbb{R}^{n})},
   \end{array}
\right.
\end{equation}
\begin{equation}\label{5.21}
\left.
\begin{array}{ll}
     \displaystyle \left\| \left(\mathbf{u}-\mathbf{q}\right)
  \Lambda^{1-s}\Delta\left(\mathbf{u}-\mathbf{q}\right)
   \right\|_{L^{2}(\mathbb{R}^{n})}&\lesssim  \left\| \mathbf{u}-\mathbf{q} \right\|_{L^{\frac{n}{2s-1}}(\mathbb{R}^{n})}\left\| \Lambda^{1-s}\Delta\left(\mathbf{u}-\mathbf{q}\right) \right\|_{L^{\frac{2n}{n-2(2s-1)}}(\mathbb{R}^{n})}
   \\[0.3cm]
  &\lesssim\left\| \mathbf{u}-\mathbf{q} \right\|_{L^{\frac{2n}{n-2(n/2+1-2s)}}(\mathbb{R}^{n})}\left\| \Lambda^{s}\Delta\left(\mathbf{u}-\mathbf{q}\right) \right\|_{L^{2}(\mathbb{R}^{n})}
    \\[0.3cm]
  &\lesssim\left\| \Lambda^{n/2+1-2s}\left(\mathbf{u}-\mathbf{q}\right) \right\|_{L^{2}(\mathbb{R}^{n})}\left\| \Lambda^{s}\Delta\left(\mathbf{u}-\mathbf{q}\right) \right\|_{L^{2}(\mathbb{R}^{n})}
  \\[0.3cm]
  &\lesssim\left\| \nabla\left(\mathbf{u}-\mathbf{q}\right) \right\|_{L^{2}(\mathbb{R}^{n})}\left\| \Lambda^{s}\Delta\left(\mathbf{u}-\mathbf{q}\right) \right\|_{L^{2}(\mathbb{R}^{n})},
   \end{array}
\right.
\end{equation}
\begin{equation}\label{5.22}
\left.
\begin{array}{ll}
     \displaystyle \left\| \Lambda^{1-s}\left(\mathbf{u}-\mathbf{q}\right)
   \Delta\left(\mathbf{u}-\mathbf{q}\right)
   \right\|_{L^{2}(\mathbb{R}^{n})}&\lesssim  \left\| \Lambda^{1-s}\left(\mathbf{u}-\mathbf{q}\right) \right\|_{L^{\frac{n}{s}}(\mathbb{R}^{n})}\left\| \Delta\left(\mathbf{u}-\mathbf{q}\right) \right\|_{L^{\frac{2n}{n-2s}}(\mathbb{R}^{n})}
   \\[0.3cm]
  &\lesssim \displaystyle \left\| \nabla\left(\mathbf{u}-\mathbf{q}\right) \right\|_{L^{2}(\mathbb{R}^{n})}\left\| \Lambda^{s}\Delta\left(\mathbf{u}-\mathbf{q}\right) \right\|_{L^{2}(\mathbb{R}^{n})}.
   \end{array}
\right.
\end{equation}
Combining \eqref{5.17}, \eqref{5.18}, \eqref{5.19}, \eqref{5.20}, \eqref{5.21} and \eqref{5.22} gives rise to
\begin{equation}\label{5.23}
\left.
\begin{array}{ll}
     \displaystyle A_{1}&\lesssim \left\| \nabla u
   \right\|^{2}_{L^{2}(\mathbb{R}^{n})} \left\| \nabla\left(\mathbf{u}-\mathbf{q}\right)
     \right\|^{2}_{L^{2}(\mathbb{R}^{n})}+\frac{\nu}{16}\left\|\Lambda^{s} \Delta\left(\mathbf{u}-\mathbf{q}\right) \right\|_{L^{2}(\mathbb{R}^{n})}
     \\[0.3cm]
  &\quad\displaystyle  +\left\| \Lambda^{s} \Delta  u
   \right\|_{L^{2}(\mathbb{R}^{n})} \left\| \nabla\left(\mathbf{u}-\mathbf{q}\right)
       \right\|_{L^{2}(\mathbb{R}^{n})}\left\| \Lambda^{s}\Delta\left(\mathbf{u}-\mathbf{q}\right) \right\|_{L^{2}(\mathbb{R}^{n})}
   \\[0.3cm]
  &\lesssim\displaystyle \left(\left\| \nabla u
   \right\|^{2}_{L^{2}(\mathbb{R}^{n})} +\left\| \Lambda^{s} \Delta u
   \right\|^{2}_{L^{2}(\mathbb{R}^{n})} \right)\left\| \nabla\left(\mathbf{u}-\mathbf{q}\right)
   \right\|^{2}_{L^{2}(\mathbb{R}^{n})}+\frac{\nu}{8}\left\|\Lambda^{s} \Delta\left(\mathbf{u}-\mathbf{q}\right) \right\|^{2}_{L^{2}(\mathbb{R}^{n})}.
   \end{array}
\right.
\end{equation}
In the same manner, recall \eqref{2.3}, the third equation in \eqref{5.1} and Lemma \ref{l2.5}, $A_{2}$, $A_{3}$ and $A_{4}$ can be bounded as follows, respectively.
 \begin{equation}\label{5.24}
\left.
\begin{array}{ll}
A_{2}=\left|\displaystyle\int_{\mathbb{R}^{n}}\mathbf{q}\nabla
\left[\left(\mathbf{u}-\mathbf{q}\right)-\alpha^{2}\Delta\left(\mathbf{u}-\mathbf{q}\right)\right]
  \Delta \left(\mathbf{u}-\mathbf{q}\right)dx\right|
    \\[0.4cm]
 \qquad\lesssim
 \left| \displaystyle\int_{\mathbb{R}^{n}}\mathbf{q}\nabla \left(\mathbf{u}-\mathbf{q}\right)
 \Delta\left(\mathbf{u}-\mathbf{q}\right)dx\right|+\alpha^{2}\left| \displaystyle\int_{\mathbb{R}^{n}}\mathbf{q}\nabla \Delta\left(\mathbf{u}-\mathbf{q}\right)
 \Delta\left(\mathbf{u}-\mathbf{q}\right)dx\right|
  \\[0.4cm]
\qquad\lesssim\left| \displaystyle\int_{\mathbb{R}^{n}}\mathbf{q}\nabla\left(\mathbf{u}-\mathbf{q}\right)
\Delta\left(\mathbf{u}-\mathbf{q}\right)dx\right|
    \\[0.4cm]
 \qquad\lesssim  \left\|\nabla\left(\mathbf{u}-\mathbf{q}\right)\right\|_{L^{2}(\mathbb{R}^{n})}
 \left\|\mathbf{q}\Delta\left(\mathbf{u}-\mathbf{q}\right)\right\|_{L^{2}(\mathbb{R}^{n})}
  \\[0.4cm]
\qquad\lesssim\left\|\nabla \left(\mathbf{u}-\mathbf{q}\right)\right\|_{L^{2}(\mathbb{R}^{n})}
\left\|  \mathbf{q} \right\|_{L^{\frac{n}{s}}(\mathbb{R}^{n})}\left\|  \Delta\left(\mathbf{u}-\mathbf{q}\right)\right\|_{L^{\frac{2n}{n-2s}}(\mathbb{R}^{n})}
 \\[0.4cm]
\qquad\lesssim\left\| \nabla \mathbf{q}
   \right\|^{2}_{L^{2}(\mathbb{R}^{n})}\left\| \nabla\left(\mathbf{u}-\mathbf{q}\right)
   \right\|^{2}_{L^{2}(\mathbb{R}^{n})}+\frac{\nu}{8}\left\|\Lambda^{s} \Delta\left(\mathbf{u}-\mathbf{q}\right) \right\|^{2}_{L^{2}(\mathbb{R}^{n})}.
    \end{array}
\right.
\end{equation}
\begin{equation}\label{5.25}
\left.
\begin{array}{ll}
  A_{3}=\left|\displaystyle\int_{\mathbb{R}^{n}}\left(\mathbf{u}-\alpha^{2}\Delta \mathbf{u}\right) \nabla \left(\mathbf{u}^{T}-\mathbf{q}^{T}\right)\Delta\left(\mathbf{u}-\mathbf{q}\right)dx\right|\\
 \vspace{5pt}
  \qquad\lesssim\left|\displaystyle\int_{\mathbb{R}^{n}} \mathbf{u}  \nabla \left(\mathbf{u}^{T}-\mathbf{q}^{T}\right)\Delta\left(\mathbf{u}-\mathbf{q}\right)dx
  \right|+\alpha^{2}\left|\displaystyle\int_{\mathbb{R}^{n}} \Delta \mathbf{u}  \nabla \left(\mathbf{u}^{T}-\mathbf{q}^{T}\right)\Delta
  \left(\mathbf{u}-\mathbf{q}\right)dx\right|\\
\vspace{5pt}
  \qquad\lesssim\displaystyle\left\|\nabla \left(\mathbf{u}-\mathbf{q}\right) \right\|_{L^{2}(\mathbb{R}^{n})}\left\|\mathbf{u}\Delta \left(\mathbf{u}-\mathbf{q}\right) \right\|_{L^{2}(\mathbb{R}^{n})}+\left\|\nabla \left(\mathbf{u}-\mathbf{q}\right) \right\|_{L^{2}(\mathbb{R}^{n})}\left\|\Delta\mathbf{u}\Delta \left(\mathbf{u}-\mathbf{q}\right) \right\|_{L^{2}(\mathbb{R}^{n})}\\
\vspace{5pt}
  \qquad\lesssim\displaystyle\left\|\nabla \left(\mathbf{u}-\mathbf{q}\right) \right\|_{L^{2}(\mathbb{R}^{n})}\left\| \mathbf{u}  \right\|_{L^{\frac{n}{s}}(\mathbb{R}^{n})}\left\|\Delta \left(\mathbf{u}-\mathbf{q}\right) \right\|_{L^{\frac{2n}{n-2s}}(\mathbb{R}^{n})}\\
  \vspace{5pt}
   \qquad\quad+
  \left\|\nabla \left(\mathbf{u}-\mathbf{q}\right) \right\|_{L^{2}(\mathbb{R}^{n})}\left\| \Delta\mathbf{u}  \right\|_{L^{\frac{n}{s}}(\mathbb{R}^{n})}\left\|\Delta \left(\mathbf{u}-\mathbf{q}\right) \right\|_{L^{\frac{2n}{n-2s}}(\mathbb{R}^{n})}\\
  \vspace{5pt}
    \qquad\lesssim\displaystyle\left\|\Lambda^{s}  \mathbf{u} \right\|^{2}_{L^{2}(\mathbb{R}^{n})}\left\|\nabla \left(\mathbf{u}-\mathbf{q}\right) \right\|^{2}_{L^{2}(\mathbb{R}^{n})}+\frac{\nu}{16}\left\|\Lambda^{s} \Delta\left(\mathbf{u}-\mathbf{q}\right) \right\|^{2}_{L^{2}(\mathbb{R}^{n})}\\
  \vspace{5pt}
   \qquad\quad +\displaystyle\left\|\nabla \left(\mathbf{u}-\mathbf{q}\right) \right\|^{2}_{L^{2}(\mathbb{R}^{n})}\left\|\Lambda^{s} \Delta \mathbf{u}  \right\|^{2}_{L^{2}(\mathbb{R}^{n})} +\frac{\nu}{16}\left\|\Lambda^{s} \Delta\left(\mathbf{u}-\mathbf{q}\right) \right\|^{2}_{L^{2}(\mathbb{R}^{n})}\\
  \qquad\lesssim\displaystyle\left(\left\|\Lambda^{s} \mathbf{u} \right\|^{2}_{L^{2}(\mathbb{R}^{n})}+\left\|\Lambda^{s} \Delta \mathbf{u}  \right\|^{2}_{L^{2}(\mathbb{R}^{n})}\right)\left\|\nabla \left(\mathbf{u}-\mathbf{q}\right) \right\|^{2}_{L^{2}(\mathbb{R}^{n})}  +\frac{\nu}{8}\left\|\Lambda^{s} \Delta\left(\mathbf{u}-\mathbf{q}\right) \right\|^{2}_{L^{2}(\mathbb{R}^{n})}.
      \end{array}
\right.
\end{equation}
In view of interpolation inequality, we further make the similar argument for $A_{4}$ to those employed in $A_{1}$:

\begin{equation}\label{5.26}
\left.
\begin{array}{ll}
\vspace{5pt}
A_{4}=\left|\displaystyle\int_{\mathbb{R}{n}}\left[\left(\mathbf{u}-\mathbf{q}\right)-\alpha^{2}\Delta \left(\mathbf{u}-\mathbf{q}\right)\right] \nabla \mathbf{q}^{T} \Delta\left(\mathbf{u}-\mathbf{q}\right)dx\right|\\
  \vspace{5pt}
  \qquad\lesssim\displaystyle\left| \int_{\mathbb{R}{n}} \left(\mathbf{u}-\mathbf{q}\right)  \nabla \mathbf{q}^{T} \Delta\left(\mathbf{u}-\mathbf{q}\right)dx\right|+\left|\int_{\mathbb{R}{n}}
  \Delta \left(\mathbf{u}-\mathbf{q}\right)  \nabla \mathbf{q}^{T} \Delta\left(\mathbf{u}-\mathbf{q}\right)dx\right|\\
  \vspace{5pt}
  \qquad\lesssim\displaystyle  \left\|\Lambda^{s}\mathbf{q}\right\|_{L^{2}(\mathbb{R}^{n})}\left\| \nabla\left(\mathbf{u}-\mathbf{q}\right) \right\|_{L^{2}(\mathbb{R}^{n})}
 \left\|\Lambda^{s}\Delta\left(\mathbf{u}-\mathbf{q}\right) \right\|_{L^{2}(\mathbb{R}^{n})}+\left\|\nabla\mathbf{q}\right\|_{L^{2}(\mathbb{R}^{n})}
    \left\| \Delta\left(\mathbf{u}-\mathbf{q}\right) \right\|^{2}_{L^{4}(\mathbb{R}^{n})}\\
  \vspace{5pt}
  \qquad\lesssim\displaystyle   \left\|\Lambda^{s}\mathbf{q}\right\|_{L^{2}(\mathbb{R}^{n})}\left\| \nabla\left(\mathbf{u}-\mathbf{q}\right) \right\|_{L^{2}(\mathbb{R}^{n})}
 \left\|\Lambda^{s}\Delta\left(\mathbf{u}-\mathbf{q}\right) \right\|_{L^{2}(\mathbb{R}^{n})}\\
  \vspace{5pt}
  \qquad\quad\displaystyle
 +\left\|\nabla\mathbf{q}\right\|_{L^{2}(\mathbb{R}^{n})}\left(C(\varepsilon) \left\| \Delta\left(\mathbf{u}-\mathbf{q}\right) \right\|^{2}_{L^{2}(\mathbb{R}^{n})}+\varepsilon \left\|\Lambda^{s}\Delta\left(\mathbf{u}-\mathbf{q}\right) \right\|^{2}_{L^{2}(\mathbb{R}^{n})}\right)\\
  \vspace{5pt}
  \qquad\lesssim\displaystyle   \left\|\Lambda^{s}\mathbf{q}\right\|^{2}_{L^{2}(\mathbb{R}^{n})}\left\| \nabla\left(\mathbf{u}-\mathbf{q}\right) \right\|^{2}_{L^{2}(\mathbb{R}^{n})}
  +\frac{\nu}{16}\left\|\Lambda^{s} \Delta\left(\mathbf{u}-\mathbf{q}\right) \right\|^{2}_{L^{2}(\mathbb{R}^{n})}\\
  \vspace{5pt}
  \qquad\quad+\displaystyle C(\varepsilon)\left\|\nabla\mathbf{q}\right\| _{L^{2}(\mathbb{R}^{n})}\left\| \Delta\left(\mathbf{u}-\mathbf{q}\right) \right\|^{2}_{L^{2}(\mathbb{R}^{n})}
 +\varepsilon\left\|\nabla\mathbf{q}\right\| _{L^{2}(\mathbb{R}^{n})}
  \left\|\Lambda^{s}\Delta\left(\mathbf{u}-\mathbf{q}\right) \right\|^{2}_{L^{2}(\mathbb{R}^{n})} \\
  \vspace{5pt}
  \qquad\lesssim\displaystyle   \left(\left\|\Lambda^{s}\mathbf{q}\right\|^{2}_{L^{2}(\mathbb{R}^{n})}
  +\left\|\Lambda^{s}\nabla\mathbf{q}\right\| _{L^{2}(\mathbb{R}^{n})}\right)\left(\left\| \nabla\left(\mathbf{u}-\mathbf{q}\right) \right\|^{2}_{L^{2}(\mathbb{R}^{n})}+\alpha^{2}\left\| \Delta\left(\mathbf{u}-\mathbf{q}\right) \right\|^{2}_{L^{2}(\mathbb{R}^{n})}\right)\\
  \vspace{5pt}
  \qquad\quad+\displaystyle   \frac{\nu}{16}\left\|\Lambda^{s} \Delta\left(\mathbf{u}-\mathbf{q}\right) \right\|^{2}_{L^{2}(\mathbb{R}^{n})}+\varepsilon\left\|\nabla\mathbf{q}\right\| _{L^{2}(\mathbb{R}^{n})}\left\|\Lambda^{s} \Delta\left(\mathbf{u}-\mathbf{q}\right) \right\|^{2}_{L^{2}(\mathbb{R}^{n})}.
   \end{array}
\right.
\end{equation}
In view of \eqref{2.3}, choosing $\varepsilon $ sufficiently small such that $\displaystyle\varepsilon\left\|\nabla\mathbf{q}\right\|_{L^{2}(\mathbb{R}^{n})}\leq  \frac{\nu}{16}$, \eqref{5.26} then leads to
\begin{equation}\label{5.27}
\left.
\begin{array}{ll}
A_{4}  \lesssim\displaystyle   \left\|\nabla\mathbf{q}\right\|^{2}_{L^{2}(\mathbb{R}^{n})}\left(\left\| \nabla\left(\mathbf{u}-\mathbf{q}\right) \right\|^{2}_{L^{2}(\mathbb{R}^{n})}+\alpha^{2}\left\| \Delta\left(\mathbf{u}-\mathbf{q}\right) \right\|^{2}_{L^{2}(\mathbb{R}^{n})}\right)+\frac{\nu}{8}\left\|\Lambda^{s} \Delta\left(\mathbf{u}-\mathbf{q}\right) \right\|^{2}_{L^{2}(\mathbb{R}^{n})}.
   \end{array}
\right.
\end{equation}
Recall \eqref{1.9}, \eqref{2.3} and \eqref{2.5}, combining \eqref{5.16} with \eqref{5.17}-\eqref{5.27} yields that for $\displaystyle \frac{n}{4}< s<1$ with $n=2,3$,
\begin{equation}\label{5.28}
\left.
\begin{array}{ll}
\vspace{5pt}
\displaystyle\frac{1}{2}\frac{d}{dt} \left(\left\| \nabla\left(\mathbf{u}-\mathbf{q}\right)\right\|^{2}_{L^{2}(\mathbb{R}^{n})}+\alpha^{2}\left\| \Delta\left(\mathbf{u}-\mathbf{q}\right)\right\|^{2}_{L^{2}(\mathbb{R}^{n})}\right)\\
  \vspace{5pt}
  \qquad \quad \displaystyle+\frac{\nu}{2}\left(\left\| \Lambda^{s}\nabla \left(\mathbf{u}-\mathbf{q}\right)\right\|^{2}_{L^{2}(\mathbb{R}^{n})}+\alpha^{2}\left\| \Lambda^{s} \Delta\left (\mathbf{u}-\mathbf{q}\right)\right\|^{2}_{L^{2}(\mathbb{R}^{n})}\right)\\
 \vspace{5pt}
 \qquad\lesssim c^{*}\left(\left\|\nabla\left( \mathbf{u}-\mathbf{q}\right ) \right \|^{2}_{L^{2}(\mathbb{R}^{n})}+\alpha^{2}\left\|\Delta\left( \mathbf{u}-\mathbf{q}\right ) \right \|^{2}_{L^{2}(\mathbb{R}^{n})}\right),
 \end{array}
\right.
\end{equation}
where
$$c^{*}=\left\| \Lambda^{s}   \mathbf{q} \right\|^{2}_{L^{2}(\mathbb{R}^{n})}+\left\| \nabla \mathbf{q} \right\|^{2}_{L^{2}(\mathbb{R}^{n})}+\left\| \nabla \mathbf{q} \right\| _{L^{2}(\mathbb{R}^{n})}+
\left\| \Lambda^{s}   \mathbf{u} \right\|^{2}_{L^{2}(\mathbb{R}^{n})}+\left\| \nabla \mathbf{u} \right\|^{2}_{L^{2}(\mathbb{R}^{n})}+\left\| \Lambda^{s}  \Delta \mathbf{u} \right\| _{L^{2}(\mathbb{R}^{n})}.
$$
We then handle  {\bf Case 2} ~~ $\displaystyle s=\frac{n}{4}$ for $n=2,3$.
\\[0.3cm]
\indent In this case, note that \eqref{5.16}, with the help of Lemma \ref{l2.5}, we have the following a priori estimates for $A_{i}~(i=1,2,3,4)$.
\begin{equation}\label{5.29}
\left.
\begin{array}{ll}
\vspace{5pt}
A_{1}=\left|\displaystyle\int_{\mathbb{R}^{n}}\left(\mathbf{u}-\mathbf{q}\right) \nabla\left(\mathbf{u}-\alpha^{2}\Delta \mathbf{u}\right) \Delta\left(\mathbf{u}-\mathbf{q}\right)dx\right|\\
 \vspace{5pt}
 \qquad\lesssim  \left|\displaystyle\int_{\mathbb{R}^{n}}\left(\mathbf{u}-\mathbf{q}\right) \nabla \mathbf{u}  \Delta\left(\mathbf{u}-\mathbf{q}\right)dx\right|+\left|\displaystyle\int_{\mathbb{R}^{n}}
 \left(\mathbf{u}-\mathbf{q}\right) \nabla \Delta \mathbf{u} \Delta\left(\mathbf{u}-\mathbf{q}\right)dx\right|.
\end{array}
\right.
\end{equation}
A straightforward calculation gives
\begin{equation}\label{5.30}
  \left|\displaystyle\int_{\mathbb{R}^{n}}\left(\mathbf{u}-\mathbf{q}\right) \nabla \mathbf{u}  \Delta\left(\mathbf{u}-\mathbf{q}\right)dx\right| \lesssim \left\|\Lambda^{\frac{n}{4}} \mathbf{u} \right\|_{L^{2}(\mathbb{R}^{n})}\left\|\Lambda^{1-\frac{n}{4}} \left[\left(\mathbf{u}-\mathbf{q}\right)\Delta\left(\mathbf{u}-\mathbf{q}\right)\right] \right\|_{L^{2}(\mathbb{R}^{n})}.
\end{equation}
Indeed, $\displaystyle \left\|\Lambda^{1-\frac{n}{4}} \left[\left(\mathbf{u}-\mathbf{q}\right)\Delta\left(\mathbf{u}-\mathbf{q}\right)\right] \right\|_{L^{2}(\mathbb{R}^{n})}$ can be bounded as
 \begin{equation}\label{5.31}
\left.
\begin{array}{ll}
\vspace{5pt}
\displaystyle  \left\|\Lambda^{1-\frac{n}{4}} \left[\left(\mathbf{u}-\mathbf{q}\right)\Delta\left(\mathbf{u}-\mathbf{q}\right)\right] \right\|_{L^{2}(\mathbb{R}^{n})}\\
\vspace{5pt}
 \qquad\lesssim\displaystyle  \left\|\Lambda^{1-\frac{n}{4}} \left[\left(\mathbf{u}-\mathbf{q}\right)\Delta\left(\mathbf{u}-\mathbf{q}\right)\right]
\right.
 \\
 \vspace{5pt}
 \qquad\qquad\left. -
 \left(\mathbf{u}-\mathbf{q}\right)\Lambda^{1-\frac{n}{4}}\Delta
 \left(\mathbf{u}-\mathbf{q}\right)
 -\Lambda^{1-\frac{n}{4}}\left(\mathbf{u}-\mathbf{q}\right)\Delta\left(\mathbf{u}-\mathbf{q}\right) \right\|_{L^{2}(\mathbb{R}^{n})}\\
\vspace{5pt}
 \qquad\quad \displaystyle +\left\| \left(\mathbf{u}-\mathbf{q}\right)\Lambda^{1-\frac{n}{4}}\Delta\left(\mathbf{u}-\mathbf{q}\right)
  \right\|_{L^{2}(\mathbb{R}^{n})}+\left\| \Lambda^{1-\frac{n}{4}}\left(\mathbf{u}-\mathbf{q}\right)\Delta\left(\mathbf{u}-\mathbf{q}\right) \right\|_{L^{2}(\mathbb{R}^{n})}.
\end{array}
\right.
\end{equation}
  Thanks to the Gagliardo-Nirenberg-Sobolev inequality and H\"{o}lder's inequality, by virtue of (I) of Lemma \ref{l2.7}, for $\displaystyle 0<s_{1}<1-\frac{n}{4}$, one gets
  \begin{equation}\label{5.32}
\left.
\begin{array}{ll}
\vspace{5pt}
\displaystyle   \left\|\Lambda^{1-\frac{n}{4}} \left[\left(\mathbf{u}-\mathbf{q}\right)\Delta\left(\mathbf{u}-\mathbf{q}\right)\right] \right.
 \\
 \vspace{5pt}
 \qquad\qquad\left.-
 \left(\mathbf{u}-\mathbf{q}\right)\Lambda^{1-\frac{n}{4}}\Delta
 \left(\mathbf{u}-\mathbf{q}\right)
 -\Lambda^{1-\frac{n}{4}}\left(\mathbf{u}-\mathbf{q}\right)\Delta\left(\mathbf{u}-\mathbf{q}\right) \right\|_{L^{2}(\mathbb{R}^{n})}\\
\vspace{5pt}
 \qquad\lesssim\displaystyle  \left\|\Lambda^{1-\frac{n}{4}-s_{1}}\left(\mathbf{u}-\mathbf{q}\right)
 \right\|_{L^{\frac{2n}{n-2(n/4+s_{1})}}(\mathbb{R}^{n})}
 \left\|\Lambda^{ s_{1}}\Delta\left(\mathbf{u}-\mathbf{q}\right)
 \right\|_{L^{\frac{2n}{n-2(n/4-s_{1})}}(\mathbb{R}^{n})}\\
\vspace{5pt}
 \qquad\lesssim\displaystyle  \left\|\nabla \left(\mathbf{u}-\mathbf{q}\right)\right\|_{L^{2}(\mathbb{R}^{n})}\left\|\Lambda^{ \frac{n}{4}}\Delta \left(\mathbf{u}-\mathbf{q}\right)\right\|_{L^{2}(\mathbb{R}^{n})},
 \end{array}
\right.
\end{equation}
 where $\displaystyle s_{1}\in \left(0,\frac{1}{2}\right)$,  $\displaystyle 1-\frac{n}{4}-s_{1}\in \left(0,\frac{1}{2}\right)$,  $\displaystyle \frac{1}{2}=\frac{1}{p_{1}}+\frac{1}{p_{2}}  $ with $p_{1}, p_{2}\in (1,\infty)$, $\displaystyle p_{1}=\frac{2n}{n-2(n/4+s_{1})}$,  $\displaystyle p_{2}=\frac{2n}{n-2(n/4-s_{1})}$.\\
 \indent In the same manner, in view of \eqref{2.3} and \eqref{5.15}, Lemma \ref{l2.10}, Agmon's inequality and interpolation inequality, the second  and third terms on the right hand side of inequality \eqref{5.31} can be bounded as follows:
  \begin{equation}\label{5.33}
\left.
\begin{array}{ll}
\vspace{5pt} \left\| \left(\mathbf{u}-\mathbf{q}\right)\Lambda^{1-\frac{n}{4}}
\Delta\left(\mathbf{u}-\mathbf{q}\right)
  \right\|_{L^{2}(\mathbb{R}^{n})}+\left\| \Lambda^{1-\frac{n}{4}}\left(\mathbf{u}-\mathbf{q}\right)\Delta\left(\mathbf{u}-\mathbf{q}\right) \right\|_{L^{2}(\mathbb{R}^{n})}\\
  \vspace{5pt}
 \qquad\lesssim\displaystyle\left\| \mathbf{u}-\mathbf{q}
  \right\|_{L^{\infty}(\mathbb{R}^{n})}\left\| \Lambda^{1-\frac{n}{4}}
\Delta\left(\mathbf{u}-\mathbf{q}\right)
  \right\|_{L^{2}(\mathbb{R}^{n})}+\left\| \Lambda^{1-\frac{n}{4}}\left(\mathbf{u}-\mathbf{q}\right) \right\|_{L^{\infty}(\mathbb{R}^{n})}\left\| \Delta\left(\mathbf{u}-\mathbf{q}\right) \right\|_{L^{2}(\mathbb{R}^{n})}\\
  \vspace{5pt}
 \qquad\lesssim\displaystyle\left\| \mathbf{u}-\mathbf{q}  \right\|^{\frac{1}{2}}_{H^{1}(\mathbb{R}^{n})}\left\| \mathbf{u}-\mathbf{q}  \right\|^{\frac{1}{2}}_{H^{2}(\mathbb{R}^{n})}
 \left\|
\Delta\left(\mathbf{u}-\mathbf{q}\right)
  \right\|^\frac{1}{2}_{L^{2}(\mathbb{R}^{n})}
 \left\| \Lambda^{\frac{n}{4}}
\Delta\left(\mathbf{u}-\mathbf{q}\right)
  \right\|^\frac{1}{2}_{L^{2}(\mathbb{R}^{n})}\\
  \vspace{5pt}
 \qquad\quad+\displaystyle
 \displaystyle\left\| \Lambda^{1-\frac{n}{4}}
 \left(\mathbf{u}-\mathbf{q}\right)  \right\|^{\frac{1}{2}}_{H^{1}(\mathbb{R}^{n})}\left\| \Lambda^{1-\frac{n}{4}}
 \left(\mathbf{u}-\mathbf{q}\right)\right\|^{\frac{1}{2}}_{H^{2}(\mathbb{R}^{n})}
 \left\|
\Delta\left(\mathbf{u}-\mathbf{q}\right)
  \right\| _{L^{2}(\mathbb{R}^{n})} \\
  \vspace{5pt}
 \qquad\lesssim\displaystyle\left\|\nabla\left(\mathbf{u}-\mathbf{q}\right)
  \right\|^{\frac{1}{2}}_{L^{2}(\mathbb{R}^{n})}\left\|\Delta\left(\mathbf{u}-\mathbf{q}\right)
  \right\| _{L^{2}(\mathbb{R}^{n})}\left\| \Lambda^{ \frac{n}{4}}
 \Delta\left(\mathbf{u}-\mathbf{q}\right)
  \right\|^{\frac{1}{2}}_{L^{2}(\mathbb{R}^{n})}\\
  \vspace{5pt}
 \qquad\quad\displaystyle
  +\left\|\nabla\left(\mathbf{u}-\mathbf{q}\right)
  \right\|^{\frac{1}{4}} _{L^{2}(\mathbb{R}^{n})} \left\| \Lambda^{ \frac{n}{4}}
 \nabla\left(\mathbf{u}-\mathbf{q}\right)
  \right\|^{\frac{1}{4}}_{L^{2}(\mathbb{R}^{n})}
  \left\|\Delta\left(\mathbf{u}-\mathbf{q}\right)
  \right\|^{\frac{1}{4}} _{L^{2}(\mathbb{R}^{n})}\\
  \vspace{5pt}
 \qquad\qquad\qquad\displaystyle \times\left\| \Lambda^{ \frac{n}{4}}
 \Delta\left(\mathbf{u}-\mathbf{q}\right)
  \right\|^{\frac{1}{4}}_{L^{2}(\mathbb{R}^{n})} \left\|
 \Delta\left(\mathbf{u}-\mathbf{q}\right)
  \right\| _{L^{2}(\mathbb{R}^{n})}.
\end{array}
\right.
\end{equation}
By the aid of H\"{o}lder's inequality, combining \eqref{5.29} with \eqref{5.30}, \eqref{5.31}, \eqref{5.32} and \eqref{5.33} leads to
  \begin{equation}\label{5.34}
\left.
\begin{array}{ll}
\vspace{5pt}
\left|\displaystyle\int_{\mathbb{R}^{n}} \left(\mathbf{u}-\mathbf{q}\right)\nabla \mathbf{u}\Delta
  \left(\mathbf{u}-\mathbf{q}\right)
  dx\right|\\
  \vspace{5pt}
  \qquad\lesssim\displaystyle \left\|\Lambda^{ \frac{n}{4}} \mathbf{u}
  \right\|  _{L^{2}(\mathbb{R}^{n})}\left\|\nabla\left(\mathbf{u}-\mathbf{q}\right)
  \right\| _{L^{2}(\mathbb{R}^{n})}\left\|\Lambda^{ \frac{n}{4}}\Delta\left(\mathbf{u}-\mathbf{q}\right)
  \right\| _{L^{2}(\mathbb{R}^{n})}\\
  \vspace{5pt}
  \qquad\quad +\displaystyle \left\|\Lambda^{ \frac{n}{4}} \mathbf{u}
  \right\|  _{L^{2}(\mathbb{R}^{n})}\left\|\nabla\left(\mathbf{u}-\mathbf{q}\right)
  \right\|^{\frac{1}{2}} _{L^{2}(\mathbb{R}^{n})}\left\| \Delta\left(\mathbf{u}-\mathbf{q}\right)
  \right\| _{L^{2}(\mathbb{R}^{n})}\left\|\Lambda^{ \frac{n}{4}}\Delta\left(\mathbf{u}-\mathbf{q}\right)
  \right\| ^{\frac{1}{2}}_{L^{2}(\mathbb{R}^{n})}\\
  \vspace{5pt}
  \qquad\quad +\displaystyle \left\|\Lambda^{ \frac{n}{4}} \mathbf{u}
  \right\|  _{L^{2}(\mathbb{R}^{n})}\left\|\nabla\left(\mathbf{u}-\mathbf{q}\right)
  \right\|^{\frac{1}{4}} _{L^{2}(\mathbb{R}^{n})}\left\|\Lambda^{ \frac{n}{4}}\nabla\left(\mathbf{u}-\mathbf{q}\right)
  \right\|^{\frac{1}{4}} _{L^{2}(\mathbb{R}^{n})}\\
  \vspace{5pt}
  \qquad\quad +\left\|\Delta\left(\mathbf{u}-\mathbf{q}\right)
  \right\|^{\frac{1}{4}} _{L^{2}(\mathbb{R}^{n})}\left\|\Lambda^{ \frac{n}{4}}\Delta\left(\mathbf{u}-\mathbf{q}\right)
  \right\| ^{\frac{1}{4}}_{L^{2}(\mathbb{R}^{n})}
  \left\| \Delta\left(\mathbf{u}-\mathbf{q}\right)
  \right\| _{L^{2}(\mathbb{R}^{n})}\\
  \vspace{5pt}
   \qquad\lesssim\displaystyle \left(\left\|\nabla\left(\mathbf{u}-\mathbf{q}\right)
  \right\|^{2} _{L^{2}(\mathbb{R}^{n})}+\alpha^{2}\left\|\Delta\left(\mathbf{u}-\mathbf{q}\right)
  \right\|^{2} _{L^{2}(\mathbb{R}^{n})}\right)\\
  \vspace{5pt}
 \qquad\qquad\displaystyle\cdot\left\|\Lambda^{ \frac{n}{4}} \mathbf{u}
  \right\|^{2} _{L^{2}(\mathbb{R}^{n})}\left(\left\|\Lambda^{ \frac{n}{4}}\nabla\left(\mathbf{u}-\mathbf{q}\right)
  \right\|^{2} _{L^{2}(\mathbb{R}^{n})}+\left\|\Lambda^{ \frac{n}{4}}\Delta\left(\mathbf{u}-\mathbf{q}\right)
  \right\|^{2} _{L^{2}(\mathbb{R}^{n})}\right).
\end{array}
\right.
\end{equation}
Similarly, note that
$$\left\|\Delta\left(\mathbf{u}-\mathbf{q}\right)
  \right\| ^{\frac{1}{2}}_{L^{2}(\mathbb{R}^{n})}\lesssim \left\|\Lambda^{ \frac{n}{4}} \nabla\left(\mathbf{u}-\mathbf{q}\right)
  \right\|^{\frac{1}{4}}  _{L^{2}(\mathbb{R}^{n})}\left\|\Lambda^{ \frac{n}{4}} \Delta\left(\mathbf{u}-\mathbf{q}\right)
  \right\|^{\frac{1}{4}}  _{L^{2}(\mathbb{R}^{n})},$$
using H\"{o}lder's inequality repeatedly, the second term on the right hand side of \eqref{5.29} can be bounded as follows:
 \begin{equation}\label{5.35}
\left.
\begin{array}{ll}
\vspace{5pt}
\left|\displaystyle\int_{\mathbb{R}^{n}} \left(\mathbf{u}-\mathbf{q}\right)\nabla \Delta \mathbf{u}\Delta
  \left(\mathbf{u}-\mathbf{q}\right)
  dx\right|\\
  \vspace{5pt}
    \qquad\lesssim\displaystyle \left\|\Lambda^{ \frac{n}{4}} \Delta\mathbf{u}
  \right\|  _{L^{2}(\mathbb{R}^{n})}\left\|\nabla\left(\mathbf{u}-\mathbf{q}\right)
  \right\| _{L^{2}(\mathbb{R}^{n})}\left\|\Lambda^{ \frac{n}{4}}\Delta\left(\mathbf{u}-\mathbf{q}\right)
  \right\| _{L^{2}(\mathbb{R}^{n})}\\
  \vspace{5pt}
  \qquad\quad +\displaystyle \left\|\Lambda^{ \frac{n}{4}} \Delta\mathbf{u}
  \right\|  _{L^{2}(\mathbb{R}^{n})}\left\|\nabla\left(\mathbf{u}-\mathbf{q}\right)
  \right\|^{\frac{1}{2}} _{L^{2}(\mathbb{R}^{n})}\left\| \Delta\left(\mathbf{u}-\mathbf{q}\right)
  \right\| _{L^{2}(\mathbb{R}^{n})}\left\|\Lambda^{ \frac{n}{4}}\Delta\left(\mathbf{u}-\mathbf{q}\right)
  \right\| ^{\frac{1}{2}}_{L^{2}(\mathbb{R}^{n})}\\
  \vspace{5pt}
  \qquad\quad +\displaystyle \left\|\Lambda^{ \frac{n}{4}}\Delta \mathbf{u}
  \right\|  _{L^{2}(\mathbb{R}^{n})}\left\|\nabla\left(\mathbf{u}-\mathbf{q}\right)
  \right\|^{\frac{1}{4}} _{L^{2}(\mathbb{R}^{n})}\left\|\Lambda^{ \frac{n}{4}}\nabla\left(\mathbf{u}-\mathbf{q}\right)
  \right\|^{\frac{1}{4}} _{L^{2}(\mathbb{R}^{n})}\\
  \vspace{5pt}
  \qquad\qquad \cdot\left\|\Delta\left(\mathbf{u}-\mathbf{q}\right)
  \right\|^{\frac{1}{4}} _{L^{2}(\mathbb{R}^{n})}\left\|\Lambda^{ \frac{n}{4}}\Delta\left(\mathbf{u}-\mathbf{q}\right)
  \right\| ^{\frac{1}{4}}_{L^{2}(\mathbb{R}^{n})}
  \left\| \Delta\left(\mathbf{u}-\mathbf{q}\right)
  \right\| _{L^{2}(\mathbb{R}^{n})}\\
  \vspace{5pt}
    \qquad\lesssim\displaystyle \left\|\Lambda^{ \frac{n}{4}} \Delta\mathbf{u}
  \right\|^{2}  _{L^{2}(\mathbb{R}^{n})}\left\|\nabla\left(\mathbf{u}-\mathbf{q}\right)
  \right\|^{2} _{L^{2}(\mathbb{R}^{n})}+\frac{\nu}{32}\left\|\Lambda^{ \frac{n}{4}}\Delta\left(\mathbf{u}-\mathbf{q}\right)
  \right\|^{2} _{L^{2}(\mathbb{R}^{n})}\\
  \vspace{5pt}
    \qquad\quad  \displaystyle +\left\|\nabla\left(\mathbf{u}-\mathbf{q}\right)
  \right\| _{L^{2}(\mathbb{R}^{n})}\left\|\Lambda^{ \frac{n}{4}}\Delta\left(\mathbf{u}-\mathbf{q}\right)
  \right\| _{L^{2}(\mathbb{R}^{n})}+\left\|\Lambda^{ \frac{n}{4}} \Delta\mathbf{u}
  \right\|^{2}  _{L^{2}(\mathbb{R}^{n})}\left\|\Delta\left(\mathbf{u}-\mathbf{q}\right)
  \right\|^{2} _{L^{2}(\mathbb{R}^{n})}\\
  \vspace{5pt}
 \qquad\quad\displaystyle + \left\|\Lambda^{ \frac{n}{4}} \Delta\mathbf{u}
  \right\|  _{L^{2}(\mathbb{R}^{n})}\left\|\nabla\left(\mathbf{u}-\mathbf{q}\right)
  \right\|^{\frac{1}{2}} _{L^{2}(\mathbb{R}^{n})}\left\|\Lambda^{ \frac{n}{4}}\nabla\left(\mathbf{u}-\mathbf{q}\right)
  \right\| _{L^{2}(\mathbb{R}^{n})}\left\|\Delta\left(\mathbf{u}-\mathbf{q}\right)
  \right\| ^{\frac{1}{2}}_{L^{2}(\mathbb{R}^{n})}\\
  \vspace{5pt}
 \qquad\quad\displaystyle +\left\|\Lambda^{ \frac{n}{4}} \Delta\mathbf{u}
  \right\|  _{L^{2}(\mathbb{R}^{n})}\left\|\Lambda^{ \frac{n}{4}} \Delta\left(\mathbf{u}-\mathbf{q}\right)
  \right\| _{L^{2}(\mathbb{R}^{n})}\left\|\Delta\left(\mathbf{u}-\mathbf{q}\right)
  \right\| _{L^{2}(\mathbb{R}^{n})}
    \\
  \vspace{5pt}
 \qquad\lesssim\displaystyle \left(1+\left\|\Lambda^{ \frac{n}{4}} \Delta\mathbf{u}
  \right\|^{2}  _{L^{2}(\mathbb{R}^{n})}\right)\left(\left\|\nabla\left(\mathbf{u}-\mathbf{q}
  \right)\right\|^{2} _{L^{2}(\mathbb{R}^{n})}+\alpha^{2}\left\|\Delta\left(\mathbf{u}-\mathbf{q}\right)
  \right\|^{2} _{L^{2}(\mathbb{R}^{n})}\right)\\
  \vspace{5pt}
 \qquad\quad\displaystyle+\frac{\nu}{8}\left(\left\|\Lambda^{ \frac{n}{4}}\nabla\left(\mathbf{u}-\mathbf{q}\right)
  \right\|^{2} _{L^{2}(\mathbb{R}^{n})}+\left\|\Lambda^{ \frac{n}{4}}\Delta\left(\mathbf{u}-\mathbf{q}\right)
  \right\|^{2} _{L^{2}(\mathbb{R}^{n})}\right).
\end{array}
\right.
\end{equation}
Collecting \eqref{5.29}, \eqref{5.34} and \eqref{5.35} together leads to
 \begin{equation}\label{5.36}
\left.
\begin{array}{ll}
\vspace{5pt}
 A_{1}\lesssim\displaystyle\left(1+\left\|\Lambda^{ \frac{n}{4}} \Delta\mathbf{u}
  \right\|^{2}  _{L^{2}(\mathbb{R}^{n})}\right) \left(\left\|\nabla\left(\mathbf{u}-\mathbf{q}\right)
  \right\|^{2} _{L^{2}(\mathbb{R}^{n})}+\alpha^{2}\left\|\Delta\left(\mathbf{u}-\mathbf{q}\right)
  \right\|^{2} _{L^{2}(\mathbb{R}^{n})}\right)\\
  \vspace{5pt}
 \qquad\displaystyle+ \left(\left\|\Lambda^{ \frac{n}{4}} \mathbf{u}
  \right\|^{2} _{L^{2}(\mathbb{R}^{n})}+\frac{\nu}{8}\right)\left(\left\|\Lambda^{ \frac{n}{4}}\nabla\left(\mathbf{u}-\mathbf{q}\right)
  \right\|^{2} _{L^{2}(\mathbb{R}^{n})}+\left\|\Lambda^{ \frac{n}{4}}\Delta\left(\mathbf{u}-\mathbf{q}\right)
  \right\|^{2} _{L^{2}(\mathbb{R}^{n})}\right).
\end{array}
\right.
\end{equation}
Making the similar a priori estimates to those employed in \eqref{5.29}-\eqref{5.36}, note that \eqref{5.16} and the third equation in \eqref{5.1}, we can conclude the bounds of $A_{2}$, $A_{3}$ and $A_{4}$ in \eqref{5.16} for $\displaystyle s=\frac{n}{4}$:
\begin{equation}\label{5.37}
\left.
\begin{array}{ll}
\vspace{5pt}
A_{2}&=\left|\displaystyle\int_{\mathbb{R}^{n}} \mathbf{q}\nabla\left[\left(\mathbf{u}-\mathbf{q}\right)-\alpha^{2}\Delta
  \left(\mathbf{u}-\mathbf{q}\right)\right]\Delta
  \left(\mathbf{u}-\mathbf{q}\right)
  dx\right|\\
  \vspace{5pt}
 &=\left|\displaystyle\int_{\mathbb{R}^{n}} \mathbf{q}\nabla \left(\mathbf{u}-\mathbf{q}\right) \Delta
  \left(\mathbf{u}-\mathbf{q}\right)
  dx\right| \\
  \vspace{5pt}
&\lesssim\displaystyle\left\| \mathbf{q}
  \right\| _{L^{2}(\mathbb{R}^{n})}\left\| \nabla \left(\mathbf{u}-\mathbf{q}\right) \Delta
  \left(\mathbf{u}-\mathbf{q}\right)
  \right\| _{L^{2}(\mathbb{R}^{n})} \\
  \vspace{5pt}
&\lesssim\displaystyle\left\| \mathbf{q}
  \right\| _{L^{2}(\mathbb{R}^{n})}\left\| \nabla \left(\mathbf{u}-\mathbf{q}\right)
  \right\| _{L^{4}(\mathbb{R}^{n})}\left\| \Delta \left(\mathbf{u}-\mathbf{q}\right)
  \right\| _{L^{4}(\mathbb{R}^{n})} \\
  \vspace{5pt}
&\lesssim\displaystyle\left\| \mathbf{q}
  \right\| _{L^{2}(\mathbb{R}^{n})}\left\| \nabla \left(\mathbf{u}-\mathbf{q}\right)
  \right\|^{2} _{L^{4}(\mathbb{R}^{n})}+\left\| \mathbf{q}
  \right\| _{L^{2}(\mathbb{R}^{n})}\left\| \Delta \left(\mathbf{u}-\mathbf{q}\right)
  \right\|^{2} _{L^{4}(\mathbb{R}^{n})}\\
  \vspace{5pt}
&\lesssim\displaystyle\left\| \mathbf{q}
  \right\| _{L^{2}(\mathbb{R}^{n})}\left(\left\| \nabla \left(\mathbf{u}-\mathbf{q}\right)
  \right\|^{2} _{L^{2}(\mathbb{R}^{n})}+\left\|\Lambda^{\frac{n}{4}} \nabla \left(\mathbf{u}-\mathbf{q}\right)
  \right\|^{2} _{L^{2}(\mathbb{R}^{n})}\right)\\
  \vspace{5pt}
&\quad\displaystyle+\left\| \mathbf{q}
  \right\| _{L^{2}(\mathbb{R}^{n})}\left(\left\| \Delta\left(\mathbf{u}-\mathbf{q}\right)
  \right\|^{2} _{L^{2}(\mathbb{R}^{n})}+\left\|\Lambda^{\frac{n}{4}}\Delta \left(\mathbf{u}-\mathbf{q}\right)
  \right\|^{2} _{L^{2}(\mathbb{R}^{n})}\right)\\
  \vspace{5pt}
&\lesssim\displaystyle\left\| \mathbf{q}
  \right\| _{L^{2}(\mathbb{R}^{n})}\left(\left\| \nabla \left(\mathbf{u}-\mathbf{q}\right)
  \right\|^{2} _{L^{2}(\mathbb{R}^{n})}+\alpha^{2}\left\|\Delta \left(\mathbf{u}-\mathbf{q}\right)
  \right\|^{2} _{L^{2}(\mathbb{R}^{n})}\right)\\
  \vspace{5pt}
&\quad\displaystyle+\left\| \mathbf{q}
  \right\| _{L^{2}(\mathbb{R}^{n})}\left(\left\| \Lambda^{\frac{n}{4}}\nabla\left(\mathbf{u}-\mathbf{q}\right)
  \right\|^{2} _{L^{2}(\mathbb{R}^{n})}+\left\|\Lambda^{\frac{n}{4}}\Delta \left(\mathbf{u}-\mathbf{q}\right)
  \right\|^{2} _{L^{2}(\mathbb{R}^{n})}\right).
\end{array}
\right.
\end{equation}
\begin{equation}\label{5.38}
\left.
\begin{array}{ll}
\vspace{5pt}
A_{3}&=\left|\displaystyle\int_{\mathbb{R}^{n}} \left(\mathbf{u}-\alpha^{2}\Delta \mathbf{u}\right)\nabla \left(\mathbf{u}^{T}-\mathbf{q}^{T}\right)\Delta\left(\mathbf{u}-\mathbf{q}\right)
  dx\right|  \\
  \vspace{5pt}
&\lesssim\displaystyle\left(\left\| \mathbf{u}
  \right\| _{L^{2}(\mathbb{R}^{n})}+\left\| \Delta\mathbf{u}
  \right\| _{L^{2}(\mathbb{R}^{n})}\right) \left\| \nabla \left(\mathbf{u}-\mathbf{q}\right) \Delta
  \left(\mathbf{u}-\mathbf{q}\right)
  \right\| _{L^{2}(\mathbb{R}^{n})} \\
  \vspace{5pt}
&\lesssim\displaystyle\left(\left\| \mathbf{u}
  \right\| _{L^{2}(\mathbb{R}^{n})}+\left\| \Delta\mathbf{u}
  \right\| _{L^{2}(\mathbb{R}^{n})}\right)\left\| \nabla \left(\mathbf{u}-\mathbf{q}\right)
  \right\| _{L^{4}(\mathbb{R}^{n})}\left\| \Delta \left(\mathbf{u}-\mathbf{q}\right)
  \right\| _{L^{4}(\mathbb{R}^{n})} \\
  \vspace{5pt}
&\lesssim\displaystyle\left(\left\| \mathbf{u}
  \right\| _{L^{2}(\mathbb{R}^{n})}+\left\| \Delta\mathbf{u}
  \right\| _{L^{2}(\mathbb{R}^{n})}\right)\left\| \nabla \left(\mathbf{u}-\mathbf{q}\right)
  \right\|^{2} _{L^{4}(\mathbb{R}^{n})}\\
  \vspace{5pt}
&\quad\displaystyle
  +\left(\left\| \mathbf{u}
  \right\| _{L^{2}(\mathbb{R}^{n})}+\left\| \Delta\mathbf{u}
  \right\| _{L^{2}(\mathbb{R}^{n})}\right)\left\| \Delta \left(\mathbf{u}-\mathbf{q}\right)
  \right\|^{2} _{L^{4}(\mathbb{R}^{n})}\\
  \vspace{5pt}
&\lesssim\displaystyle\left(\left\| \mathbf{u}
  \right\| _{L^{2}(\mathbb{R}^{n})}+\left\| \Delta\mathbf{u}
  \right\| _{L^{2}(\mathbb{R}^{n})}\right)\left(\left\| \nabla \left(\mathbf{u}-\mathbf{q}\right)
  \right\|^{2} _{L^{2}(\mathbb{R}^{n})}+\left\|\Lambda^{\frac{n}{4}} \nabla \left(\mathbf{u}-\mathbf{q}\right)
  \right\|^{2} _{L^{2}(\mathbb{R}^{n})}\right)\\
  \vspace{5pt}
&\quad\displaystyle+\left(\left\| \mathbf{u}
  \right\| _{L^{2}(\mathbb{R}^{n})}+\left\| \Delta\mathbf{u}
  \right\| _{L^{2}(\mathbb{R}^{n})}\right)\left(\left\| \Delta\left(\mathbf{u}-\mathbf{q}\right)
  \right\|^{2} _{L^{2}(\mathbb{R}^{n})}+\left\|\Lambda^{\frac{n}{4}}\Delta \left(\mathbf{u}-\mathbf{q}\right)
  \right\|^{2} _{L^{2}(\mathbb{R}^{n})}\right)\\
  \vspace{5pt}
&\lesssim\displaystyle\left(\left\| \mathbf{u}
  \right\| _{L^{2}(\mathbb{R}^{n})}+\left\| \Delta\mathbf{u}
  \right\| _{L^{2}(\mathbb{R}^{n})}\right)\left(\left\| \nabla \left(\mathbf{u}-\mathbf{q}\right)
  \right\|^{2} _{L^{2}(\mathbb{R}^{n})}+\alpha^{2}\left\|\Delta \left(\mathbf{u}-\mathbf{q}\right)
  \right\|^{2} _{L^{2}(\mathbb{R}^{n})}\right)\\
  \vspace{5pt}
&\quad\displaystyle+\left(\left\| \mathbf{u}
  \right\| _{L^{2}(\mathbb{R}^{n})}+\left\| \Delta\mathbf{u}
  \right\| _{L^{2}(\mathbb{R}^{n})}\right)\left(\left\| \Lambda^{\frac{n}{4}}\nabla\left(\mathbf{u}-\mathbf{q}\right)
  \right\|^{2} _{L^{2}(\mathbb{R}^{n})}+\left\|\Lambda^{\frac{n}{4}}\Delta \left(\mathbf{u}-\mathbf{q}\right)
  \right\|^{2} _{L^{2}(\mathbb{R}^{n})}\right).
\end{array}
\right.
\end{equation}
Making the similar argument employed in estimating $A_{1}$, $A_{2}$ and $A_{3}$, we obtain
 \begin{equation}\label{5.39}
\left.
\begin{array}{ll}
\vspace{5pt}
 A_{4}&=\displaystyle\left|\int_{\mathbb{R}^{n}}\left[\left(\mathbf{u}-\mathbf{q}\right)
  -\alpha^{2}\Delta \left(\mathbf{u}-\mathbf{q}\right) \right]\nabla \mathbf{q}^{T}\Delta\left(\mathbf{u}-\mathbf{q}\right)dx\right|\\
  \vspace{5pt}
  &\lesssim\displaystyle \left(\left\|\nabla\left(\mathbf{u}-\mathbf{q}\right)
  \right\|^{2} _{L^{2}(\mathbb{R}^{n})}+\alpha^{2}\left\|\Delta\left(\mathbf{u}-\mathbf{q}\right)
  \right\|^{2} _{L^{2}(\mathbb{R}^{n})}\right)\\
  \vspace{5pt}
 &\quad\displaystyle+ \left(\left\|\Lambda^{ \frac{n}{4}} \mathbf{q}
  \right\|^{2} _{L^{2}(\mathbb{R}^{n})}+\left\| \nabla\mathbf{q}
  \right\|^{2} _{L^{2}(\mathbb{R}^{n})}\right)\\
  \vspace{5pt}
&\qquad\displaystyle\times\left(\left\|\Lambda^{ \frac{n}{4}}\nabla\left(\mathbf{u}-\mathbf{q}\right)
  \right\|^{2} _{L^{2}(\mathbb{R}^{n})}+\left\|\Lambda^{ \frac{n}{4}}\Delta\left(\mathbf{u}-\mathbf{q}\right)
  \right\|^{2} _{L^{2}(\mathbb{R}^{n})}\right).
\end{array}
\right.
\end{equation}
In addition, by interpolation inequality, for $n=2,3$, there holds
\begin{equation*}
\left.
\begin{array}{ll}
    \left\|\Lambda^{\frac{n}{4}} \mathbf{u}
  \right\|^{2} _{L^{2}(\mathbb{R}^{n})}\lesssim  \left\|  \mathbf{u}
  \right\|  _{L^{2}(\mathbb{R}^{n})} \left\|\nabla \mathbf{u}
  \right\| _{L^{2}(\mathbb{R}^{n})}\lesssim  \left\|  \mathbf{u}
  \right\|^{2}  _{L^{2}(\mathbb{R}^{n})}+ \left\|\nabla \mathbf{u}
  \right\|^{2} _{L^{2}(\mathbb{R}^{n})},
  \\[0.3cm]
   \left\|\Lambda^{\frac{n}{4}} \mathbf{q}
  \right\|^{2} _{L^{2}(\mathbb{R}^{n})}\lesssim  \left\|  \mathbf{q}
  \right\|  _{L^{2}(\mathbb{R}^{n})} \left\|\nabla \mathbf{q}
  \right\| _{L^{2}(\mathbb{R}^{n})}\lesssim  \left\|  \mathbf{q}
  \right\|^{2}  _{L^{2}(\mathbb{R}^{n})}+ \left\|\nabla \mathbf{q}
  \right\|^{2} _{L^{2}(\mathbb{R}^{n})},
\end{array}
\right.
\end{equation*}
recall  \eqref{5.16},  \eqref{5.36}, \eqref{5.38} and \eqref{5.39}, note that by \eqref{1.11} and  \eqref{2.3}
\begin{equation*}
    \left\| \mathbf{q}
  \right\|^{2} _{L^{2}(\mathbb{R}^{n})}+\left\| \mathbf{u}
  \right\|^{2} _{L^{2}(\mathbb{R}^{n})}+\left\| \nabla\mathbf{q}
  \right\|^{2} _{L^{2}(\mathbb{R}^{n})}+\left\| \nabla\mathbf{u}
  \right\|^{2} _{L^{2}(\mathbb{R}^{n})}+\left\| \Delta\mathbf{u}
  \right\|^{2} _{L^{2}(\mathbb{R}^{n})}\lesssim\displaystyle  \left\|  \mathbf{v}
  \right\|^{2} _{L^{2}(\mathbb{R}^{n})}\lesssim\displaystyle \left\|  \mathbf{v}_{0}
  \right\|^{2} _{L^{2}(\mathbb{R}^{n})},
\end{equation*}
 and the assumptions of Theorem \ref{t1.1} that for $\displaystyle s=\frac{n}{4}$, $\displaystyle \left\|  \mathbf{v}_{0}
  \right\|^{2} _{L^{2}(\mathbb{R}^{n})}\lesssim \varepsilon^{*}$ for $\varepsilon^{*}$ sufficiently small, in particular, choosing $\displaystyle\varepsilon^{*}\lesssim \frac{3\nu}{8}$ and using interpolation inequality, we obtain
\begin{equation}\label{5.40}
\left.
\begin{array}{ll}
\vspace{5pt}
 & \displaystyle\frac{1}{2}\frac{d}{dt}\left(\left\|\nabla\left(\mathbf{u}-\mathbf{q}\right)
  \right\|^{2} _{L^{2}(\mathbb{R}^{n})}+\alpha^{2}\left\|\Delta\left(\mathbf{u}-\mathbf{q}\right)
  \right\|^{2} _{L^{2}(\mathbb{R}^{n})}\right)\\
  \vspace{5pt}
  &\qquad\qquad\displaystyle+\frac{\nu}{2}\left(\left\|\Lambda^{ \frac{n}{4}}\nabla\left(\mathbf{u}-\mathbf{q}\right)
  \right\|^{2} _{L^{2}(\mathbb{R}^{n})}+\left\|\Lambda^{ \frac{n}{4}}\Delta\left(\mathbf{u}-\mathbf{q}\right)
  \right\|^{2} _{L^{2}(\mathbb{R}^{n})}\right)\\
  \vspace{5pt}
  &\qquad\lesssim\displaystyle\left(1+\left\|\Lambda^{ \frac{n}{4}} \Delta\mathbf{u}
  \right\|^{2} _{L^{2}(\mathbb{R}^{n})}\right)\left(\left\|\nabla\left(\mathbf{u}-\mathbf{q}\right)
  \right\|^{2} _{L^{2}(\mathbb{R}^{n})}+\alpha^{2}\left\|\Delta\left(\mathbf{u}-\mathbf{q}\right)
  \right\|^{2} _{L^{2}(\mathbb{R}^{n})}\right) \\
  \vspace{5pt}
  &\qquad\lesssim\displaystyle\left(1+\left\|\Lambda^{ \frac{n}{4}} \mathbf{v}
  \right\|^{2} _{L^{2}(\mathbb{R}^{n})}\right)\left(\left\|\nabla\left(\mathbf{u}-\mathbf{q}\right)
  \right\|^{2} _{L^{2}(\mathbb{R}^{n})}+\alpha^{2}\left\|\Delta\left(\mathbf{u}-\mathbf{q}\right)
  \right\|^{2} _{L^{2}(\mathbb{R}^{n})}\right).
\end{array}
\right.
\end{equation}
Notice that the initial data of both solutions $\mathbf{v}$ and $\mathbf{w}$ coincide, it is easy to check that  $\displaystyle\left\|\nabla\left(\mathbf{u}_{0}-\mathbf{q}_{0}\right)
  \right\|^{2} _{L^{2}(\mathbb{R}^{n})}=0$ and $\displaystyle\left\|\Delta\left(\mathbf{u}_{0}-\mathbf{q}_{0}\right)
  \right\|^{2} _{L^{2}(\mathbb{R}^{n})}=0$ from  \eqref{5.3}. By \eqref{1.11} and Gronwall's inequality, combining  \eqref{5.28} with  \eqref{5.40} yields that for $\displaystyle \frac{n}{4}\leq s<1$ and for any $t\in [0,T]$,
\begin{equation}\label{5.41}
\left.
\begin{array}{ll}
 \displaystyle\left\|\nabla\left(\mathbf{u}-\mathbf{q}\right)
  \right\|^{2} _{L^{2}(\mathbb{R}^{n})}+\alpha^{2}\left\|\Delta\left(\mathbf{u}-\mathbf{q}\right)
  \right\|^{2} _{L^{2}(\mathbb{R}^{n})}=0.
\end{array}
\right.
\end{equation}
Note that \eqref{5.3}, it follows from \eqref{5.15} and \eqref{5.41} that for $\displaystyle \frac{n}{4}\leq s<1$ and for any $t\in [0,T]$, $\displaystyle\left\| \mathbf{v}-\mathbf{w}
  \right\|^{2} _{L^{2}(\mathbb{R}^{n})}=0$.
  \\[0.3cm]
  \indent This finishes the proof of Theorem \ref{t1.3}.\hfill$\Box$\\

At the end of this section, we present the proof of Corollary \ref{c1.4}.
 \\[0.3cm]
{\bf Proof.}\quad Differentiating the second equation in \eqref{1.1} with respect to $x$ and $t$ shows
\begin{equation*}
\partial^{k}_{t} \nabla^{m}\mathbf{u}-\alpha^{2}\partial^{k}_{t} \nabla^{m}\Delta \mathbf{u}=\partial^{k}_{t} \nabla^{m}\mathbf{v}.
\end{equation*}
Squaring this equation and integrating in space yields, after some integration by parts,
\begin{equation*}
\left\|\partial^{k}_{t} \nabla^{m}\mathbf{u}\right\|^{2}_{L^{2}(\mathbb{R}^{n})}+2\alpha^{2}\left\|\partial^{k}_{t} \nabla^{m+1} \mathbf{u}\right\|^{2}_{L^{2}(\mathbb{R}^{n})}+\alpha^{4}\left\|\partial^{k}_{t} \nabla^{m+2} \mathbf{u}\right\|^{2}_{L^{2}(\mathbb{R}^{n})}= \left\|\partial^{k}_{t} \nabla^{m} \mathbf{v}\right\|^{2}_{L^{2}(\mathbb{R}^{n})}.
\end{equation*}
This is the identity \eqref{1.14}.\\
\indent Note that $\frac{n}{4}\leq s<1$ for $n=2,3$, thanks to the Gagliardo-Nirenberg-Sobolev inequality and interpolation inequality, it follows from \eqref{1.14} that
 \begin{equation}\label{5.42}
\left\{
\begin{array}{ll}
 \left\| \partial^{k}_{t} \nabla^{m}\mathbf{u}\right\|^{2}_{L^{n}(\mathbb{R}^{n})}&\leq C\left\| \partial^{k}_{t} \nabla^{m}\mathbf{v}\right\|^{2}_{L^{2}(\mathbb{R}^{n})},
 \\[0.3cm]
\left\| \partial^{k}_{t} \nabla^{m+1}\mathbf{u}\right\|^{2}_{L^{n}(\mathbb{R}^{n})}&\leq C\left\| \partial^{k}_{t} \nabla^{m}\mathbf{v}\right\|^{2}_{L^{2}(\mathbb{R}^{n})},
 \\[0.3cm]
\left\| \partial^{k}_{t} \nabla^{m}\Lambda^{s}\mathbf{u}\right\|^{2}_{L^{n}(\mathbb{R}^{n})}& \leq C\left\| \partial^{k}_{t} \nabla^{m}\mathbf{v}\right\|^{2} _{L^{2}(\mathbb{R}^{n})}.
   \end{array}
\right.
\end{equation}
 This yields the estimate \eqref{1.15}.\\
 \indent Finally, combining the first inequality in \eqref{5.42} with the regularity bounds \eqref{1.13} yields the estimate \eqref{1.16}.\\
 \indent This finishes the proof of Corollary \ref{1.4}.\hfill$\Box$
\section*{Acknowledgments }
Zaihui Gan is partially supported by the National Science Foundation
of China under grants (No. 11571254).

\end{document}